\documentclass[12pt]{article}

\usepackage{amssymb,amsmath,amsthm,amsfonts,mathrsfs}
\usepackage{subfigure,epsfig,psfrag,epstopdf}
\usepackage{fullpage,graphicx,color}
\usepackage{algorithm,algorithmic}
\usepackage{tabularx,array,booktabs,bm,multirow}
\usepackage{longtable}
\usepackage{cite,url}
\usepackage[bookmarks=false]{hyperref}
\usepackage[labelfont=bf]{caption}
\usepackage{bbm}
\usepackage{tikz}
\usetikzlibrary{positioning, shapes, arrows.meta}
\theoremstyle{definition}
\newtheorem{definition}{Definition}
\theoremstyle{theorem}
\newtheorem{theorem}{Theorem}
\theoremstyle{proposition}
\newtheorem{proposition}{Proposition}
\theoremstyle{corollary}
\newtheorem{corollary}{Corollary}

\newcommand{\mB}{\mathcal{B}}
\newcommand{\mD}{\mathcal{D}}
\newcommand{\mF}{\mathcal{F}}
\newcommand{\mG}{\mathcal{G}}
\newcommand{\mI}{\mathcal{I}}
\newcommand{\mJ}{\mathcal{J}}
\newcommand{\mK}{\mathcal{K}}
\newcommand{\mL}{\mathcal{L}}
\newcommand{\mM}{\mathcal{M}}
\newcommand{\mO}{\mathcal{O}}
\newcommand{\mP}{\mathcal{P}}
\newcommand{\mR}{\mathcal{R}}

\newcommand{\mX}{\mathcal{X}}
\newcommand{\mY}{\mathcal{Y}}
\newcommand{\RR}{\mathbb{R}}

\newcommand{\ba}{{\bf a}}
\newcommand{\bB}{{\bf B}}
\newcommand{\bC}{{\bf C}}

\newcommand{\bF}{{\bf F}}
\newcommand{\bL}{{\bf L}}
\newcommand{\bn}{{\bf n}}
\newcommand{\bS}{{\bf S}}
\newcommand{\bu}{{\bf u}}
\newcommand{\bv}{{\boldsymbol{\upsilon}}}
\newcommand{\bx}{{\bf x}}
\newcommand{\by}{{\bf y}}
\newcommand{\bz}{{\bf z}}
\def\bnabla{\boldsymbol{\nabla}}
\def\blambda{\boldsymbol{\lambda}}

\newcommand{\argmin}{\operatorname{argmin}}

\begin{document}

\title{Variational State-Dependent Inverse Problems in PDE-Constrained Optimization:
A Survey of Contemporary Computational Methods and Applications}
\date{}
\author{{\bf \normalsize Vladislav Bukshtynov}\footnote{Corresponding author: \url{Vladislav.Bukshtynov@colorado.edu}} \\
{\it \small Department of Mathematics, University of Colorado Boulder, Boulder, CO 80309, USA}}
\maketitle

\begin{abstract}
State-dependent parameter identification, where unknown model parameters depend on one or more state variables in partial
differential equations~(PDEs) or coupled PDE systems, is fundamental to a wide range of problems in physics, engineering,
and materials science. This review surveys PDE-constrained optimization approaches for such inverse problems, emphasizing
the underlying mathematical theory and key computational advances developed since 2011. We discuss variational formulations,
adjoint-based gradient methods, regularization strategies, and modern computational frameworks, and highlight representative
applications, with particular emphasis on identifiability, ill-posedness, and structural limits of state-dependent inverse
problems. The review concludes with major open challenges and emerging research directions related to nonconvexity, identifiability,
regularization, adjoint computation, data limitations, and model-class dependence.

\end{abstract}

\section{Introduction}
\label{sec:intro}

Reliable mathematical and computational modeling of complex physical phenomena critically depends on accurate knowledge of
the underlying material properties~\cite{Bukshtynov2011,Bukshtynov2013}. In a broad range of applications in physics,
engineering, and materials science, these properties (or constitutive relations) are not fixed constants but vary with
the \emph{state(s)} of the modeled system~\cite{ChaventLemonnier1974,Bukshtynov2011}. Examples include temperature-dependent
thermal conductivity and viscosity in heat transfer and fluid mechanics; strain- or deformation-dependent constitutive
parameters in solid mechanics; concentration-dependent diffusivities in porous media; and more generally, transport or
reaction coefficients that evolve with phase, saturation, or chemical composition~\cite{Bukshtynov2011,Bukshtynov2013}.
Accurately identifying such \emph{state-dependent} material properties is essential for predictive simulations, quantitative
design, and control of multiphysics systems. Yet this task is challenging, as the functional dependence of these coefficients
on state variables is often unknown, highly nonlinear, and difficult to measure directly, particularly for complex materials
such as liquid metals and metal alloys~\cite{VolkovProtasLiaoGlander2009}, newly created polymers and composites, or
high-temperature fluids~\cite{Bukshtynov2011,Bukshtynov2013,BukshtynovPhD2012}. Moreover, constitutive laws inferred from
inverse problems should often be interpreted as effective, model-consistent descriptions valid within a prescribed regime,
rather than intrinsic material properties independent of modeling assumptions.

Inverse problems of parameter estimation for partial differential equations (PDEs) offer a principled framework for inferring
unknown material properties from indirect, incomplete, and possibly noisy observations (measurements). Classical formulations
have focused primarily on \emph{space-dependent} parameters, i.e., material properties viewed as functions of the spatial
coordinates. Such problems are relatively well studied and form the basis of many imaging technologies in
medicine~\cite{NashedScherzer2002,AdlerHolder2022} and geosciences~\cite{GottliebDuChateau1996,OliverReynoldsLiu2008}.
They often lead to PDE-constrained optimization problems with well-developed mathematical theory~\cite{BanksKunisch1989,Isakov2017},
regularization strategies~\cite{EnglHankeNeubauer1996,Kugler2000}, and computational algorithms~\cite{Tarantola2005,NashedScherzer2002}.
In contrast, the identification of \emph{state-dependent} parameters (e.g., functions of temperature, pressure, concentration,
deformation, or other dependent variables) poses a fundamentally different type of inverse problem. Here, the goal is to recover
a functional relationship that holds \emph{uniformly} across the entire spatial domain, coupling the unknown parameter directly
to the solution of the governing PDEs. This leads to a strongly nonlinear inverse problem and introduces computational and
analytical challenges absent in the classical space-dependent setting. A defining feature of such problems is that identifiability
of the constitutive relation is typically restricted to the range of state values actually explored by the system, leading to
fundamental underdetermination outside this regime. As this review emphasizes, these challenges are not merely computational,
but reflect intrinsic structural limitations of the inverse problem, including restricted identifiability, non-injectivity of
the inverse map, and objective-dependent nonexistence of solutions. This review focuses exclusively on \emph{deterministic and
variational} formulations of state-dependent inverse problems within a PDE-constrained optimization framework. Bayesian and
probabilistic approaches to inverse problems, while important and well developed in their own right, address fundamentally
different questions related to uncertainty quantification and statistical inference and are therefore beyond the scope of
the present survey.

Despite its importance, state-dependent parameter identification has historically received limited attention. Foundational
analytical work was established in 1974 by Chavent and Lemonnier~\cite{ChaventLemonnier1974}, who introduced a variational
framework, proved existence results, and derived adjoint expressions for the gradients of least-squares cost functionals. Early
computational efforts, such as those then surveyed and extended in~\cite{Bukshtynov2011,Bukshtynov2013}, explored simplified
model problems~\cite{DuChateau2004,CannonDuChateau1980a,DuChateau1997}, problems posed in an infinite domain~\cite{Luo2003},
discrete formulations~\cite{HankeScherzer1999}, and spline-based parameterizations~\cite{Alifanov2004,Alifanov2007}, often
employing regularization~\cite{Kugler2000,Kugler2003,Neubauer2008}, linearization~\cite{Tai1995}, or Green's function
techniques~\cite{JanickiKindermann2009} to address the ill-posed nature of the problem. These pioneering contributions laid
the groundwork for modern approaches, but they preceded the development of mature adjoint-based PDE-constrained optimization
methods that now underpin contemporary large-scale computations.

Over the past decade, significant progress has been made in applying \emph{full-length} PDE-constrained optimization techniques
to state-dependent parameter identification. Advances in continuous, so-called ``optimize--then--discretize'', formulations,
automatic differentiation, adjoint solvers, multilevel discretizations, regularization strategies, and numerical methods for
computing gradients defined on level sets have enabled robust and efficient algorithms applicable to multiphysics systems.
These developments have broadened the scope of problems that can be addressed, including time-dependent, nonlinear, and coupled
PDE models arising in fluid mechanics, heat transfer, solid-fluid interaction, and other areas.

\medskip
\noindent\textbf{Structure and Scope.}
The length and level of detail of this review reflect the historically fragmented development of inverse problems with state-dependent
parameters across analysis, optimization, and applications. Rather than providing a tutorial introduction, the goal is to synthesize
foundational analytical results, modern computational formulations, and representative applications into a unified reference that makes
explicit the structural limitations, identifiability mechanisms, and modeling assumptions common to this class of problems. This paper
provides a comprehensive examination of the theory and computation underlying state-dependent parameter identification using
PDE-constrained optimization. Specifically, we:

\begin{enumerate}
  \item Summarize the mathematical foundations established in the seminal work of Chavent and Lemonnier~\cite{ChaventLemonnier1974}
    (to the best of our knowledge, never appeared in the English language), including variational formulations, adjoint-based
    gradient derivations, and existence results.
  \item Review early computational approaches developed prior to 2011, including those in~\cite{Bukshtynov2011,Bukshtynov2013} and
    related works, that introduced optimization-based frameworks, functional parameterizations, and discrete formulations for
    simplified model problems.
  \item Focus on developments since 2011, highlighting continuous PDE-constrained optimization formulations, adjoint-based
    gradient evaluation techniques, regularization strategies, computational frameworks, and algorithms capable of treating
    realistic multiphysics systems with state-dependent coefficients.
\end{enumerate}

\noindent\textbf{Objectives.}
The goals of this review are threefold:
\begin{itemize}
  \item to synthesize the mathematical foundations of state-dependent parameter identification;
  \item to survey the major computational advances and algorithmic strategies developed over the past decade; and
  \item to identify open challenges, unresolved analytical questions, and emerging research directions that define the frontier
    of this growing field.
\end{itemize}

By synthesizing recent advances in theory, computation, and applications, this review aims to provide a unified perspective on
state-dependent parameter identification and to serve both as a comprehensive reference and an invitation to further research
in PDE-constrained optimization for complex material systems.

\section{Mathematical Foundations}
\label{sec:math}

\subsection{Problem Formulation}
\label{sec:problem_formulation}

This section develops a general mathematical framework for inverse problems involving \emph{state-dependent} parameters (i.e.,
material properties) governed by PDEs. The proposed formulation builds upon and generalizes the approaches
of~\cite{ChaventLemonnier1974,Bukshtynov2011,Bukshtynov2013,BukshtynovPhD2012}, providing a unified framework capable of addressing
a broad spectrum of problems arising in multiphysics modeling, including continuum mechanics and nonequilibrium thermodynamics.

\subsubsection{Governing equations with state-dependent parameters}
\label{sec:gov_eqn}

Let $\Omega \subset \RR^d$, $d=1,2,3$, be a bounded open domain with sufficiently smooth boundary $\partial\Omega$, and let $[0,t_f]$
denote a time interval (in the stationary case, time dependence is omitted). We consider a state variable
\begin{equation}
  \bu : \Omega \times (0,t_f] \rightarrow \RR^{d+1},
  \label{eq:state}
\end{equation}
governed by a system of PDEs of the abstract form
\begin{equation}
  \mG \bigl( \bu, p(\bu) \bigr) = f \quad \text{in } \Omega \times (0,t_f],
  \label{eq:general_PDE}
\end{equation}
supplemented with boundary and initial conditions
\begin{equation}
  \begin{aligned}
    \mB \bigl( \bu, p(\bu) \bigr) &= g \quad \text{on } \partial\Omega \times (0,t_f],\\
    \bu(\cdot,0) &= \bu_0 \quad \text{in } \Omega.
  \end{aligned}
  \label{eq:general_BC}
\end{equation}
Here, $\mG$ denotes a (possibly nonlinear) differential operator, $\mB$ represents boundary operators (Dirichlet, Neumann,
Robin, or mixed), $f$ and $g$ are known source terms and boundary data, and $p(\bu)$ is an \emph{unknown constitutive relation}
which depends on one or more components of the state variable $\bu$.

The defining feature of the problems considered here is that the unknown parameter $p$ depends on the \emph{dependent variable}
$\bu$, rather than on spatial or temporal coordinates. This distinction separates state-dependent inverse problems from classical
space-dependent parameter identification and gives rise to fundamentally different analytical and computational challenges.
The abstract formulation \eqref{eq:general_PDE}--\eqref{eq:general_BC} encompasses a wide range of models. Two canonical examples,
studied extensively in~\cite{Bukshtynov2011,Bukshtynov2013}, are briefly recalled for concreteness.

\paragraph{Heat conduction with temperature-dependent conductivity~\cite{Bukshtynov2011}.}
Let $\bu = T(\bx)$, $\bx \in \Omega$, denote the temperature field. Assuming a nonlinear, temperature-dependent thermal conductivity
$k(T)$, the governing equation obtained from energy conservation takes the form
\begin{subequations}
  \label{eq:heat}
  \begin{align}
    -\bnabla \cdot \bigl(k(T)\bnabla T\bigr) &= f \quad \text{in } \Omega, \\
    T &= T_0 \quad \text{on } \partial\Omega.
  \end{align}
\end{subequations}

\paragraph{Incompressible flow with temperature-dependent viscosity~\cite{Bukshtynov2013}.}
Let $\bu=(\bv, p, T)$ denote velocity, pressure, and temperature. The viscosity $\mu(T)$ enters the incompressible
Navier--Stokes equations\footnote{$(\cdot)^{\top}$ denotes the adjoint (transpose) of a linear mapping.} as
\begin{subequations}
  \label{eq:NS_heat}
  \begin{align}
    \partial_t \bv + (\bv \cdot \bnabla) \bv + \bnabla p - \bnabla \cdot
    \left[ \mu(T) \bigl( \bnabla \bv + (\bnabla \bv)^{\top} \bigr) \right] &= 0 \quad \text{in } \Omega, \label{eq:NS_heat_a} \\
    \bnabla \cdot \bv &= 0 \quad \text{in } \Omega, \label{eq:NS_heat_b} \\
    \partial_t T + (\bv \cdot \bnabla)T - \bnabla \cdot (k \bnabla T) &= 0 \quad \text{in } \Omega, \label{eq:NS_heat_c}
  \end{align}
\end{subequations}
subject to appropriate boundary and initial conditions.

\subsubsection{Admissible set of constitutive relations}
\label{sec:adm_set}

Let $\zeta$ denote a scalar state variable (e.g., temperature, concentration, strain invariant) extracted from $\bu$.
The constitutive relation
\begin{equation}
  \label{eq:scalar_state}
  p = p(\zeta)
\end{equation}
is assumed to be defined on a state (e.g., desired temperature) interval $\mD \subset \RR$.

\begin{definition}[Identifiability interval]
\label{def:intI}
Let $\bu$ be a solution of \eqref{eq:general_PDE}--\eqref{eq:general_BC}. Following~\cite{Kugler2003},
the \emph{identifiability interval} is defined as
\begin{equation}
  \label{eq:intI}
  \mI := [\zeta_{\alpha}, \zeta_{\beta}], \qquad
  \zeta_{\alpha} = \min_{(\bx,t) \in \overline{\Omega} \times [0, t_f]} \zeta(\bx,t), \qquad
  \zeta_{\beta} = \max_{(\bx,t) \in \overline{\Omega} \times [0,t_f]} \zeta (\bx,t).
\end{equation}
Only values of $p(\zeta)$ for $\zeta \in \mI$ can be informed by the data~\cite{ChaventLemonnier1974,Bukshtynov2011,Bukshtynov2013}.
\end{definition}

\begin{definition}[Reconstruction interval]
Let
\begin{equation}
  \label{eq:intD}
  \mD := [\zeta_a, \zeta_b], \qquad \zeta_a \le \inf \mI, \qquad \zeta_b \ge \sup \mI.
\end{equation}
The interval $\mD$ is the \emph{interval of (desired) reconstruction} on which the constitutive relation $p(\zeta)$ is sought.
In general, $\mI \subseteq \mD$.
\end{definition}

\begin{definition}[Measurement span]
Let $\tilde{\zeta}$ denote $M$ measured values of state variable $\zeta$. The \emph{measurement span} is defined as
\begin{equation}
  \label{eq:intM}
  \mM := \bigl[ \, \min_{1 \le i \le M} \min_{0 < t \le t_f} \tilde{\zeta}_i(t), \;
  \max_{1 \le i \le M} \max_{0 < t \le t_f} \tilde{\zeta}_i(t) \, \bigr], \qquad \mM \subseteq \mI.
\end{equation}
\end{definition}

Following~\cite{Bukshtynov2011,Bukshtynov2013}, admissible constitutive relations $p(\zeta)$ are assumed to belong to
the following \emph{admissible set}
\begin{equation}
  \mP = \left\{ p(\zeta) \in C^1(\mD) \; \big| \; 0 < m_p \le p(\zeta) \le M_p < \infty, \ \forall \zeta \in \mD \right\},
  \qquad \mP \subset \mX,
\label{eq:Pset}
\end{equation}
where bounds $m_p$ and $M_p$ reflect admissibility for the underlying physical phenomenon (e.g., the Clausius--Duhem inequality
for thermodynamic admissibility~\cite{Muschik1989,Triano2008}) and ensure well-posedness of the governing PDEs~\cite{ColemanNoll1963,Liu1972,TruesdellNoll2004}.
We also assume that the admissible set $\mP$ is embedded into an infinite-dimensional Hilbert space $\mX$, chosen according
to the requirements of the particular application.

\begin{figure}[!h]
\centering
\begin{tikzpicture}[
  scale=0.8,
  >=Stealth,
  axis/.style={->, thick},
  interval/.style={|<->|, ultra thick},
  Bcurve/.style={thick, blue},
  Rcurve/.style={thick, red, dashed},
  RRcurve/.style={thick, red, dotted},
  point/.style={circle, blue, fill=red, inner sep=1.6pt},
  label/.style={font=\small},
  fillD/.style={fill=gray!20, opacity=0.5},
  fillI/.style={fill=blue!20, opacity=0.5},
  fillM/.style={fill=red!20, opacity=0.5},
  proj/.style={-, dashed, gray},
  projB/.style={-, dashed, blue},
  projR/.style={-, dashed, red}
]

% Horizontal shift for right plot
\def\xshift{7.0cm}

% Vertical bounds
\def\ymin{-4.0}
\def\ymax{5.0}

% Tick height
\def\tickh{0.1} % increased tick height

% =================================================
% ALL SPANS
% =================================================
\fill[fillD] (-5,-3.6) rectangle (14, 4.2);
\fill[fillI] (-4,-2.6) rectangle (13, 3.1);
\fill[fillM] (-3,-2.0) rectangle (12, 2.0);

\draw[interval] (14, -3.6) -- (14, 4.2)
  node[midway, rotate=90, yshift=-10pt, label] {reconstruction interval $\mD$};

\draw[interval, blue] (13, -2.6) -- (13, 3.1)
  node[midway, rotate=90, yshift=-11pt, label] {identifiability interval $\mI$};

\draw[interval, red] (12, -2.0) -- (12, 2.0)
  node[midway, rotate=90, yshift=-11pt, label] {measurement span $\mM$};

% =================================================
% LEFT PLANE
% =================================================

% Axes
\draw[axis] (-5,0) -- (5,0) node[below] {$x$};
\draw[axis] (0,\ymin) -- (0,\ymax) node[right] {$\zeta$};

% Blue curve
\draw[Bcurve]
  plot[smooth] coordinates {
    (-4, 3.1) (-3, 2.0) (-2, 1.3) (-1, 0.9)
    (-0.5, 0.7) (0.5, 0.3) (2, -0.5) (3, -2.0) (4, -2.6)
  };

% Vertical dashed lines
\draw[proj] (-4, 3.1) -- (-4, 0);
\draw[proj] (-3, 2.0) -- (-3, 0);
\draw[proj] (-2, 1.3) -- (-2, 0);
\draw[proj] (-1, 0.9) -- (-1, 0);
\draw[proj] (3, -2.0) -- (3, 0);
\draw[proj] (4, -2.6) -- (4, 0);

% X-axis ticks
\foreach \x/\xlabel/\above in {-4/$-1$/below,-3/$x_1$/below,-2/$x_2$/below,3/$x_M$/above,4/$1$/above} {
  \draw (\x,\tickh) -- (\x,-\tickh);
  \ifthenelse{\equal{\above}{above}}
    {\node[label, above] at (\x,0) {\xlabel};}
    {\node[label, below] at (\x,0) {\xlabel};}
}

% Horizontal dashed lines
\draw[proj] (0, 4.2) -- (14, 4.2);
\draw[proj] (0, -3.6) -- (14, -3.6);
\draw[projB] (-4, 3.1) -- (13, 3.1);
\draw[projR] (-3, 2.0) -- (12, 2.0);
\draw[projR] (3, -2.0) -- (12, -2.0);
\draw[projB] (0, -2.6) -- (13, -2.6);

% Endpoint Y-labels
\node[label, left] at (0, 4.2) {$\zeta_b$};
\node[label, left, yshift=-6pt] at (0, 3.1) {$\zeta_{\beta}$};
\node[label, left] at (0, -2.6) {$\zeta_{\alpha}$};
\node[label, left] at (0, -3.6) {$\zeta_a$};

% Points on curve
\node[point] at (-3, 2.0) {};
\node[point] at (-2, 1.3) {};
\node[point] at (-1, 0.9) {};
\node[point] at (0.5, 0.3) {};
\node[point] at (3, -2.0) {};

% =================================================
% RIGHT PLANE
% =================================================
\begin{scope}[xshift=\xshift]

% Axes
\draw[axis] (-1,0) -- (4.5,0) node[below] {$p$};
\draw[axis] (0,\ymin) -- (0,\ymax) node[right] {$\zeta$};

% Red curve
\draw[Rcurve]
  plot[smooth] coordinates {
    (1.0, -2.6) (1.5, -0.7) (3.5, 1.2) (4.0, 3.1)
  };
\draw[RRcurve]
  plot[smooth] coordinates {
    (4.0, 3.1) (4.25, 3.8) (4.5, 4.2)
  };
\draw[RRcurve]
  plot[smooth] coordinates {
    (0.85, -3.6) (0.9, -3.3) (1.0, -2.6)
  };

\end{scope}

\end{tikzpicture}

\caption{Schematic showing (left) the solution $\zeta(x, t_0)$ at some fixed time $t_0$ and (right) the corresponding constitutive
  relation $p(\zeta)$ defined over their respective domains, i.e., $\Omega = (-1,1)$ and the identifiability region $\mI$. The dotted
  line represents an extension of the constitutive relation $p(\zeta)$ from $\mI$ to the (desired) reconstruction interval $\mD$.
  In the Figure on the right, the horizontal axis is to be interpreted as the ordinate. The concept is adopted from~\cite{Bukshtynov2013}.}

\label{fig:intervals}
\end{figure}
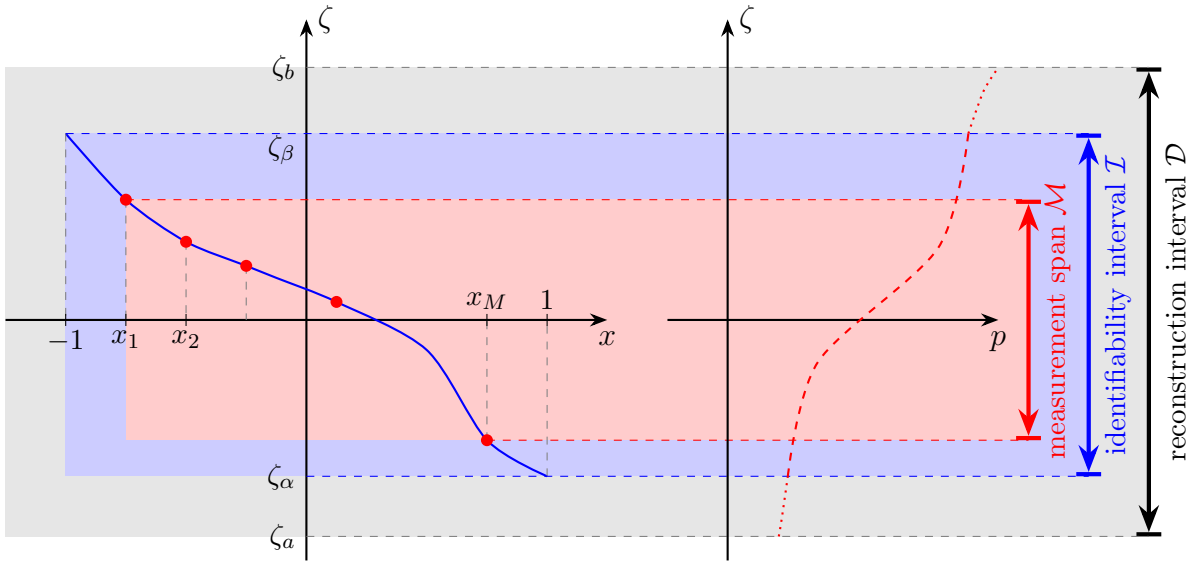

To provide a visual summary of these concepts, Figure~\ref{fig:intervals} illustrates the reconstruction interval $\mD$,
the identifiability interval $\mI$, and the measurement span $\mM$ along with a sample state-variable profile.

\subsubsection{PDE-constrained optimization formulation as inverse problem}
\label{sec:PDE_opt}

Let $\mO$ denote an observation operator mapping the state $\bu$ to measurable outputs (observations) $\by_{\text{obs}}$. Depending on
the application, $\mO$ may represent pointwise sensors placed inside domain $\Omega$, boundary measurements, time traces, or integral
quantities. Abstractly,
\begin{equation}
  \label{eq:observ}
  \mO : \bu \mapsto \by_{\text{obs}} \in \mY,
\end{equation}
where $\mY$ is a suitable Hilbert space.

The forward (parameter--to--observation) map
\begin{equation}
  \mF : p \mapsto \mO \bigl( \bu(p) \bigr)
  \label{eq:Fmap}
\end{equation}
defines an operator equation
\begin{equation}
  \mF(p) = \by_{\text{obs}},
  \label{eq:inverse_eq}
\end{equation}
which is typically nonlinear and ill-posed~\cite{Tarantola2005,NashedScherzer2002}.

Following~\cite{ChaventLemonnier1974,BanksKunisch1989,Bukshtynov2011,Bukshtynov2013}, the inverse problem \eqref{eq:inverse_eq}
is recast as a PDE-constrained optimization problem:
\begin{equation}
  \min_{p(\bu) \in \mP} \mJ(\bu, p(\bu))
  \label{eq:opt_general}
\end{equation}
subject to \eqref{eq:general_PDE}--\eqref{eq:general_BC}. A prototypical cost functional is
\begin{equation}
  \mJ(\bu, p(\bu)) = \frac{1}{2} \left\| \mO(\bu(p)) - \by_{\text{obs}} \right\|_{\mY}^2 + \mR(p),
  \label{eq:cost_general}
\end{equation}
where $\mR(p)$ is a regularization term, e.g.,
\begin{equation*}
  \mR(p) = \frac{\lambda}{2} \| p - \bar p \|_{H^1(\mD)}^2,
\end{equation*}
introduced to stabilize the reconstruction in the presence of noise in the observed data
$\by_{\text{obs}}$~\cite{TikhonovArsenin1977,Neubauer2008} by choosing a proper computational scheme with regularization weight
$\lambda$ and reference solution $\bar p$.

The optimization variable $p=p(\bu)$ is infinite-dimensional, and the corresponding gradient $\bnabla_p \mJ$ represents
the sensitivity of the cost functional $\mJ(\bu, p(\bu))$ with respect to perturbations of the constitutive relation $p(\bu)$ itself.
As emphasized in~\cite{Bukshtynov2011,Bukshtynov2013}, this structure differs fundamentally from standard PDE-constrained
optimization problems, where controls depend on independent variables~\cite{Biegler2007,Gunzburger2003}. To clarify the generality
of the framework introduced above, Table~\ref{tab:mapping} maps the abstract quantities appearing in \eqref{eq:state}--\eqref{eq:general_PDE}
and \eqref{eq:scalar_state}--\eqref{eq:cost_general} to representative model problems \eqref{eq:heat}--\eqref{eq:NS_heat} studied
in~\cite{Bukshtynov2011,Bukshtynov2013}.

\begin{table}[!htbp]
  \centering
  \caption{Mapping between abstract formulation and representative model problems.}
  \label{tab:mapping}
  \footnotesize
  \begin{tabular}{p{4.0cm} p{5.25cm} p{5.25cm}}
    \toprule
    \multicolumn{1}{c}{\textbf{Abstract object}} &
    \multicolumn{1}{c}{\textbf{Heat conduction~\cite{Bukshtynov2011}} by \eqref{eq:heat}} &
    \multicolumn{1}{c}{\textbf{Navier--Stokes flow~\cite{Bukshtynov2013}} by \eqref{eq:NS_heat}}\\
    \midrule
    Time dependency &
    Stationary problem &
    $t \in [0,t_f]$ \\[0.3em]

    State variable $\bu$ &
    $T(\bx)$ &
    $(\bv,p,T)(\mathbf{x},t)$ \\[0.3em]

    Scalar state variable $\zeta$ &
    $T$ &
    $T$ \\[0.3em]

    Constitutive relation $p(\zeta)$ &
    Thermal conductivity $k(T)$ &
    Dynamic viscosity $\mu(T)$ \\[0.3em]

    Governing operator $\mG$ &
    $-\bnabla \cdot (k(T) \bnabla T)$ &
    Navier--Stokes + energy equation \\[0.3em]

    Observation operator $\mO$ &
    $\tilde T(\bx), \ \bx \in \Sigma \subseteq \Omega$ &
    $\tilde T(\bx,t), \ \bx \in \Sigma \subseteq \Omega, \ t \in [0,t_f]$ \\[0.3em]

    Identifiability interval $\mI$ &
    $[\min T,\max T]$ in $\Omega$ &
    $[\min T,\max T]$ in $\Omega\times[0,t_f]$ \\[0.3em]

    Measurement span $\mM$ &
    Range of measured temperatures &
    Range of measured $T$-traces \\[0.3em]

    Regularization $\mR(p)$ &
    $\frac{\lambda}{2} \| k - \bar k \|_{L_2(\mI)}^2$ and $\frac{\lambda}{2} \| k - \bar k \|_{\dot{H}^1(\mI)}^2$&
    $\frac{\lambda}{2} \|\theta - \bar \theta\|_{\dot{H}^1(\mI)}^2, \ \mu(T) = \theta^2(T) + m_{\mu}$ \\
    \bottomrule
  \end{tabular}
\end{table}

We now state the inverse problem in a mathematically precise form. Given a spatial domain $\Omega \subset \RR^d$ and time interval
$[0,t_f]$, a PDE system of the form \eqref{eq:general_PDE}--\eqref{eq:general_BC}, noisy observations $y_{\mathrm{obs}} \in \mY$,
and an admissible set of constitutive relations $\mP$, find an \emph{optimal solution}
\begin{equation}
  \hat{p} = \underset{p \in \mP}{\argmin} \, \left\{ \frac{1}{2} \left\| \mO(\bu(p)) - y_{\mathrm{obs}} \right\|_{\mY}^2
  + \mR(p) \right\}.
  \label{eq:final_opt}
\end{equation}
Optimization problem \eqref{eq:final_opt} is characterized by the following features~\cite{ChaventLemonnier1974,Bukshtynov2011,Bukshtynov2013}:
\begin{enumerate}
  \item The control variable $p(\bu)$ depends on the solution $\bu$ itself.
  \item The forward operator $\mF(p)$ in \eqref{eq:inverse_eq} is nonlinear and compact.
  \item The problem \eqref{eq:final_opt} is generally nonconvex and locally ill-posed~\cite{Tarantola2005}.
  \item Admissibility constraints on $p(\bu)$ encode physics-related consistency.
\end{enumerate}

As shown in~\cite{Bukshtynov2011,Bukshtynov2013,Luenberger1976}, first-order optimality conditions can be derived by introducing
adjoint variables associated with \eqref{eq:general_PDE}, leading to gradient expressions defined on the identifiability interval
$\mI$ rather than on the physical domain $\Omega$. The derivation of these adjoint equations, the structure of the resulting
gradients, and their numerical realization form the subject of the following sections.

\subsection{Analytical Properties}
\label{sec:anal_properties}

\subsubsection{Analytical foundations by Chavent \& Lemonnier~(1974)}
\label{sec:ChaventLemonnier1974}

The analytical framework for inverse problems involving state-dependent parameters originates in the seminal work of Chavent and
Lemonnier~\cite{ChaventLemonnier1974}. Their analysis provides a rigorous foundation for the problem formulation adopted in
Section~\ref{sec:problem_formulation} and anticipates several concepts that have since become standard in PDE-constrained
optimization and nonlinear inverse problems. Since this work was never published in English, we provide here a comparatively
detailed account of its main analytical contributions, exceeding the level of exposition typically found in review articles.

Specifically, \cite{ChaventLemonnier1974} considers a quasilinear parabolic equation of the form
\begin{subequations}
  \label{eq:CL_forward}
  \begin{align}
    \partial_t \bu - \bnabla \cdot \bigl( p(\bu)\, \bnabla \bu \bigr) &= f \quad \text{in } Q := \Omega \times (0,T),
    \label{eq:CL_forward_u} \\
    \bnabla \bu \cdot \bn &= 0 \quad \text{on } \partial\Omega \times (0,T), \\
    \bu(\cdot,0) &= \bu_0 \quad \text{in } \Omega.
  \end{align}
\end{subequations}
The unknown is the \emph{constitutive law} $p : \RR \to \RR^+$, which depends on the state variable $\bu(\bx,t;p)$.
\footnote{The unknown $p$ is a scalar constitutive law $p(\zeta):\RR\to\RR^+$, which is evaluated along the state via the composition
$p(\bu(\bx,t;p))$; the resulting coefficient $p(\bu)$ is therefore defined on $Q$, although $p$ itself is not a space-time dependent
function as also stated in Section~\ref{sec:adm_set}.} The admissible set of parameters is defined through pointwise bounds
\begin{equation}
  0 < p_{\min} \le p(\zeta) \le p_{\max}, \qquad \zeta \in \RR,
  \label{eq:CL_bounds}
\end{equation}
which ensure uniform ellipticity of the differential operator. The inverse problem is formulated as the minimization of
a least-squares functional of output-tracking type,
\begin{equation}
  \mJ(p) = \sum_{i=1}^M \int_{0}^{T} \left( \frac{1}{|\Omega_i|} \int_{\Omega_i} \bu(\bx, t; p)\, d\bx - \bz_i(t) \right)^2 \, dt,
  \label{eq:Jp_CL}
\end{equation}
where $\bu(\bx,t;p)$ denotes the solution of the state equation~\eqref{eq:CL_forward_u} associated with the state-dependent
parameter $p(\cdot)$, $\{\Omega_i\}_{i=1}^M \subset \Omega$ are prescribed measurement subdomains with Lebesgue measure
$|\Omega_i|$, and $\bz_i(t)$ represents the time-dependent observation corresponding to the $i$-th measurement region.
The normalization factor $1/|\Omega_i|$ reflects the fact that the available data correspond to spatial averages of the
state over $\Omega_i$, rather than pointwise measurements.

\begin{theorem}[Well-posedness of the state equation {\cite{ChaventLemonnier1974,Lions1968,Lions1969,JosephLundgren1972}}]
  \label{thm:CL_wellposed}
  Let $\Omega \subset \mathbb{R}^d$ be a bounded domain with sufficiently smooth boundary, and let $T>0$. Assume that the
  constitutive law $p:\mathbb{R}\to\mathbb{R}^+$ satisfies the uniform ellipticity condition~\eqref{eq:CL_bounds}. Let the
  source term satisfy $f \in L^\infty(Q)$ and the initial condition satisfy $\bu(\cdot,0) \in L^\infty(\Omega)$. Impose
  homogeneous Neumann boundary conditions on $\partial\Omega$.
  
  Then the quasilinear parabolic problem~\eqref{eq:CL_forward_u} admits a unique weak solution
  \begin{equation*}
    \bu \in L^2(0,T; H^1(\Omega)) \cap L^\infty(Q), \qquad \partial_t \bu \in L^2(0,T; H^{-1}(\Omega)).
  \end{equation*}
  Moreover, the solution depends continuously on the data $(f,\bu_0)$ and on the constitutive law $p$ with respect to
  the corresponding weak topologies.
\end{theorem}

Although the spatial dimension is not specified explicitly in~\cite{ChaventLemonnier1974}, assumptions~(4)--(12) therein rely
only on standard integration by parts, ellipticity, and maximum principle arguments. Consequently, the analysis applies
to bounded domains $\Omega \subset \RR^d$ for arbitrary $d \ge 1$. Also, the assumption $\bu(\cdot,0) \in L^\infty(\Omega)$
is imposed to guarantee uniform bounds on the solution via a maximum principle. Related well-posedness results for nonlinear
diffusion equations can be found in the modern monograph of V\'azquez~\cite{Vazquez2007}.

\begin{corollary}[Identifiable support of the constitutive law]
  \label{cor:CL_identifiable_support}
  Under the assumptions of Theorem~\ref{thm:CL_wellposed}, the solution $\bu$ satisfies uniform bounds
  \begin{equation*}
    u_{\min} \le \bu(\bx,t) \le u_{\max} \qquad \text{for a.e. } (\bx,t)\in Q,
  \end{equation*}
  where $u_{\min}$ and $u_{\max}$ depend only on the initial data, the source term, and the constants $p_{\min}$ and $p_{\max}$.
  As a consequence, only values of the constitutive law $p(\zeta)$ for $\zeta \in [u_{\min},u_{\max}]$ can influence the observations
  and, hence, be identified from data.
\end{corollary}

\noindent
This structural restriction underlies the inverse problem formulation adopted in Section~\ref{sec:problem_formulation} and
plays a central role in the analysis of identifiability and regularization discussed later in this review. This observation
becomes particularly transparent at the level of the first variation.

\begin{proposition}[Support of parameter sensitivity {\cite{ChaventLemonnier1974}}]
  \label{prop:CL_sensitivity_support}
  Let $[u_{\min},u_{\max}]$ denote the essential range of $\bu$ in $Q$. Then the first variation of the objective functional
  $\mJ(p)$ vanishes for any perturbation $\delta p$ supported outside the attained state range, satisfying
  \begin{equation*}
    \operatorname{supp}(\delta p) \cap [u_{\min},u_{\max}] = \emptyset.
  \end{equation*}
  Equivalently, the functional $\mJ(p)$ is insensitive to perturbations of the constitutive law $p(\zeta)$ for \(\zeta \notin [u_{\min},u_{\max}]\).
\end{proposition}

\noindent
This property can be interpreted, in modern terms, as a support restriction on the adjoint-weighted sensitivity of the objective
functional with respect to the constitutive law.

It is important to emphasize that the contribution of Chavent and Lemonnier~\cite{ChaventLemonnier1974} does not lie in the
derivation of new existence or uniqueness results for quasilinear parabolic equations. Rather, their analysis relies explicitly
on classical well-posedness theory to guarantee the existence, uniqueness, and boundedness of the forward solution associated
with a given state-dependent parameter $p(\bu)$~\cite{Lions1968,Lions1969,JosephLundgren1972,Vazquez2007}. Under standard
assumptions of uniform ellipticity \eqref{eq:CL_bounds} and appropriate regularity of the data, the forward problem admits
a unique weak solution $\bu$ satisfying uniform $L^\infty(Q)$ bounds. These properties, established in earlier works on nonlinear
parabolic equations, constitute essential prerequisites for the inverse problem framework developed in~\cite{ChaventLemonnier1974}.

One of the most influential contributions of~\cite{ChaventLemonnier1974} is the rigorous analysis of differentiability of the
parameter-to-state map with respect to the nonlinear constitutive law $p$. This result provides the analytical foundation for
sensitivity analysis and gradient-based optimization in inverse problems with state-dependent parameters.

\begin{theorem}[Differentiability of the parameter--to--state map {\cite{ChaventLemonnier1974}}]
  \label{thm:CL_frechet}
  Under the assumptions of Theorem~\ref{thm:CL_wellposed}, and assuming in addition sufficient spatial regularity of the state variable
  (in particular $\nabla \bu \in L^\infty(Q)$ and spatial dimension $d<3$)\footnote{These are standard sufficient conditions in modern
  analyses but not stated explicitly in~\cite{ChaventLemonnier1974}.}, the parameter--to--state map
  \begin{equation*}
    p \mapsto \bu(p)
  \end{equation*}
  is G\^ateaux differentiable\footnote{Under the additional regularity assumptions stated above, this directional differentiability
  may be interpreted as Fr\'echet differentiability in appropriate Banach spaces; however, such functional-analytic refinements are
  not made explicit in~\cite{ChaventLemonnier1974}.} as a map into $L^2(Q)$, with a derivative that is linear and continuous in $\delta p$.

  For any admissible perturbation $\delta p$, the directional derivative $\delta \bu := D\bu(p)[\delta p]$ exists and is characterized
  as the unique weak solution of the linearized parabolic problem
  \begin{subequations}
    \label{eq:CL_linearized_full}
    \begin{align}
      \partial_t \delta \bu - \bnabla \cdot \bigl( p(\bu)\, \bnabla \delta \bu \bigr) &= \bnabla \cdot \bigl( \delta p(\bu)\, \bnabla \bu \bigr)
      && \text{in } Q, \label{eq:CL_linearized} \\
      \bnabla \delta \bu \cdot \bn &= 0 && \text{on } \partial\Omega \times (0,T), \\
      \delta \bu(\cdot,0) &= 0 && \text{in } \Omega.
    \end{align}
  \end{subequations}
  Moreover, the map $\delta p \mapsto \delta \bu$ is linear and continuous from $L^\infty(\mathbb{R})$ into $L^2(Q)$.
\end{theorem}

\noindent
The derivation of~\eqref{eq:CL_linearized}\footnote{This tangent equation is not written explicitly in~\cite{ChaventLemonnier1974} but is
reconstructed from the variational differentiability proof in Section~II of that work.} relies crucially on the uniform $L^\infty(Q)$ bounds
established in Theorem~\ref{thm:CL_wellposed}, which ensure that the compositions $p(\bu)$ and $\delta p(\bu)$ are well defined and bounded.
Moreover, in view of Corollary~\ref{cor:CL_identifiable_support}, the sensitivity $\delta \bu$ depends only on the values
of $\delta p$ over the attained state range $[u_{\min},u_{\max}]$.

This result extends earlier work on nonlinear heat conduction problems~\cite{CannonDuChateau1973} and represents one of
the earliest rigorous treatments of differentiation with respect to an infinite-dimensional, nonparametric constitutive law.
In modern terminology, \eqref{eq:CL_linearized} defines the tangent equation associated with the inverse problem and serves
as the natural starting point for adjoint-based gradient formulations.

Having established well-posedness, identifiability, and differentiability of the forward map, we turn to the adjoint formulation
associated with the linearized problem~\eqref{eq:CL_linearized}, which provides an efficient characterization of the first variation
of the objective functional with respect to the constitutive law. By recasting the inverse problem as an optimal control
problem~\cite{Lions1968,Cea1971}, Chavent and Lemonnier~\cite{ChaventLemonnier1974} derive an adjoint equation for a backward-in-time
variable $\blambda$, whose introduction eliminates the explicit dependence on the state sensitivity $\delta \bu$.

\begin{theorem}[Adjoint equation and gradient representation {\cite{ChaventLemonnier1974,Lions1968,Cea1971}}]
  \label{thm:CL_adjoint}
  Under the assumptions of Theorems~\ref{thm:CL_wellposed} and~\ref{thm:CL_frechet}, let $\bu=\bu(\bx,t;p)$ denote the solution of the
  state equation~\eqref{eq:CL_forward_u}. Then there exists a unique adjoint variable
  \begin{equation*}
    \blambda \in L^2(0,T; H^1(\Omega)) \cap C([0,T]; H^{-1}(\Omega)),
  \end{equation*}
  satisfying the backward (in time) parabolic problem\footnote{The adjoint equation is understood in a weak (or very weak) sense, consistent
  with the regularity assumptions required for differentiability.}
  \begin{equation}
    \begin{aligned}
      - \partial_t \blambda - \bnabla \cdot \bigl( p(\bu)\, \bnabla \blambda \bigr) &= -2 \sum_{i=1}^M
      \left(\frac{1}{|\Omega_i|} \int_{\Omega_i} \bu(\bx,t)\, d\bx - \bz_i(t) \right)
      \frac{\mathbf{1}_{\Omega_i}}{|\Omega_i|} && \text{in } Q, \\
      \bnabla \blambda \cdot \bn \, &= \, 0 && \text{on } \partial\Omega \times (0,T), \\
      \blambda(\cdot, T) \, &= \, 0 && \text{in } \Omega.
    \end{aligned}
    \label{eq:CL_adjoint}
  \end{equation}
  Moreover, the first variation of the objective functional admits the representation
  \begin{equation}
    \mJ'(p)[\delta p] = \int_Q \delta p\bigl(\bu(\bx,t)\bigr) \, \bnabla \bu(\bx,t)\cdot \bnabla \blambda(\bx,t)\, d\bx\,dt.
    \label{eq:CL_gradient_u}
  \end{equation}
\end{theorem}

\begin{corollary}[State-space structure of the gradient]
  \label{cor:CL_state_gradient}
  The gradient of the objective functional is naturally defined with respect to the state variable $\bu$, rather than the physical
  variables $(\bx,t)$. In particular, perturbations of the constitutive law $p(\zeta)$ outside the attained state range
  $[u_{\min},u_{\max}]$ do not contribute to the first variation $\mJ'(p)$.
\end{corollary}

The adjoint representation~\eqref{eq:CL_gradient_u} shows that the sensitivity of the objective functional depends only on the
values of the constitutive law along the solution trajectory $\bu(Q)$. This observation is consistent with the identifiability
restriction established in Corollary~\ref{cor:CL_identifiable_support} and provides a precise analytical explanation for why
the inverse problem is intrinsically posed in a functional space over the state variable.

In modern terminology, the adjoint problem~\eqref{eq:CL_adjoint} enables the evaluation of the derivative of the reduced objective
functional without explicitly computing the state sensitivity $\delta \bu$. Adjoint-based gradient formulations of this type have
since become standard in PDE-constrained optimization and inverse problems~\cite{Gunzburger2003,Hinze2009,Troltzsch2010,BorziSchulz2011}.
From this perspective, the work of Chavent and Lemonnier~\cite{ChaventLemonnier1974} anticipates by several decades the
reduced-space methodology that now underpins most contemporary analytical and numerical approaches to nonlinear inverse
problems governed by PDEs.

To counteract the intrinsic ill-posedness of functional parameter identification, Chavent and Lemonnier introduce admissible sets
incorporating positivity and smoothness constraints, including bounds on higher-order derivatives of $p$. This anticipates later
developments in variational regularization~\cite{TikhonovArsenin1977,EnglHankeNeubauer1996} and convex analysis-based formulations
of inverse problems~\cite{BurgerOsher2004}.

While~\cite{ChaventLemonnier1974} includes numerical illustrations, its enduring impact lies in the analytical structure it establishes.
Modern approaches—such as reduced-space methods, adjoint-based gradient optimization, and operator-learning frameworks—can be viewed
as computational realizations or extensions of the same principles. In particular, contemporary learning-based methods discussed later
in this review differ primarily in how the functional space for $p(\bu)$ is parameterized, while relying implicitly on the same adjoint
sensitivity structure identified in \eqref{eq:CL_gradient_u}. From a PDE-theoretic standpoint, extensions to degenerate or nonlinear
diffusion regimes, such as porous-medium-type operators, require substantially different analytical tools~\cite{Vazquez2007}.

A further conceptual contribution of~\cite{ChaventLemonnier1974} lies in its clear separation between identifiability and stability.
While the analysis establishes that only values of the constitutive law over the attained state range can influence the data, it does
not claim unconditional stability of the inverse problem. Instead, admissible sets are introduced both as a regularization mechanism
and as a means of encoding prior physical knowledge about the expected structure of $p$, an interpretation that resonates with modern
physics-informed and learning-based inverse approaches.

{\bf Relation to identifiability intervals.} The essential range $[u_{\min},u_{\max}]$ appearing in the analysis of Chavent and
Lemonnier~\cite{ChaventLemonnier1974} is closely related to the \emph{identifiability interval} $\mI$ introduced in Definition~\ref{def:intI}.
Both characterize the subset of the state space over which the constitutive law $p(\zeta)$ can influence the data and, consequently,
be identified. The distinction lies primarily in the mode of definition: whereas $[u_{\min},u_{\max}]$ is defined in an essential
(almost everywhere) sense, reflecting measure-theoretic properties of weak solutions, the interval $\mI$ is defined via pointwise
extrema and is tailored to the analysis of reconstruction and measurement design. Under additional regularity assumptions ensuring
continuity of the state variable, the two intervals coincide; in general, $\mI$ provides a conservative, operational approximation
of the essential range relevant for numerical identification.

\subsubsection{Further analytical properties for elliptic problems}
\label{sec:BukshtynovProtasVolkov2011}

Beyond the general parabolic framework established by Chavent and Lemonnier~\cite{ChaventLemonnier1974}, additional analytical
properties have been rigorously derived for specific classes of PDEs involving \emph{state-dependent} constitutive relations.
An example is the stationary heat conduction problem with temperature-dependent conductivity studied in~\cite{Bukshtynov2011},
which serves as a good elliptic prototype within the general framework introduced in Section~\ref{sec:problem_formulation}.
More specifically, \cite{Bukshtynov2011} considers the boundary value problem \eqref{eq:heat} where the unknown constitutive
relation $k = k(T)$ represents a temperature-dependent thermal conductivity. In contrast to classical coefficient identification
problems, the parameter $k$ depends on the solution itself, placing~\eqref{eq:heat} squarely within the class of state-dependent
inverse problems considered in this review. A key analytical contribution of~\cite{Bukshtynov2011} is the establishment of
a \emph{minimum principle} for the state variable, which plays a role analogous to the $L^\infty$ bounds employed
by~\cite{ChaventLemonnier1974} in the parabolic setting.

\begin{theorem}[Minimum principle for temperature-dependent conductivity {\cite{Bukshtynov2011,Grisvard1985}}]
\label{thm:BPV_minimum}
Let $T \in C^2(\Omega)\cap C^0(\overline{\Omega})$ be a solution of the boundary value problem \eqref{eq:heat}. Assume that the
constitutive relation satisfies $k(T)>0$ for all $T \in \mI$, and that the source term satisfies $f(\bx) > 0$ for all
$\bx \in \Omega$. Then
\begin{equation*}
  \min_{\bx\in\overline{\Omega}} T(\bx) = \min_{\bx\in\partial\Omega} T_0(\bx),
\end{equation*}
that is, the minimum of the temperature field $T(\bx)$ is attained on the boundary $\partial\Omega$.
\end{theorem}

\noindent
Theorem~\ref{thm:BPV_minimum} establishes a one-sided bound on the range of the state variable for elliptic problems with
state-dependent conductivity. In the presence of a nonzero source term $f$, the classical maximum principle applies only in
one direction: for $f>0$, the minimum of the solution is attained on the boundary, whereas interior maxima may exceed the
boundary values. As a consequence, the constitutive law $k(T)$ is identifiable only for values of $T$ satisfying
\begin{equation*}
  T \ge \min_{\partial\Omega} T_0,
\end{equation*}
with no a priori upper bound determined solely by boundary data. This illustrates that, for stationary problems,
identifiability intervals may be inherently one-sided unless additional structural assumptions (e.g., $f\equiv 0$ or
sign-changing sources) are imposed.

From this perspective, Theorem~\ref{thm:BPV_minimum} provides an elliptic analogue of the identifiability restriction
identified in~\cite{ChaventLemonnier1974}, while highlighting a fundamental difference between elliptic and parabolic
settings. In the terminology of Section~\ref{sec:adm_set}, the identifiability interval $\mI$ is therefore lower-bounded
by $\min_{\partial\Omega} T_0$, but its upper extent must be determined either from interior measurements or from
additional a priori information.

A second fundamental analytical result of~\cite{Bukshtynov2011} concerns the differentiability of the map
\begin{equation}
  \mF : k \mapsto \left. T(\cdot;k) \right|_{\Sigma},
  \label{eq:FmapT}
\end{equation}
where $\Sigma \subseteq \Omega$ denotes the measurement subdomain on which temperature observations are available, and
the solution $T(\cdot;k)$ is treated as an element of $L^2(\Sigma)$.

\begin{theorem}[Differentiability of the forward map {\cite{Bukshtynov2011}}]
\label{thm:BPV_frechet}
Let $\Omega \subset \RR^d$ be a bounded domain with $C^{1,1}$ boundary. Assume that the conductivity $k \in C^1(\mD)$
satisfies the uniform ellipticity condition
\begin{equation*}
  0 < m_k \le k(\zeta) \le M_k < \infty \qquad \text{for all } \zeta \in \mD,
\end{equation*}
forming admissible set $\mK$, with $m_k$ sufficiently large, and that the corresponding solution $T = T(\bx; k)$
of~\eqref{eq:heat} satisfies
\begin{equation*}
  \|\nabla T\|_{L^\infty(\Omega)} < \infty.
\end{equation*}
Then the map \eqref{eq:FmapT} from $\mK$ into $L^2(\Sigma)$ is Fr\'echet differentiable in the norm $H^1(\mI)$.
\end{theorem}

\noindent
The proof, given in Appendix~B of~\cite{Bukshtynov2011}, follows the variational framework introduced
in~\cite{ChaventLemonnier1974}, but adapts it to a stationary elliptic setting. As in the parabolic case, the use of
appropriate variational and dual formulations (analogous to those employed in~\cite{ChaventLemonnier1974}) shows that
the state sensitivity depends only on values of $\delta k(\zeta)$ for $\zeta$ within the attained state range.

\subsubsection{Other key analytical works on parameter identification}
\label{sec:other_analytical}

The literature addressing inverse problems with \emph{state-dependent} parameters is extensive. In this subsection,
we collect additional analytical contributions (articles and monographs) that establish uniqueness, identifiability,
differentiability, stability, or related structural properties for coefficients of the form $p(\bu)$. Each reference
is accompanied by a concise annotation highlighting its main analytical contribution to the reconstruction of
state-dependent parameters and its relationship to the frameworks discussed in the preceding sections.
Note that throughout this subsection, the term \emph{state-dependent parameter} is used broadly to include diffusion,
reaction, or source terms whose dependence on the solution introduces genuine nonlinearity into the inverse problem.

\textbf{Cannon \& DuCh\^ateau~(1970s--1980s).} A sequence of papers
\cite{CannonDuChateau1973,CannonDuChateau1980a,CannonDuChateau1980b,CannonDuChateau1983,CannonDuChateau1987}
constitutes some of the earliest rigorous analytical work on inverse problems involving coefficients depending explicitly
on the state. These contributions introduced monotonicity arguments, auxiliary inverse problems, and constructive uniqueness
techniques for nonlinear PDEs, thereby laying much of the conceptual groundwork for later analytical developments in
state-dependent parameter identification.

In~\cite{CannonDuChateau1973}, the authors analyze a nonlinear parabolic equation in which both the temporal and diffusive
coefficients depend on the solution itself. Under a structural proportionality assumption between these coefficients, they
prove uniqueness and existence for the coupled inverse problem of reconstructing the constitutive laws from boundary flux
measurements, providing one of the earliest rigorous identifiability results for state-dependent diffusion coefficients.
The paper~\cite{CannonDuChateau1980a} addresses the identification of an unknown nonlinear diffusion coefficient in a quasilinear
parabolic equation from overspecified boundary data. The inverse problem is formulated as an auxiliary variational identification
problem over an admissible class of coefficients, for which existence of minimizers and convergence of finite-dimensional
approximations are established, together with illustrative numerical results.
In~\cite{CannonDuChateau1980b}, attention is shifted to inverse source problems for nonlinear diffusion equations in which
the unknown source term depends on the state variable. Rather than adopting a variational or control-based formulation, the
authors develop a form-free identification approach that treats the source as a genuinely infinite-dimensional unknown. They
establish local-in-time existence and uniqueness of a state-dependent source function consistent with over-specified data,
while explicitly highlighting extreme sensitivity and numerical ill-conditioning of the method. Although the primary unknown
is a source term rather than a diffusion coefficient, the analysis exhibits many structural features common to state-dependent
inverse problems, including identifiability restrictions and strong dependence on the attained state range.
The work~\cite{CannonDuChateau1983} studies a nonlinear hyperbolic inverse problem in which the unknown source term depends
on the solution of a one-dimensional wave equation. Using boundary observations, the authors establish existence results and
analyze qualitative properties of the reconstruction, extending state-dependent inverse problem methodology beyond the parabolic
setting.
Finally, \cite{CannonDuChateau1987} examines the reconstruction of state-dependent diffusion coefficients through the lens of
experimental design and problem formulation. Rather than assuming a parametric form or employing an output least-squares strategy,
the authors propose a form-free reconstruction approach driven by over-specified transient data. Through analytical considerations
and numerical experiments, they show that while exact coefficient recovery is possible in noise-free settings, even small data
perturbations lead only to weak (integral-norm) convergence. This leads to the central conclusion that classical strong-solution
frameworks are generally ill-suited for such inverse problems and that a weak formulation (placing the coefficient in an
integral-norm space and admitting weak PDE solutions) is more appropriate.

\textbf{Egger, Pietschmann \& Schlottbom.}
In a modern analytical treatment, Egger \emph{et al.}~\cite{Egger2014} revisit and substantially extend classical results of Cannon
and DuCh\^ateau for nonlinear heat conduction laws. While the earlier works established identifiability and uniqueness of
state-dependent conductivities (typically relying on monotonicity arguments and comparison principles), \cite{Egger2014} provides
a new energy-based proof of Cannon's uniqueness result for the stationary quasilinear elliptic problem and derives corresponding
\emph{quantitative stability estimates}. The analysis is then extended to the quasilinear parabolic case by treating the
time-dependent problem as a perturbation of the stationary setting under appropriate experimental design assumptions. This work
forms a rigorous bridge between classical maximum-principle-based identifiability theory and modern stability-oriented analysis for
PDE-constrained inverse problems with state-dependent coefficients.

Building on this stability-oriented reformulation of Cannon's uniqueness paradigm, the subsequent work~\cite{Egger2017} investigates
the robustness of state-dependent parameter identifiability under model extensions. In particular, the authors establish uniqueness
results for state-dependent diffusion coefficients in quasilinear elliptic and parabolic problems that include additional lower-order
terms. Using monotonicity arguments and energy-based comparison techniques, they extend classical identifiability theory to models
in which transport or reaction terms are present alongside nonlinear diffusion. A key analytical contribution is the demonstration
that uniqueness of the constitutive law $p(\bu)$ can persist under controlled lower-order perturbations, provided suitable structural
conditions are satisfied. This work shows that identifiability of state-dependent parameters is not restricted to pure diffusion
prototypes and remains valid for a broader class of quasilinear models relevant in applications.

\textbf{Additional analytical developments~(1980s–2000s).}
Following the foundational framework established by~\cite{ChaventLemonnier1974}, a broad analytical literature emerged addressing
inverse problems for quasilinear and nonlinear partial differential equations in which unknown coefficients depend explicitly on the
solution. Early variational and monotonicity-based uniqueness results for nonlinear conductivity identification were developed by
Kohn and Vogelius~\cite{KohnVogelius1984,KohnVogelius1987}, introducing comparison principles and relaxed variational formulations
that later proved influential in nonlinear coefficient reconstruction.

A parallel analytical tradition developed from the late 1980s onward in the study of inverse problems for quasilinear and nonlinear
parabolic equations, with particular emphasis on coefficients entering the governing equations \emph{nonlinearly} (often through the solution).
Monographic treatments and research articles by Prilepko, Yamamoto, and collaborators~\cite{Prilepko2000,ChoulliYamamoto1996,ChoulliYamamoto1997}
established local uniqueness and conditional stability results for nonlinear inverse problems using linearization techniques, Carleman
estimates, and continuation arguments. Although these works are not always phrased explicitly in terms of constitutive laws of the form
$p(\bu)$, they rigorously address inverse problems in which the unknown coefficient depends nonlinearly on the state variable, thereby
extending classical inverse problem theory to genuinely quasilinear settings.

More recently, \cite{Brunk2023} investigated identifiability and uniqueness for nonlinear parabolic equations of Cahn--Hilliard type with
\emph{state-dependent mobility} and \emph{chemical potential} terms. Focusing on the simultaneous reconstruction of multiple constitutive
functions entering the governing equation nonlinearly, the authors analyze scaling invariances, conditional uniqueness, and observability
requirements under various measurement scenarios. Rather than adopting a reduced-space variational formulation, the work is framed in terms
of equation-error inversion and operator identifiability, emphasizing structural limitations that persist independently of numerical
discretization or optimization strategy. These results complement the classical theory of Cannon and DuCh\^ateau and the stability-oriented
analysis of Egger \emph{et al.} by clarifying when state-dependent constitutive laws are uniquely recoverable \emph{in principle}, even before
considering adjoint-based optimization or regularization.

Related analytical investigations by Cannon and DuCh\^ateau~\cite{CannonDuChateau1980a,CannonDuChateau1980b} examined identifiability
and instability phenomena for nonlinear diffusion coefficients from overdetermined boundary and interior observations, providing early
demonstrations of the severe ill-posedness inherent in such inverse problems. Complementary functional-analytic studies by Engl, Plato,
and others~\cite{EnglKunischNeubauer1989,Plato1995} further clarified the instability mechanisms and regularization requirements for
nonlinear ill-posed operator equations, supplying an abstract framework for understanding the limitations of coefficient reconstruction
beyond the PDE-specific setting. In parallel, developments in shape calculus and PDE-constrained
optimization~\cite{DelfourZolesio2011,SokolowskiZolesio1992,HaslingerNeittaanmaki1996} established rigorous adjoint-based differentiation
and sensitivity analysis tools applicable to inverse problems with state-dependent material laws. Together, these contributions complement
the variational framework of~\cite{ChaventLemonnier1974} and form an important bridge between early identifiability theory and later
quantitative stability and regularization analyses for state-dependent inverse problems.

\textbf{V\'azquez~(2007).}
The monograph~\cite{Vazquez2007} develops the fundamental PDE theory for degenerate and nonlinear diffusion models, including porous-medium
and fast-diffusion equations in which the diffusion coefficient depends on the state variable. Although inverse problems and parameter
reconstruction are not addressed explicitly, this body of work provides the essential analytical foundation for any study of
$p(\bu)$-reconstruction in degenerate regimes. In particular, the monograph characterizes existence, uniqueness, regularity,
finite-speed-of-propagation, extinction phenomena, and qualitative behavior of solutions when the constitutive law $p(\bu)$ may vanish
or become singular. These features differ fundamentally from the uniformly elliptic setting assumed in most classical inverse problem
analyses. Consequently, any identifiability, differentiability, or stability analysis for state-dependent parameter reconstruction
in degenerate diffusion models must necessarily build upon the theory developed in~\cite{Vazquez2007}.

\textbf{Isakov (2017).}
The monograph~\cite{Isakov2017} is a comprehensive and authoritative reference on the analytical theory of inverse problems for partial
differential equations. It systematically develops the principal tools used in uniqueness and stability analysis, including Carleman
estimates, energy methods, linearization techniques, and conditional stability frameworks, for a wide range of elliptic, parabolic,
and hyperbolic PDEs. While the book is not devoted exclusively to coefficients depending on the state variable, it treats nonlinear and
quasilinear inverse problems throughout and provides many results directly relevant to the identification of coefficients entering the
governing equations in a nonlinear manner. In the context of $p(\bu)$-reconstruction, Isakov’s work serves as a foundational methodological
source. Many of the uniqueness arguments, continuation principles, and conditional stability estimates presented in the book have been
adapted, explicitly or implicitly, in the analysis of state-dependent inverse problems, including models involving nonlinear diffusion
or reaction terms. As such,~\cite{Isakov2017} complements the specialized analytical results discussed above by providing a unifying
framework within which $p(\bu)$-reconstruction problems can be rigorously formulated and analyzed.

\textbf{Methodological and regularization literature.}
Although not focused exclusively on state-dependent coefficients of the form $p(\bu)$, the monographs and papers
\cite{TikhonovArsenin1977,EnglHankeNeubauer1996,BurgerOsher2004,Hinze2009,Troltzsch2010,BorziSchulz2011} provide the foundational
regularization, convex-analysis, and PDE-constrained optimization frameworks that underpin much of the modern analysis and numerical
treatment of inverse problems with state-dependent parameters. These works justify admissible set formulations, variational
regularization strategies, and reduced-space optimization approaches that appear, either explicitly or implicitly, in the analyses of
state-dependent parameter identification. They also supply practical computational methodologies, including Tikhonov regularization,
Sobolev-gradient techniques, and adjoint-based reduced formulations, which are routinely combined with the analytical results surveyed
above to obtain stable numerical reconstructions of nonlinear constitutive laws.

\textbf{Scope and open analytical challenges.}
The references surveyed in Section~\ref{sec:anal_properties} demonstrate that a substantial analytical theory now exists for inverse
problems with coefficients depending on the state variable in scalar elliptic and parabolic equations. This theory addresses key
structural questions including well-posedness of the forward problem, identifiability intervals determined by the attained state range,
Fr\'echet differentiability of the parameter-to-state map, adjoint-based sensitivity representations, uniqueness, and, more recently,
quantitative stability estimates. Related, though more limited, analytical results are also available for certain hyperbolic problems
with state-dependent source terms or coefficients.

By contrast, a comparable analytical framework for inverse problems with state-dependent parameters in \emph{systems} of partial
differential equations, or in strongly nonlinear models such as the Navier--Stokes equations in their general form (even in forward
well-posed regimes), is largely absent. In particular, problems in which the constitutive law depends on one or more components of
the solution vector, as in coupled transport, flow, or reactive systems of the type considered in~\cite{Bukshtynov2013}, pose substantial
unresolved challenges. These include the lack of maximum principles, limited regularity of solutions, nontrivial coupling between state
variables, and the breakdown of classical monotonicity and comparison arguments. At present, the development of rigorous identifiability,
differentiability, and stability theory for state-dependent parameter reconstruction in such systems remains an open and largely unexplored
direction of research, motivating the developments discussed in the remainder of this review.

\subsection{Variational and Adjoint Framework}
\label{sec:var_adj_framework}

As discussed in Sections~\ref{sec:problem_formulation} and \ref{sec:anal_properties}, inverse problems involving
\emph{state-dependent parameters} are characterized by the fact that the unknown constitutive quantity $p = p(\bu)$ enters the governing
PDE system \eqref{eq:general_PDE}--\eqref{eq:general_BC} through a functional dependence on the solution $\bu(\bx, t; p)$. This coupling
induces a nonlinear and typically compact parameter--to--observation map \eqref{eq:Fmap} and leads, in general, to an ill-posed inverse
problem.

From a computational perspective, several strategies have been developed to address state-dependent parameter reconstruction. A widely
used approach is to introduce an explicit \emph{parameterization} of the constitutive relation $p(\bu)$, for example via low-order
polynomials, spline representations, or other finite-dimensional bases, and to estimate the corresponding coefficients,
cf.~\cite{Alifanov2004,Alifanov2007}. This results in a finite-dimensional optimization problem that can be treated with standard numerical
techniques, but may introduce approximation bias and limit resolution unless a sufficiently rich representation is employed.

An alternative formulation treats the constitutive relation $p(\bu)$ itself as an element of an infinite-dimensional function space.
In this setting, the admissible set $\mP \subset \mX$, defined in \eqref{eq:Pset}, is chosen to reflect both analytical regularity and
physical constraints, such as positivity or boundedness, as discussed in Section~\ref{sec:anal_properties}. This viewpoint avoids
an \emph{a priori} reduction of dimensionality and allows the reconstruction to adapt to the information content of the data. Recent
physics-informed and operator-learning approaches may be viewed as data-driven realizations of this paradigm, although the present
review focuses mainly on variational formulations grounded in classical PDE analysis.

Most computational approaches for state-dependent parameter estimation ultimately reduce to the solution of an optimization problem
constrained by the governing PDE system. Two principal paradigms are commonly distinguished. In the \emph{discretize--then--optimize}
approach, the PDE is first discretized and sensitivities are derived for the resulting finite system. In contrast, the
\emph{optimize--then--discretize} approach derives first-order optimality conditions at the continuous level prior to discretization.
The distinction between these approaches has been recognized in the PDE-constrained optimization community since at least the early 2010s
(e.g., \cite{HinzeTroltzsch2010,Hinze2009}), and the terminology has been formalized and widely discussed in the literature on adjoint
computation for discretized PDEs (e.g., \cite{Wilcox2015,Bukshtynov2011,Bukshtynov2013}). The latter strategy preserves the variational
structure of the problem and yields adjoint formulations that are independent of the spatial and temporal discretization. This
framework is adopted in what follows and forms the basis for the adjoint-based gradient evaluations reviewed below.

\subsubsection{Parameter estimation as an optimization problem}
\label{sec:opt_problem}

Let $\bu$ denote the state variable governed by the PDE system \eqref{eq:general_PDE}--\eqref{eq:general_BC}, and let $p(\bu) \in \mP$
denote a state-dependent parameter to be reconstructed. For a given $p$, the corresponding state $\bu = \bu(p)$ is defined implicitly
through the governing equations, and the associated observations $\by_{\text{obs}} \in \mY$ are obtained via the observation operator
$\mO$, as defined in~\eqref{eq:observ}. Together, these define the parameter--to--observation map $\mF : \mP \rightarrow \mY$,
introduced by~\eqref{eq:Fmap}--\eqref{eq:inverse_eq}. The inverse problem consists in recovering a parameter $p$ such that
$\mF(p) \approx \by_{\text{obs}}$. Following the standard variational formulation introduced in Section~\ref{sec:PDE_opt}, this problem
is recast as the minimization of the least-squares cost functional $\mJ(\bu, p(\bu))$, cf.~\eqref{eq:cost_general}, leading to the
optimization problem \eqref{eq:final_opt}.

In many applications, the admissible set $\mP  \subset \mX$ incorporates inequality constraints reflecting physical principles, such as
positivity of material coefficients or thermodynamic admissibility. These constraints may be enforced explicitly or incorporated through
variable transformations~\cite{BoydVandenberghe2004}, barrier methods~\cite{Ruszczynski2006}, or penalty formulations~\cite{Vogel2002}.
Independently of the chosen strategy, problem~\eqref{eq:final_opt} is generally nonconvex, and numerical algorithms are therefore expected
to converge only to local minimizers. The first-order optimality condition associated with~\eqref{eq:final_opt} requires the directional
derivative $\mJ'(p; \delta p)$ of the cost functional $\mJ(p)$ with respect to the control variable $p$ to vanish at $\hat{p}$ for all
perturbations $\delta p$~\cite{Luenberger1976}, i.e.,
\begin{equation}
  \forall \delta p \in \mX, \qquad \mJ'(\hat{p}; \delta p) = 0.
  \label{eq:opt_cond_general}
\end{equation}

In computational practice, the \emph{(local) optimizer} $\hat{p}$ is obtained using iterative gradient-based algorithms of the form
\begin{equation}
  \label{eq:gd_general}
  \begin{aligned}
    p^{(n+1)} &= p^{(n)} - \tau^{(n)} \, \bnabla_p \mJ(p^{(n)}), \\
    p^{(0)} &= p_0,
  \end{aligned}
\end{equation}
where $\bnabla_p \mJ(p)$ denotes the gradient of the cost functional $\mJ(p)$ with respect to $p$, $\tau^{(n)}$ is a suitably chosen step
length, and $p_0$ is an \emph{initial guess} chosen based on prior knowledge or modeling assumptions. In large-scale settings,
conjugate-gradient and quasi-Newton methods are commonly employed to reduce the number of forward and adjoint PDE solves per
iteration~\cite{Nocedal2006,Vogel2002}.

The dominant computational cost in gradient-based algorithms is the evaluation of the gradient $\bnabla_p \mJ$. For state-dependent
parameters, this quantity represents the sensitivity of the misfit functional with respect to perturbations of an entire function, rather
than a finite-dimensional parameter vector. Efficient and accurate computation of $\bnabla_p \mJ$ is therefore critical for the practical
solution of problem~\eqref{eq:final_opt}, especially in large-scale or time-dependent settings.

The evaluation of cost functional gradients is based on a variational formulation of the governing PDE system and the introduction of adjoint
variables as Lagrange multipliers enforcing the state equations. Within the \emph{optimize--then--discretize} paradigm adopted here, all
differentiation is carried out at the continuous level, preserving the variational structure of the problem. For state-dependent parameters,
variations of the control induce implicit variations of the state through the constitutive relation $p(\bu)$, leading to adjoint systems that
differ from those arising in problems with state-independent controls. The resulting adjoint-based gradient expressions are independent of the
dimension of the observation space and have a computational cost comparable to that of a single forward solve. Their derivations are reviewed
in the next subsection.

\subsubsection{Cost functional gradients via variational (adjoint-based) analysis}
\label{sec:var_adj_grads}

The derivation of cost functional gradients for inverse problems with \emph{state-dependent parameters} is far from straightforward and differs
fundamentally from the case of parameters that depend explicitly on space or time. In the state-dependent setting, the structure of the adjoint
system and the resulting gradient expression depend sensitively on the analytical structure of the differential operator $\mG$ appearing in the
governing equations~\eqref{eq:general_PDE}--\eqref{eq:general_BC}, as well as on the manner in which the constitutive relation $p(\bu)$
enters the PDE system. By contrast, for spatially dependent parameters $p(\bx)$ and space-time dependent parameters $p(\bx,t)$, adjoint-based
gradient derivations follow a standard and largely problem-independent pattern and are treated in many tutorial and textbook
presentations~\cite{BukshtynovBook2023,Hinze2009,Troltzsch2010,Nocedal2006,Gunzburger2003}.

As discussed in detail in~\cite{Bukshtynov2011,Bukshtynov2013}, the analytical derivation of gradients for inverse problems with state-dependent
parameters begins with the introduction of an adjoint variable $\blambda$ and an associated \emph{Lagrangian functional} defined as
\begin{equation}
  \label{eq:Lagrangian_state_dep}
  \mL(\bu,p(\bu),\blambda) = \mJ(\bu,p(\bu)) + \big\langle \blambda, \, \mG(\bu,p(\bu)) - f \big\rangle_{\mX},
\end{equation}
where $\langle\cdot,\cdot\rangle_{\mX}$ denotes the appropriate duality pairing between the adjoint space and the state space, here taken to be
a Hilbert space $\mX$ (here taken, for simplicity, as $\mX = L_2$). Assuming that the partial Fr\'echet derivatives
$\mJ_{\bu} = \partial \mJ / \partial \bu$, $\mJ_{p} = \partial \mJ / \partial p$, $\mG_{\bu} = \partial \mG / \partial \bu$, and
$\mG_{p} = \partial \mG / \partial p$ exist in the appropriate dual spaces, the first variation of the Lagrangian \eqref{eq:Lagrangian_state_dep}
can be written as
\begin{equation}
  \label{eq:first_variation_L}
  \delta \mL(\bu,p(\bu),\blambda)
  = \left\langle \mJ_{\bu}, \delta\bu \right\rangle_{\mX}
  + \left\langle \mJ_{p}, \delta p \right\rangle_{\mX}
  + \left\langle \blambda,\, \mG_{\bu}\,\delta\bu
  + \mG_{p}\,\delta p \right\rangle_{\mX}
  + \left\langle \mG(\bu,p(\bu)) - f,\, \delta\blambda \right\rangle_{\mX}.
\end{equation}
Here $\delta\bu$, $\delta p$, and $\delta\blambda$ denote admissible variations in the state, control, and adjoint variables, respectively.
Although $p$ depends on the state through the constitutive relation $p = p(\zeta(\bu))$, variations $\delta p$ are interpreted in the sense
of Riesz representation~\cite{Berger1977} in the control space, as discussed below.

The central analytical difficulty in deriving explicit expressions for the gradient $\bnabla_p \mJ$ lies in reshaping
\eqref{eq:first_variation_L} into the Riesz form
\begin{equation}
  \label{eq:first_variation_L_reshaped}
  \delta \mL(\bu,p(\bu),\blambda)
  = \left\langle \bnabla_{p} \mJ, \delta p \right\rangle_{\mX}
  + \left\langle \bnabla_{\bu} \mJ, \delta\bu \right\rangle_{\mX}
  + \left\langle \bnabla_{\blambda} \mJ, \delta\blambda \right\rangle_{\mX},
\end{equation}
where $\bnabla_{\bu} \mJ$, $\bnabla_{p} \mJ$, and $\bnabla_{\blambda} \mJ$ denote the Riesz representers of the corresponding partial
derivatives of the Lagrangian. First-order optimality conditions are obtained by requiring \eqref{eq:first_variation_L} to vanish for
all admissible variations $\delta\bu$, $\delta p$, and $\delta\blambda$. Comparison of \eqref{eq:first_variation_L} and
\eqref{eq:first_variation_L_reshaped} shows that setting $\bnabla_{\blambda} \mJ = 0$ recovers the state equation, whereas
$\bnabla_{\bu} \mJ = 0$ yields the adjoint system. In contrast to inverse problems with state-independent controls, the
adjoint equations generally contain additional terms arising from the implicit dependence of the state on the control through the
constitutive relation $p(\bu)$. As a consequence, the precise form of the adjoint system depends sensitively on the analytical structure
of the operator $\mG$ and on the manner in which $p(\bu)$ enters the governing equations.

Once the adjoint equation has been introduced, the gradient $\bnabla_{p} \mJ$ of the reduced cost functional $\mJ(p)$ can be identified.
Specifically, for constitutive relations of the form $p = p(\zeta)$, where $\zeta \in \mD$ is a scalar state variable introduced in
Section~\ref{sec:adm_set}, the gradient evaluated at $\zeta_0 \in \mD$ admits a representation of the form
\begin{equation}
  \label{eq:state_dep_gradient_general}
  \bnabla_p \mJ(\zeta_0) = \int_0^{t_f} \int_{\Omega} \delta\!\left(\zeta(\bx,t) - \zeta_0\right)\,
  F\!\left( \bu,\blambda,\bx,t;\, p(\zeta) \right) \, d\bx \, dt,
\end{equation}
where $\delta(\cdot)$ denotes the Dirac delta distribution and $F$ is a problem-dependent integrand involving the state $\bu$, the
adjoint variable $\blambda$, and derivatives of the operator $\mG$ with respect to the constitutive relation.

The appearance of the Dirac delta reflects the fact that sensitivities with respect to the constitutive relation $p(\cdot)$ are
accumulated over all space-time locations at which the state variable $\zeta$ attains the value $\zeta_0$. Equivalently, the gradient
may be interpreted as an integral over \emph{codimension-one level sets} of the state field $\zeta(\bx,t)$. This structure is
characteristic of inverse problems with state-dependent parameters and fundamentally distinguishes the present formulation from
standard adjoint-based gradients for space- or time-dependent controls, cf.~\cite{Bukshtynov2011,Bukshtynov2013}.

The gradient structure \eqref{eq:state_dep_gradient_general} is consistent with the classical result of Chavent and
Lemonnier~\cite{ChaventLemonnier1974}. In the notation of Theorem~\ref{thm:CL_adjoint}, the gradient of the reduced functional with
respect to the constitutive relation $p = p(\zeta)$ can be written, in this setting, as
\begin{equation}
  \label{eq:CL_gradient_rewritten}
  \bnabla_p \mJ(\zeta_0) = \int_{\Omega \times (0,T)} \delta\!\left(\zeta(\bx,t)-\zeta_0\right)\,
  \bnabla \bu(\bx,t)\cdot \bnabla \blambda(\bx,t) \, d\bx \, dt,
\end{equation}
which is precisely of the form \eqref{eq:state_dep_gradient_general} with the integrand
$F(\bu,\blambda,\bx,t;p(\zeta)) = \bnabla \bu(\bx,t)\cdot \bnabla \blambda(\bx,t)$.\footnote{The overall sign of the integrand depends
on the convention adopted for the adjoint equation; the present form is consistent with Theorem~\ref{thm:CL_adjoint} and the original
formulation of~\cite{ChaventLemonnier1974}.} Thus, the Chavent--Lemonnier gradient may be interpreted as a particular instance of the
general adjoint-based representation derived above, with sensitivities accumulated over space-time level sets of the state variable.

\paragraph{Heat conduction with temperature-dependent conductivity.} We now turn to the time-independent inverse heat conduction
problem~\eqref{eq:heat} considered in~\cite{Bukshtynov2011}, which fits naturally into the framework of state-dependent parameters
introduced in Section~\ref{sec:problem_formulation}. Specifically, the inverse problem consists in identifying the constitutive
relation $k = k(T)$ from temperature measurements $\tilde T(\bx)$. The following result summarizes the gradient representation derived
in~\cite{Bukshtynov2011}.

\begin{theorem}[Gradient for the state-dependent conductivity problem {\cite{Bukshtynov2011}}]
\label{thm:Bukshtynov2011_gradient}
Let $\Omega$ be a smooth bounded open set and $\delta k \in L_2([T_a,T_b])$. Assume that $f$ and the corresponding state $T$ are
sufficiently smooth. Then the directional derivative of the cost functional
\begin{equation*}
  \mJ(k) = \frac{1}{2}\int_\Sigma \bigl[\tilde T(\bx) - T(\bx; k) \bigr]^2 \, d\bx
\end{equation*}
admits the Riesz representation
\begin{equation}
  \label{eq:B11_rr}
  \mJ'(k; \delta k) = \int_{T_a}^{T_b} \bnabla_k^{L_2} \mJ(s) \, \delta k(s)\, ds,
\end{equation}
with the $L_2$-gradient given by
\begin{equation}
  \label{eq:B11_grad_L2}
  \bnabla_k^{L_2} \mJ(s) = \int_{\Omega} \chi_{[T_\alpha,\,T(\bx)]}(s)\, \Delta T^*(\bx)\, d\bx
  - \int_{\partial\Omega} \chi_{[T_\alpha,\,T(\bx)]}(s)\, \frac{\partial T^*}{\partial n}(\bx)\, d\sigma.
\end{equation}
Here, $\chi_{[a,b]}$ denotes the characteristic function of the interval $[a,b]$, and the adjoint variable $T^*$ solves
\begin{equation}
  \label{eq:B11_adjoint}
  \begin{aligned}
    k(T)\,\Delta T^* &= \bigl[\tilde T(\bx) - T(\bx) \bigr]\,\chi_\Sigma(\bx) \qquad \text{in } \Omega, \\
    T^* &= 0 \qquad \text{on } \partial\Omega,
  \end{aligned}
\end{equation}
where $\chi_\Sigma$ denotes the characteristic function of the sensing domain $\Sigma \subseteq \Omega$.
\end{theorem}

Expression \eqref{eq:B11_grad_L2} can be rewritten in a form directly comparable to the Chavent--Lemonnier gradient and to
\eqref{eq:state_dep_gradient_general}. Interpreting the derivative of the characteristic function in the sense of distributions
and integrating by parts, one obtains
\begin{equation}
  \label{eq:B11_delta_form}
  \bnabla_k^{L_2} \mJ(s) = -\int_{\Omega} \delta\!\bigl(T(\bx)-s\bigr)\, \nabla T(\bx)\cdot\nabla T^*(\bx)\, d\bx,
\end{equation}
which coincides with the (time-independent) Chavent--Lemonnier representation after identifying the state variable $\zeta$ with
the temperature $T$ (up to the adjoint sign convention). In particular, \eqref{eq:B11_delta_form} is of the general form
\eqref{eq:state_dep_gradient_general} with
\begin{equation*}
  F(\bu,\blambda,\bx; p(\zeta)) = -\nabla T(\bx)\cdot\nabla T^*(\bx), \qquad \zeta = T.
\end{equation*}

While the Dirac delta representation \eqref{eq:B11_delta_form} makes the connection with level-set--based sensitivities explicit,
it is not well suited for numerical computations, as it involves differentiation of quantities localized on level sets of the
state. For this reason, \cite{Bukshtynov2011} advocates the equivalent form \eqref{eq:B11_grad_L2}, which avoids explicit delta
distributions and expresses the gradient through volume and boundary integrals involving the adjoint solution. Using the adjoint
equation \eqref{eq:B11_adjoint}, this representation can be further transformed into
\begin{equation}
  \label{eq:B11_comp_form}
  \bnabla_k^{L_2} \mJ(s) = \int_{\Sigma} \chi_{[T_\alpha,\,T(\bx)]}(s)\, \frac{\tilde T(\bx)-T(\bx)}{k(T(\bx))}\, d\bx
  - \int_{\partial\Omega} \chi_{[T_\alpha,\,T(\bx)]}(s)\, \frac{\partial T^*}{\partial n}(\bx)\, d\sigma,
\end{equation}
which is computationally preferable since numerical differentiation with respect to the state variable appears only in boundary
terms. In summary, the gradient expressions derived in~\cite{Bukshtynov2011} provide a concrete realization of the general
structure \eqref{eq:state_dep_gradient_general} for a stationary inverse problem with a state-dependent constitutive law. They
further clarify how the classical level-set--based gradient of Chavent and Lemonnier can be recast into equivalent forms that are
more amenable to stable and efficient numerical implementation.

\paragraph{Incompressible flow with temperature-dependent viscosity.} We conclude this section by considering the coupled
Navier--Stokes--heat inverse problem studied in~\cite{Bukshtynov2013}, introduced in Section~\ref{sec:problem_formulation} and
governed by \eqref{eq:NS_heat}. In this setting, the constitutive relation is given by a temperature-dependent viscosity $\mu(T)$,
and the inverse problem consists in reconstructing this dependence from temperature measurements $\tilde T(\bx, t)$ collected at
discrete sensor locations. The following theorem summarizes the adjoint-based gradient representation derived in~\cite{Bukshtynov2013}
in its most direct ($L_2$) form.

\begin{theorem}[Gradient for the Navier--Stokes problem with state-dependent viscosity {\cite{Bukshtynov2013}}]
\label{thm:Bukshtynov2013_gradient}
Let $\Omega \subset \mathbb{R}^d$ be a sufficiently regular open bounded domain and $\delta \mu \in \mX = L_2 (\mD)$. Assume that
the solutions $\bv$ and $T$ of system~\eqref{eq:NS_heat} are sufficiently smooth. Then the directional derivative of the cost
functional\footnote{here, we consider the pure form of the dependency $\mu = \mu(T)$ as the use of the slack variable $\theta$ used
in the original statement of this theorem in~\cite{Bukshtynov2013} will be discussed later.}
\begin{equation*}
  \mJ(\mu) = \frac{1}{2} \int_0^{t_f} \sum_{i=1}^M \bigl[ T(\bx_i, \tau; \mu) - \tilde T_i (\tau) \bigr]^2 \, d\tau,
\end{equation*}
with respect to the constitutive relation $\mu = \mu (T)$ admits the Riesz representation
\begin{equation}
  \label{eq:B13_rr}
  \mJ'(\mu; \delta \mu) = \int_{-\infty}^{\infty} \bnabla_\mu^{L_2} \mJ(s)\, \delta \mu (s) \, ds,
\end{equation}
where the $L_2$-gradient is given by
\begin{equation}
  \label{eq:B13_grad_L2}
  \bnabla_{\mu}^{L_2} \mJ(s) = - \int_0^{t_f} \int_{\Omega} \delta\!\bigl(T(\bx,\tau)-s \bigr)\, \bigl[\bnabla \bv(\bx,\tau)
  + (\bnabla \bv(\bx,\tau))^{\top}\bigr] : \bnabla \bv^*(\bx,\tau) \, d\bx \, d\tau.
\end{equation}
Here, $\delta(\cdot)$ denotes the Dirac delta distribution and the adjoint variable $\bv^*$ solves the adjoint system associated
with~\eqref{eq:NS_heat}, as given in~\cite{Bukshtynov2013}.
\end{theorem}

\noindent
As before, we note that the gradient representation \eqref{eq:B13_grad_L2} is fully consistent with the general structure derived
in Section~\ref{sec:var_adj_grads} as it can be written in the abstract form provided by \eqref{eq:state_dep_gradient_general} by
identifying the scalar state variable $\zeta \equiv T$ and the problem-specific integrand as
\begin{equation*}
  F(\bu,\blambda,\bx,t;p(\zeta)) = - \, \bigl[\bnabla \bv + (\bnabla \bv)^{\top}\bigr] : \bnabla \bv^*, \qquad \zeta = T.
\end{equation*}

As in the classical result of Chavent and Lemonnier and in the stationary heat conduction problem of~\cite{Bukshtynov2011}, the
Dirac delta factor in~\eqref{eq:B13_grad_L2} reflects the accumulation of sensitivities over space-time level sets of the state
variable $T(\bx,t)$. The Navier--Stokes example further demonstrates that this level-set structure persists in fully coupled,
time-dependent multiphysics systems, thereby confirming the generality of the adjoint-based gradient
representation~\eqref{eq:state_dep_gradient_general} for inverse problems with state-dependent constitutive laws.

\subsubsection{Challenges in moving to practical computations}
\label{sec:grad_challenge}

The adjoint-based framework developed in this section provides a rigorous variational characterization of gradients for inverse
problems with state-dependent constitutive relations. However, translating these continuous-level results into reliable and
efficient numerical algorithms introduces a number of additional challenges. Some of these difficulties are common to nonlinear
PDE-constrained inverse problems, while others are intrinsic to the level-set--based structure of gradients of the
form~\eqref{eq:state_dep_gradient_general}. We summarize below the most important issues that arise in practical computations
and motivate the algorithmic choices discussed in subsequent sections.

\begin{itemize}
  \item \emph{Identifiability and reconstruction range.} In practice, the constitutive relation can only be identified over the
    subset of state values actually attained by solutions of the forward problem during the optimization process. As observed
    in~\cite{Bukshtynov2011,Bukshtynov2013}, the \emph{identifiability interval} $\mI$ is typically a strict subset of the
    prescribed admissible domain $\mD$. Attempts to reconstruct $p(\zeta)$ outside $\mI$ are inherently ill-posed and lead to
    severe instability unless additional prior information is imposed, cf.~Section~\ref{sec:shifting_I}.
  \item \emph{Ill-posedness and sensitivity to noise.} Inverse problems with state-dependent parameters amplify measurement noise
    through both the solution operator and the differentiation inherent in adjoint-based gradients. This effect is particularly
    pronounced when sensitivities accumulate on level sets of the state, making explicit regularization a central component of
    any practical algorithm, cf.~Section~\ref{sec:noisy_data}.
  \item \emph{Low regularity of $L_2$-gradients.} Gradients obtained in the natural $L_2(\mD)$ setting are generally discontinuous
    and may only exist in a distributional sense. As demonstrated in~\cite{Bukshtynov2011,Bukshtynov2013}, such gradients are
    unsuitable for direct use in iterative optimization schemes and must be transformed into smoother Sobolev gradients to ensure
    numerical stability and mesh-independent convergence; cf.~Section~\ref{sec:Sobolev_gradients}.
  \item \emph{Computational complexity of level-set integrals.} The general representation \eqref{eq:state_dep_gradient_general}
    involves Dirac delta distributions supported on codimension-one level sets of the state variable. Numerically, this leads to
    a number of algorithmic and discretization challenges: (i)~accurate approximation of integrals localized near level sets,
    (ii)~sensitivity to the regularity and topology of isolines or isosurfaces, and (iii)~reliance on discrete state and adjoint
    fields that may not resolve these sets adequately on finite meshes, cf.~Section~\ref{sec:level_set_int}.
  \item \emph{Discrete consistency of state and adjoint solvers.} Accurate gradient evaluation requires strict consistency between
    the discretizations of the forward and adjoint problems. Inconsistencies in boundary conditions, stabilization operators, or
    time-integration schemes can lead to incorrect gradients and stagnation or divergence of optimization algorithms,
    cf.~Sections~\ref{sec:two_paradigms} and \ref{sec:gradient_based_OPT}.
  \item \emph{Well-posedness constraints on admissible parameters.} State-dependent constitutive relations must satisfy structural
    constraints (such as positivity, boundedness, or monotonicity) to guarantee well-posedness of the governing PDEs. Enforcing
    these constraints during optimization is essential but nontrivial, particularly in infinite-dimensional function spaces,
    cf.~Section~\ref{sec:barrier_slack}.
  \item \emph{Nonconvexity and local convergence.} The resulting optimization problems are typically nonconvex, with multiple local
    minima. Consequently, gradient-based algorithms exhibit only local convergence and are sensitive to initialization, a feature
    clearly illustrated in the numerical studies of~\cite{Bukshtynov2011,Bukshtynov2013}.
  \item \emph{Verification of adjoint gradients.} Given the analytical and implementation complexity of adjoint systems for
    state-dependent parameters, systematic verification of gradient computations (e.g., via Taylor series remainder or
    finite-difference tests) is essential to ensure correctness and robustness of the optimization process,
    cf.~Section~\ref{sec:verif_kappa}.
\end{itemize}

These challenges highlight the gap between variational gradient formulas derived at the continuous level and their practical
numerical realization. They motivate the introduction of Sobolev gradients, regularization strategies, admissible-set enforcement,
and carefully designed discretizations within the \emph{optimize--then--discretize} paradigm. These issues are addressed in the
following sections devoted to numerical methods and computational implementation.

\section{Computational Methods}
\label{sec:comp}

This section provides an overview of computational strategies for solving inverse problems with state-dependent constitutive
relations. Building on the variational and adjoint-based foundations developed in Section~\ref{sec:var_adj_grads}, we outline both
classical and modern approaches, emphasizing their connection to the gradient structure \eqref{eq:state_dep_gradient_general} and
the practical challenges identified in Section~\ref{sec:grad_challenge}.

\subsection{Classical Optimization-Based Methods}
\label{sec:class_OPT}

\subsubsection{Two principal paradigms}
\label{sec:two_paradigms}

Early computational approaches to inverse problems with state-dependent parameters are rooted in deterministic, optimization-based
formulations of PDE-constrained inverse problems. Within this framework, two principal paradigms are commonly distinguished, namely
\emph{reduced-space} and \emph{all-at-once} formulations. This distinction is already implicit in the mathematical foundations
developed in Section~\ref{sec:var_adj_grads}, where the reduced cost functional and adjoint-based gradients were introduced.

\textbf{Reduced-space formulations.} In reduced-space methods, the state variables are eliminated using the solution operator
associated with the forward problem. Specifically, for a given constitutive relation $p(\zeta)$, the governing PDE is assumed to be
well-posed and solved to the accuracy required by the optimization algorithm, yielding the state $\bu = \bu(p)$. The inverse problem
is then recast as an unconstrained optimization problem posed solely in the space of constitutive relations, i.e.,
\begin{equation*}
  \min_{p \in \mP} \; \mJ_{\mathrm{red}}(p) \equiv \mJ(\bu(p); p).
\end{equation*}
Gradients of the reduced functional are computed via adjoint equations, as developed in Section~\ref{sec:var_adj_grads}, leading
directly to representations of the form \eqref{eq:state_dep_gradient_general}.

Reduced-space approaches are conceptually appealing and form the basis of many classical algorithms for PDE-constrained optimization.
They allow for a clean separation between the forward solver and the optimization routine, and they are particularly effective when
robust and efficient solvers for the forward and adjoint problems are available. Comprehensive treatments of reduced-space methods can
be found in the monographs by Lions~\cite{Lions1971}, Glowinski \emph{et al.}~\cite{GlowinskiLionsHe2008},
Hinze \emph{et al.}~\cite{Hinze2009}, and Tr\"{o}ltzsch~\cite{Troltzsch2010}. In the context of inverse problems, this paradigm is also
closely related to the ``optimize--then--discretize'' philosophy discussed in~\cite{Biegler2010} and later adopted in~\cite{Bukshtynov2013}.
Large-scale reduced-space implementations of this paradigm, developed in deterministic history-matching and production-optimization settings
for geosciences and energy systems~\cite{Bukshtynov2015,VolkovBukshtynov2018}, further demonstrate the practical viability of adjoint-based
optimize--then--discretize formulations, even though the unknowns in these works are typically spatially distributed rather than state-dependent.
Rigorous analyses of deterministic reduced-space formulations for severely ill-posed PDE-constrained inverse problems, such as Electrical
Impedance Tomography, further demonstrate that adjoint-based optimize--then--discretize methods naturally expose structural non-injectivity
and identifiability limitations imposed by the measurement model. Importantly, these limitations are intrinsic and persist independently of
discretization or numerical error~\cite{AbdullaBukshtynovSeif2021}.

\textbf{All-at-once formulations.} In contrast, all-at-once methods treat the state variables, adjoint variables, and unknown constitutive
relations as coupled unknowns within a single optimization problem. Rather than eliminating the state via the solution operator, the
governing PDEs are enforced as constraints, typically through Lagrange multipliers or penalty terms. This leads to large but highly
structured systems involving the state, adjoint, and parameter variables simultaneously. All-at-once formulations are particularly
attractive when the forward problem is strongly nonlinear, when multiple physics are tightly coupled, or when the solution operator is
difficult to define or differentiate explicitly.

From a numerical perspective, these methods naturally align with Newton-type and interior-point strategies and can exhibit superior
robustness with respect to poor initial guesses. However, they also pose significant challenges in terms of memory requirements, solver
design, and preconditioning. Classical references on all-at-once approaches include the works of
Glowinski \emph{et al.}~\cite{GlowinskiLionsTremolieres1981} and Biegler \emph{et al.}~\cite{Biegler2003conf,Biegler2007}. In inverse
problems and data assimilation, closely related formulations are discussed in~\cite{AschBocquetNodet2016, ItoKunisch2008}.

Both paradigms ultimately rely on the same variational principles, but they differ substantially in how the gradient
structure~\eqref{eq:state_dep_gradient_general} is realized computationally. In reduced-space methods, this structure emerges after
elimination of the state variables and is expressed explicitly through adjoint solutions. In all-at-once methods, the same structure
is embedded implicitly within the coupled optimality system. The choice between these paradigms is therefore problem dependent and
closely tied to the analytical properties of the forward model, the nature of the state dependence, and the computational challenges
outlined in Section~\ref{sec:grad_challenge}.

\subsubsection{Gradient-based optimization algorithms}
\label{sec:gradient_based_OPT}

Once adjoint-based expressions for the gradient of the reduced functional are available, a broad class of gradient-based optimization
algorithms can be brought to bear on the inverse problem. This paradigm constitutes the computational backbone of most deterministic
approaches to PDE-constrained inverse problems and optimal control, including problems with state-dependent parameters.

\textbf{Infinite-dimensional viewpoint.} In contrast to classical parameter estimation problems posed in finite dimensions, inverse
problems with state-dependent constitutive relations are naturally formulated in \emph{function spaces}. As a result, optimization
algorithms are most appropriately interpreted at the \emph{infinite-dimensional} level, with discretization introduced only afterward.
This ``optimize--then--discretize'' perspective ensures that the computed gradients are consistent with the underlying variational
structure and remain meaningful under mesh refinement; see, e.g.,~\cite{Lions1971,Hinze2009,Troltzsch2010,Biegler2010}. Within this
setting, the gradient of the reduced functional $\mJ_{\mathrm{red}}(p) = \mJ(\bu(p);p)$ is an element of the dual space $\mP^*$ and
is typically represented via a Riesz map associated with a chosen inner product on $\mP$. The choice of this inner product plays a
critical role, as it implicitly defines the geometry of the optimization landscape and can significantly affect convergence rates.
Common choices include $L_2$-type inner products, Sobolev metrics, or problem-specific weighted norms;
see~\cite{Hinze2009,BorziSchulz2011,Neuberger1997,ProtasBewleyHagen2004,Protas2008} for detailed discussions.

\textbf{First- and second-order methods.} Given the gradient structure derived in Section~\ref{sec:var_adj_grads}, classical first-order
methods such as steepest descent~(SD) and nonlinear conjugate gradient~(CG) methods can be applied directly in function space. To
accelerate convergence, quasi-Newton methods, most notably BFGS and its limited-memory variants (L-BFGS), are widely used. In the
PDE-constrained setting, these methods construct approximations of the inverse Hessian using gradient information alone and often achieve
a favorable balance between computational cost and convergence speed~\cite{Nocedal2006,Hinze2009}. Full Newton and Gauss--Newton methods
are also employed, especially when the forward model is mildly nonlinear and second-order information can be computed or approximated
efficiently; see~\cite{ItoKunisch2008,Biegler2010}.

\textbf{State dependence and algorithmic implications.} For inverse problems with state-dependent parameters, the evaluation of the reduced
gradient is intrinsically more involved than in the state-independent case, as variations of the constitutive relation induce variations
of the state that must be accounted for consistently. As shown in Section~\ref{sec:var_adj_grads}, this coupling manifests itself through
additional terms in the adjoint equation and in the gradient representation~\eqref{eq:state_dep_gradient_general}. Nevertheless, once the
adjoint problem has been solved, the resulting gradient can be used within standard function-space optimization frameworks without further
modification. The choice of optimization method and metric on $\mP$ is therefore highly problem dependent and must be informed by the
analytical properties of the forward model, the degree of ill-posedness, and the nature of the state dependence discussed in
Section~\ref{sec:grad_challenge}.

\textbf{Generic gradient-based algorithm.}
A prototypical reduced-space optimization algorithm for state-dependent parameter reconstruction is summarized in
Algorithm~\ref{alg:state_dep_PDE_opt}. Its structure mirrors classical PDE-constrained optimization algorithms, cf.~\cite{BukshtynovBook2023},
but explicitly accounts for the state dependence of the unknown constitutive relation through the adjoint-based gradient evaluation. The
algorithm should be interpreted as a generic template rather than a specific implementation. Variants commonly incorporate regularization,
projections onto admissible sets, problem-adapted inner products on $\mP$, and inexact or iterative solves of the state and adjoint equations.
The same algorithmic structure underlies both reduced-space and all-at-once formulations: in the former, the state and adjoint variables
are eliminated through solution operators, whereas in the latter the corresponding optimality conditions are realized implicitly within
a coupled system of state, adjoint, and parameter variables.
\begin{algorithm}[t]
  \caption{Gradient-based optimization for state-dependent parameter reconstruction}
  \label{alg:state_dep_PDE_opt}
  \begin{algorithmic}[1]
    \small
    \STATE Discretize $\Omega$ (and time, if applicable); initialize discrete $\bu$, $p(\zeta)$, and adjoint $\blambda$.
    \STATE Load observational data (measurements) $\tilde{\bu}(\bx)$.
    \STATE Set $k \leftarrow 0$ and choose initial guess $p^{(0)} \in \mP$.
    \WHILE{stopping criterion not satisfied}
      \STATE \textbf{State solve:} Compute $\bu^{(k)} = \bu(p^{(k)})$ by solving the discretized forward
        problem~\eqref{eq:general_PDE}--\eqref{eq:general_BC} (\emph{state eliminated as optimization variable in reduced-space
        formulation}).
      \STATE Evaluate cost functional $\mJ(\bu^{(k)};p^{(k)})$, cf.~\eqref{eq:cost_general}.
      \STATE \textbf{Adjoint solve:} Compute adjoint $\blambda^{(k)}$ associated with $\bu^{(k)}$ and $p^{(k)}$ (\emph{explicit in
        reduced-space; implicit in all-at-once}).
      \STATE \textbf{Gradient evaluation:} Compute a discretized realization of $\bnabla_p \mJ(\bu^{(k)};p^{(k)})$,
        cf.~\eqref{eq:state_dep_gradient_general}.
      \STATE \textbf{Search direction:} Compute $\delta p^{(k)}$ (SD, CG, BFGS/L-BFGS, etc.).
      \STATE \textbf{Update:} Choose step $\alpha^{(k)}$ and set $p^{(k+1)} = p^{(k)} + \alpha^{(k)} \delta p^{(k)}$.
      \STATE $k \leftarrow k + 1$.
    \ENDWHILE
\end{algorithmic}
\end{algorithm}

\subsubsection{Sobolev gradient methods}
\label{sec:Sobolev_gradients}

As demonstrated in~\cite{Bukshtynov2011,Bukshtynov2013}, cost functional gradients obtained naturally in the $L_2$ setting are often
insufficient for numerical optimization in state-dependent inverse problems. In particular, $L_2$ gradients may lack the regularity
required of the unknown parameter, exhibit strong mesh dependence, and be poorly defined outside the identifiability region. These
difficulties can be alleviated by extracting gradients in smoother Sobolev spaces~\cite{Neuberger1997,ProtasBewleyHagen2004,Protas2008},
e.g., $H^1(\mD)$, endowed with the inner product
\begin{equation}
  \label{eq:H1_inner_product_review}
  \langle z_1, z_2 \rangle_{H^1(\mD)} = \int_{\mD} \left( z_1 z_2 + \ell^2 \, \partial_s z_1 \, \partial_s z_2 \right) ds,
  \qquad z_1,z_2 \in H^1(\mD),
\end{equation}
where $\ell \in \RR$ is a user-chosen length-scale parameter. The corresponding Sobolev gradients are obtained by defining appropriate
Riesz maps, i.e.,
\begin{equation}
  \label{eq:riesz_review}
  \mJ'(p;\delta p) = \langle \bnabla_p^{L_2}\mJ, \delta p \rangle_{L_2(\mD)} = \langle \bnabla_p^{H^1}\mJ, \delta p \rangle_{H^1(\mD)},
  \qquad \forall \, \delta p \in H^1(\mD),
\end{equation}
leading to improved regularity, implicit preconditioning, and mesh-independent convergence. The Sobolev $H^1$ gradients $\bnabla_p^{H^1}\mJ$
are computed as the solution of the following elliptic boundary-value problem (Helmholtz-type filter):
\begin{equation}
  \label{eq:sobolev_filter_review}
  \begin{aligned}
    \bnabla_p^{H^1}\mJ - \ell^2 \, \partial_s^2 \bnabla_p^{H^1}\mJ &= \bnabla_p^{L_2}\mJ \qquad && \text{in } \mD, \\
    \partial_s \bnabla_p^{H^1}\mJ &= 0 && \text{on } \partial\mD.
  \end{aligned}
\end{equation}

From an algorithmic standpoint, the use of Sobolev gradients augments the standard adjoint-based optimization loop by a single additional
step: after computing $\bnabla_p^{L_2}\mJ$, one solves~\eqref{eq:sobolev_filter_review} to obtain $\bnabla_p^{H^1}\mJ$, which is then
used to define the search direction. This procedure was employed in~\cite{Bukshtynov2011,Bukshtynov2013} and can be interpreted as a form
of operator-based \emph{preconditioning}. The parameter $\ell$ controls the degree of smoothing applied to the gradient: increasing $\ell$
suppresses high-frequency components, whereas $\ell\to 0$ recovers the $L_2$ gradient. As shown in~\cite{ProtasBewleyHagen2004}, this
operation is equivalent to applying a \emph{low-pass filter} to the $L_2$ gradient, with $\ell$ acting as a cut-off scale. In
reduced-space methods, the Sobolev metric is applied after elimination of the state via a Riesz map, whereas in all-at-once formulations
it is incorporated implicitly through the choice of inner product on the parameter space, without altering the optimality conditions.

\subsubsection{Dealing with noisy measurement data}
\label{sec:noisy_data}

In the presence of noisy measurements, inverse reconstructions of state-dependent constitutive relations are well known to exhibit
instabilities manifested as small-scale oscillations in the recovered parameter~\cite{TikhonovArsenin1977,Neubauer2008,Tarantola2005}.
This behavior was demonstrated explicitly for the class of problems considered here in~\cite{Bukshtynov2011,Bukshtynov2013}, where even
moderate noise levels were shown to contaminate reconstructions obtained without regularization.

Within the optimization framework, a standard stabilization strategy is \emph{Tikhonov regularization}~\cite{EnglHankeNeubauer1996},
whereby the cost functional \eqref{eq:cost_general} is augmented by a penalty term $\mR(p)$ that incorporates a priori information about
the unknown constitutive relation. Two choices were investigated in~\cite{Bukshtynov2011}, corresponding to different regularity
assumptions on $p$.

\textbf{$L_2$ regularization.} Choosing
\begin{equation}
  \mR(p) = \frac{\lambda_1}{2} \| p - \bar p \|_{L_2(\mD)}^2, \qquad \lambda_1 \in \RR^+,
\end{equation}
adds a zero-order penalty discouraging \emph{large local deviations} from a reference constitutive relation $\bar p$. The corresponding
$L_2$ gradient of the regularized cost functional becomes
\begin{equation}
  \bnabla_p^{L_2} \mJ_{\lambda_1} = \bnabla_p^{L_2} \mJ + \lambda_1 (p - \bar p),
\end{equation}
where $\bnabla_p^{L_2} \mJ$ denotes the unregularized gradient derived from the adjoint equations, e.g., obtained in the generalized
form~\eqref{eq:state_dep_gradient_general}. While this approach improves stability, it does not explicitly suppress high-frequency
oscillations.

\textbf{$H^1$ regularization.} Alternatively, imposing smoothness through
\begin{equation}
  \mR(p) = \frac{\lambda_2}{2} \| p - \bar p \|_{\dot H^1(\mD)}^2 = \frac{\lambda_2}{2} \int_{\mD} |\partial_s (p - \bar p)|^2 \, ds,
  \qquad \lambda_2 \in \RR^+
\end{equation}
penalizes \emph{rapid variations} of the constitutive relation. In this case, the $L_2$ gradient acquires an additional second-order term,
\begin{equation}
  \bnabla_p^{L_2} \mJ_{\lambda_2} = \bnabla_p^{L_2} \mJ - \lambda_2 \, \partial_s^2 p,
\end{equation}
supplemented with natural boundary conditions on $\partial\mD$ arising from integration by parts. Numerical experiments
in~\cite{Bukshtynov2011} demonstrated that $H^1$ regularization is significantly more effective than its $L_2$ counterpart in
suppressing noise-induced oscillations.

In practice, these regularized gradients are subsequently mapped into the chosen Sobolev metric using the Riesz map discussed in
Section~\ref{sec:Sobolev_gradients}. The regularization parameters $\lambda_1$ and $\lambda_2$ thus provide an explicit mechanism for
balancing data fidelity against stability, while remaining fully compatible with adjoint-based gradient evaluation; stability and
convergence of Tikhonov regularization using the Sobolev $H^1$ norm were established rigorously in~\cite{Kugler2003}. Algorithmically,
Tikhonov regularization introduces additional terms into the gradient through explicit penalties in the cost functional, while Sobolev
gradients achieve stabilization by applying a metric-dependent smoothing operator to the gradient without modifying the objective.

\subsection{Modern Computational Phase (Emerging 2011--Present)}
\label{sec:modern_phase}

More recent developments in inverse problems with state-dependent constitutive relations are driven by increased model complexity,
higher-dimensional and nonlinear parameterizations, and the availability of large-scale computational resources. Beyond improved
computational capacity, this modern phase is characterized by a broadening of the algorithmic scope, including systematic strategies
for enforcing physical constraints on constitutive relations, techniques for extending the identifiability interval $\mI$ toward a
desired reconstruction range, and more efficient evaluation of gradients in level-set--based formulations. At the same time, attention
has shifted toward time-dependent and multiphysics inverse problems, along with the incorporation of Bayesian inference and
uncertainty-quantification methodologies. Most recently, advances in physics-informed learning, operator-learning frameworks, and
hybrid optimization--data-driven approaches have begun to complement classical adjoint-based techniques, further expanding the range
of computational tools available for state-dependent inverse problems.

\subsubsection{Enforcing constraints on state-dependent constitutive relations}
\label{sec:barrier_slack}

A classical strategy for enforcing constraints in PDE-constrained optimization relies on penalty, barrier, or augmented Lagrangian
formulations. In these approaches, constraints are incorporated into the objective functional through penalty terms or via Lagrange
multipliers combined with penalty stabilization. Augmented Lagrangian methods, in particular, improve conditioning relative to pure
penalty approaches and allow constraints to be enforced accurately without excessively large penalty parameters. These techniques are
well established and widely analyzed in the context of PDE-constrained optimization and control; see,
e.g.,~\cite{GlowinskiLionsHe2008,Bertsekas2016,Nocedal2006,Hinze2009,Biegler2010,Ruszczynski2006,Vogel2002}. Despite their generality,
augmented Lagrangian and barrier methods introduce additional dual variables, algorithmic parameters, and coupling conditions, which can
substantially complicate adjoint-based implementations, particularly when the constitutive relation depends explicitly on the state.

A simpler and computationally economical alternative for enforcing pointwise inequality constraints is based on a direct change of
variables using \emph{slack variables}. This approach, commonly referred to as a slack-variable reparameterization and discussed in,
e.g.,~\cite{BoydVandenberghe2004}, was adopted in~\cite{Bukshtynov2013} to enforce the positivity constraint $\mu(T) > 0$ on
temperature-dependent viscosity $\mu(T)$ in~\eqref{eq:NS_heat_a}. Introducing an unconstrained auxiliary function $\theta(T)$, the
constitutive relation is parametrized as
\begin{equation}
  \mu(T) = \theta^2(T) + m_{\mu}, 
  \label{eq:slack}
\end{equation}
where $m_{\mu}$ corresponds to the prescribed lower bound $m_p$ in the admissible set $\mP$ defined in~\eqref{eq:Pset}. This transformation
automatically enforces $\mu(T)\ge m_{\mu} > 0$ and converts the original inequality-constrained problem~\eqref{eq:opt_general} into an
\emph{unconstrained} one,
\begin{equation}
  \hat{\theta}
  = \underset{\theta \in \mX}{\argmin} \, \mJ(\theta),
  \label{eq:minJslack}
\end{equation}
while preserving the structural form of the cost functional. A key advantage of the slack-variable formulation is its \emph{minimal impact}
on the computational pipeline: the forward and adjoint problems retain their original structure, and the constraint enters only through the
parameterization of the constitutive relation. Consequently, standard reduced-space, adjoint-based optimization algorithms can be applied
without introducing additional multipliers, penalty updates, or auxiliary optimality conditions.

When reconstruction is meaningful only on the identifiability interval $\mI\subseteq\mD$, slack-variable parameterizations provide a natural
mechanism for enforcing physical constraints globally on $\mD$ while allowing the data to inform the solution only on $\mI$. Outside $\mI$,
the reconstructed constitutive relation is determined implicitly by the chosen parametrization and regularization, without requiring explicit
constraint handling or any modification of the adjoint equations.

\subsubsection{Expanding the identifiability region $\mI$}
\label{sec:shifting_I}

In inverse problems with state-dependent constitutive relations, the sensitivity of the cost functional with respect to the unknown parameter
is typically confined to an \emph{identifiability region} $\mI \subset \mD$, determined by the range of state values attained by the forward
solution under a fixed set of sources, boundary conditions, or controls. Outside this interval, gradients may be formally extended, but such
extensions do not generate new sensitivity information and therefore do not support reliable reconstruction.

To address this limitation, an adaptive strategy was introduced in~\cite{Bukshtynov2011} to \emph{expand} the effective identifiability region
$\mI = \mI_{(0)} = [\zeta_{\alpha,(0)},\, \zeta_{\beta,(0)}]$ by modifying the experimental or computational setup so that the state solution
spans shifted intervals
\begin{equation}
  \label{eq:shift_I}
  \mI_{(j)} := [\zeta_{\alpha,(j)},\, \zeta_{\beta,(j)}] = [\zeta_{\alpha,(j-1)} + h_j,\, \zeta_{\beta,(j-1)} + h_j],
  \qquad h_j \in \RR, \quad j \ge 1.
\end{equation}
By performing a sequence of such shifts, the union of identifiability intervals $\mI_{(j)}, j \ge 0$, can be made to cover a larger target
range $\mD$ on which the constitutive relation is to be reconstructed. At a strategic level, the method proceeds as follows:
\begin{enumerate}
  \item[(1)] Modify control inputs, e.g., source terms $f_{(j)}$ and/or boundary conditions $g_{(j)}$
    in~\eqref{eq:general_PDE}--\eqref{eq:general_BC}, to induce a shift of the state solution $\bu_{(j)}$ and hence of the identifiability
    region.
  \item[(2)] Acquire new measurements $\tilde \bu_{(j)}$ corresponding to the modified inputs and solve a PDE-constrained inverse problem
    to reconstruct the constitutive relation $p_{(j)}$ on the shifted interval~\eqref{eq:shift_I}.
  \item[(3)] Iterate this process, initializing the reconstruction $p_{(j)}$ on $\mI_{(j)}$ using the most recent estimate $p_{(j-1)}$, until
    $\bigcup_{j \ge 0} \mI_{(j)} \supset \mD$.
\end{enumerate}

Because each shifted problem generates sensitivity information over a different portion of the state space, later reconstructions must not
degrade estimates obtained on earlier intervals. This is ensured by augmenting the cost functional~\eqref{eq:cost_general} at stage $j$ with
a Tikhonov-type penalty,
\begin{equation}
  \mJ_{(j)}(p) = \mJ(p) + \frac{\gamma}{2} \int_{\zeta_{\alpha,(0)}}^{\zeta_{\beta,(j-1)}-\delta}
  \bigl[p(\zeta) - p_{(j-1)}(\zeta)\bigr]^2 \, d\zeta,
  \label{eq:J_shift}
\end{equation}
where $p_{(j-1)}$ denotes the reconstruction obtained on the union of previous identifiability intervals, and $\gamma,\delta \in \RR^+$ are
the stabilization parameters. This term preserves previously identified parameter values while also providing additional regularization
against measurement noise, cf.~Section~\ref{sec:noisy_data}.

Algorithmically, the shifting strategy decouples the generation of sensitivity information from the solution of the inverse problem itself.
It preserves the adjoint-based optimization framework and can be viewed as an outer-loop mechanism that enriches the inverse problem by
progressively expanding the region of identifiability. This viewpoint naturally extends to time-dependent, multiphysics, and controlled
systems, where the accessible state range can be influenced through external inputs.

From a broader perspective, identifiability-region shifting can be viewed as a form of adaptive experimental or computational design, in
which control inputs are chosen to steer the state trajectory through previously unexplored regions of the state space. In this sense, the
procedure operates at the level of the control-to-observation mapping, actively shaping the domain on which the parameter-to-observation
operator is informative.

\subsubsection{Integration over codimension-one level sets}
\label{sec:level_set_int}

When the reduced gradient $\bnabla_p \mJ$ of the cost functional $\mJ(p)$ with respect to a state-dependent constitutive relation
$p = p(\zeta)$ can be expressed in a level-set form, cf.~\eqref{eq:state_dep_gradient_general}, its numerical evaluation requires
integration over codimension-one manifolds defined implicitly by the state variable $\zeta$. While this structure connects naturally to
classical level-set and shape-calculus ideas, it also introduces computational challenges, particularly when the forward and adjoint
solutions are available only in discrete form.

At a fixed value $\zeta_0$ of the state variable $\zeta$, the gradient $\bnabla_p \mJ(\zeta_0)$ involves integration over the level set
\begin{equation}
  \Gamma_{\zeta_0} := \{ \bx \in \Omega \; : \; \zeta(\bx) = \zeta_0 \},
  \label{eq:level_set_def_review}
\end{equation}
whose geometry is not known \emph{a priori} and must be inferred from grid-based approximations of $\zeta$. As a result, the numerical task
is not simply quadrature, but the accurate evaluation of integrals over implicitly defined, possibly complex contours (or hypersurfaces in
higher dimensions), whose location and topology depend on the solution of the forward problem. More generally, integrals over
$\Gamma_{\zeta_0}$ can be expressed using the coarea formula. As discussed at length in~\cite{Bukshtynov2013}, for sufficiently smooth
$\phi(\zeta_0,\bx) : \RR \times \Omega \rightarrow \RR$ representing the level-set function whose zero level defines the contour of
integration $\Gamma_{\zeta_0}$ and satisfying $|\nabla \phi| \neq 0$ on $\Gamma_{\zeta_0}$, one has
\begin{equation}
  \int_{\Omega} \delta\bigl(\phi(\zeta_0,\bx)\bigr)\, g(\bx)\, d\bx \; = \;
  \int_{\Gamma_{\zeta_0}} \frac{g(\bx)}{|\nabla \phi(\zeta_0,\bx)|}\, d\sigma,
  \label{eq:coarea_review}
\end{equation}
where $\phi(\zeta_0,\bx) = \zeta(\bx) - \zeta_0$  and $g(\bx) : \Omega \rightarrow \RR$ represents
$F\!\left( \bu,\blambda,\bx,t;\, p(\zeta) \right)$ from~\eqref{eq:state_dep_gradient_general} in the present context (as in~\cite{Bukshtynov2013},
the time dependence is omitted here). This identity reveals that gradient evaluation can be interpreted equivalently as either (i)~a contour
(line) integral over $\Gamma_{\zeta_0}$ or (ii)~an area integral over $\Omega$ involving a Dirac delta distribution. Both interpretations are
mathematically equivalent, but they lead to different numerical strategies and error mechanisms.

Existing methods for evaluating integrals of both forms in~\eqref{eq:coarea_review} can be broadly classified into two categories:
\begin{itemize}
  \item[(A)] \emph{Contour-based (geometric) methods}, which explicitly approximate the level set $\Gamma_{\zeta_0}$ and perform numerical
    integration directly on the reconstructed manifold;
  \item[(B)] \emph{Area-based methods}, which approximate the Dirac delta distribution $\delta(\phi)$ and evaluate the corresponding volume
    integral over $\Omega$.
\end{itemize}
This distinction mirrors classical developments in level-set methods and has been extensively studied in the numerical analysis literature
\cite{Engquist2004,ZahediTornberg2010,Mayo1984,Smereka2006,Beale2008,MinGibou2007,MinGibou2008,Towers2007,Towers2009}.

Within group~(A), geometric integration techniques developed by Min and Gibou~\cite{MinGibou2007,MinGibou2008} approximate level sets by
piecewise polynomial curves (in two dimensions) obtained from a triangulation of $\Omega$. The contour integral is then decomposed into
a collection of standard one-dimensional quadratures. These methods can achieve high accuracy but require nontrivial geometric reconstruction
and data structures. Group~(B) encompasses a variety of area-integration techniques. One class relies on regularized Dirac delta functions
$\delta_\varepsilon$, whose support is tied to the grid resolution~\cite{Engquist2004,ZahediTornberg2010}. Improved accuracy can be obtained
by adapting the regularization width to the local gradient magnitude $|\nabla \phi|$. A related family of methods, initiated by
Mayo~\cite{Mayo1984} and further developed by Smereka~\cite{Smereka2006}, constructs discrete approximations of the delta distribution with
compact support. Consistent finite-difference--based approximations were also proposed by Towers~\cite{Towers2007,Towers2009}. Area-based
methods are generally simpler to implement and robust on fixed grids, but their accuracy depends sensitively on the chosen regularization
and grid resolution.

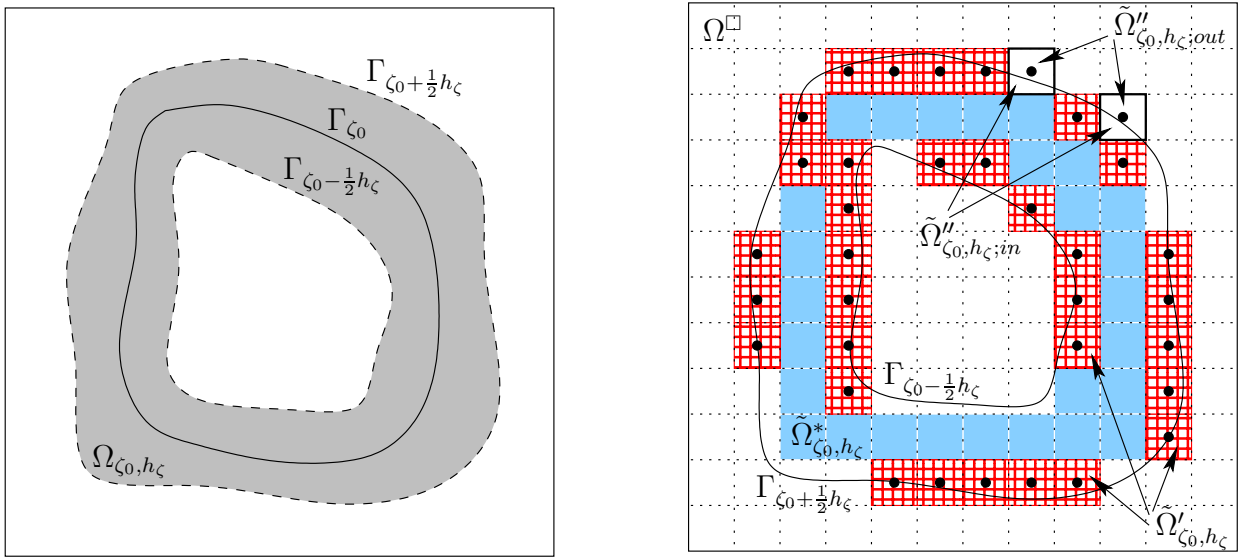
\begin{figure}[h]
  \centering
  \begin{tikzpicture}
    \node[inner sep=0] (img) at (0,0)
      {\includegraphics[width=0.45\textwidth]{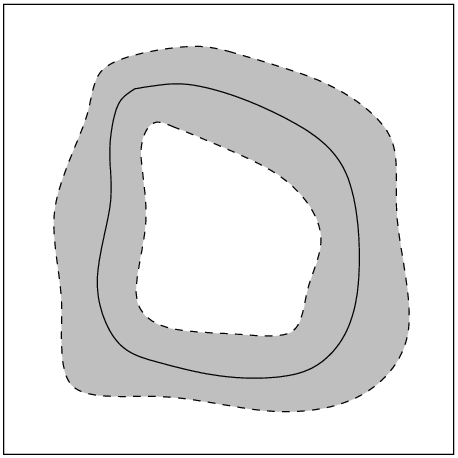}};
    \node at (0.9,2.1) {$\Gamma_{\zeta_0}$};
    \node at (0.7,1.4) {$\Gamma_{\zeta_0-\frac12 h_\zeta}$};
    \node at (1.8,2.7) {$\Gamma_{\zeta_0+\frac12 h_\zeta}$};
    \node at (-2,-2.4) {$\Omega_{\zeta_0,h_\zeta}$};
  \end{tikzpicture}
  \hfill
  \begin{tikzpicture}
    \node[inner sep=0] (img) at (0,0)
      {\includegraphics[width=0.45\textwidth]{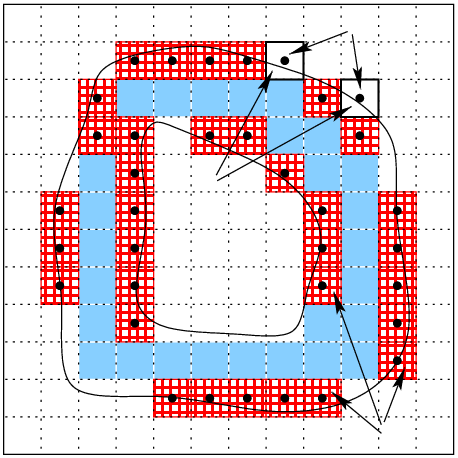}};
    \node at (-3.2,3.3) {$\Omega^{\square}$};
    \node at (-2.1,-2.9) {$\Gamma_{\zeta_0+\frac12 h_\zeta}$};
    \node at (-0.4,-1.4) {$\Gamma_{\zeta_0-\frac12 h_\zeta}$};
    \node at (-1.8,-2.15) {$\tilde{\Omega}^{*}_{\zeta_0,h_\zeta}$};
    \node at (0.1,0.5) {$\tilde{\Omega}_{\zeta_0,h_\zeta;in}''$};
    \node at (2.75,3.3) {$\tilde{\Omega}_{\zeta_0,h_\zeta;out}''$};
    \node at (3.05,-3.35) {$\tilde{\Omega}_{\zeta_0,h_\zeta}'$};
  \end{tikzpicture}
  \caption{Illustration of the hybrid area-contour integration strategy (adopted from~\cite{Bukshtynov2013}). (a)~The narrow band
    $\Omega_{\zeta_0,h_\zeta} = \{\bx \in \Omega : \zeta(\bx) \in [\zeta_0 - \tfrac12 h_\zeta,\, \zeta_0 + \tfrac12 h_\zeta]\}$
    between the neighboring level sets $\Gamma_{\zeta_0-\frac12 h_\zeta}$ and $\Gamma_{\zeta_0+\frac12 h_\zeta}$. (b)~Finite-element
    approximation of the narrow band used for area-based quadrature. The computational mesh $\Omega^{\square}$ consists of quadrilateral
    elements; shaded cells represent elements lying entirely within $\Omega_{\zeta_0,h_\zeta}$, forming $\tilde{\Omega}^{*}_{\zeta_0,h_\zeta}$,
    while checked cells indicate elements whose representative points lie inside $\Omega_{\zeta_0,h_\zeta}$ and form
    $\tilde{\Omega}'_{\zeta_0,h_\zeta}$. These subsets are used to construct the hybrid area-contour approximation of integrals over the level
    set $\Gamma_{\zeta_0}$. See Section~\ref{sec:level_set_int} and~\cite{Bukshtynov2013} for details.}
  \label{fig:hybrid_area_contour}
\end{figure}

\textbf{Hybrid ``area-contour'' approach~\cite{Bukshtynov2013}.} In the context of state-dependent inverse problems, an additional difficulty
arises because gradient information is required over an \emph{entire interval} of state values, rather than at a single level set.
\cite{Bukshtynov2013} introduced a hybrid strategy that combines features of both groups~(A) and~(B). The method approximates contour integrals
by carefully constructed area integrals that exploit local information about the level-set geometry while avoiding explicit geometric
reconstruction. This approach was designed specifically for PDE-constrained inverse problems with state-dependent parameters and provides
a compromise between accuracy, robustness, and computational complexity. The geometric idea underlying this construction is illustrated in
Figure~\ref{fig:hybrid_area_contour}, which shows the narrow band between neighboring level sets and its approximation by a collection of finite
elements used for area-based quadrature.

Let the left-hand side of~\eqref{eq:coarea_review} be denoted by $f(\zeta_0)$; the procedure consists of the following steps:
\begin{itemize}
  \item[(1)] \emph{State averaging.} For a fixed state value $\zeta_0 \in \mI$, introduce a small state increment $h_\zeta > 0$ and replace
    pointwise evaluation of $f(\zeta_0)$ by averaging over the interval $[\zeta_0 - \tfrac{1}{2}h_\zeta,\, \zeta_0 + \tfrac{1}{2}h_\zeta]$,
    i.e.,
    \begin{equation*}
      f(\zeta_0) \;\approx\; \frac{1}{h_\zeta} \int_{\zeta_0 - \frac{1}{2}h_\zeta}^{\zeta_0 + \frac{1}{2}h_\zeta} f(\zeta)\, d\zeta.
    \end{equation*}
  \item[(2)] \emph{Subdomain construction.} Define the narrow band
    \begin{equation*}
      \Omega_{\zeta_0,h_\zeta} \;:=\; \left\{ \bx \in \Omega \;:\; \zeta(\bx) \in \left[\zeta_0 - \tfrac{1}{2}h_\zeta,\,
      \zeta_0 + \tfrac{1}{2}h_\zeta \right] \right\},
    \end{equation*}
    containing all points between the two neighboring level sets. This region is approximated by the union of $N = N_{\zeta_0,h_\zeta}$ finite
    elements whose representative (center) points $\bx_j^0$ satisfy $\zeta(\bx_j^0) \in [\zeta_0 - \tfrac{1}{2}h_\zeta,\, \zeta_0 +
    \tfrac{1}{2}h_\zeta], j = 1, \ldots, N$ (see Figure~\ref{fig:hybrid_area_contour}).
  \item[(3)] \emph{Area-based quadrature.} The contour integral value $f(\zeta_0)$ is approximated by an area integral over
    $\Omega_{\zeta_0,h_\zeta}$ and evaluated using a compound midpoint rule,
    \begin{equation}
      f(\zeta_0) \;\approx\; \frac{1}{h_\zeta} \int_{\zeta_0 - \frac{1}{2}h_\zeta}^{\zeta_0 + \frac{1}{2}h_\zeta}
      \int_{\Omega} \delta\bigl(\zeta(\bx)-\zeta\bigr)\, g(\bx)\, d\bx\, d\zeta \;\approx\; \frac{h^2}{h_\zeta} \sum_{j=1}^{N} g(\bx_j^0).
      \label{eq:hybrid_approx}
    \end{equation}
    Here, $h$ denotes a characteristic spatial mesh size, $\zeta_0$ denotes a fixed level of the state variable, and $\zeta$ is used as
    the integration variable over the state interval.
\end{itemize}
The accuracy of this hybrid approximation is quantified by the following result.

\begin{theorem}[Accuracy of ``area-contour'' approach {\cite{Bukshtynov2013}}]
  \label{thm:lsm_int_acc}
  The hybrid approximation~\eqref{eq:hybrid_approx} satisfies
  \begin{equation*}
    \int_{\Omega} \delta\!\left(\zeta(\bx)-\zeta_0\right) g(\bx)\, d\bx \;=\; \frac{h^2}{h_\zeta}
    \sum_{j=1}^{N_{\zeta_0,h_\zeta}} g(\bx_j^0) \;+\; \mO(h^2) \;+\; \mO(h_\zeta),
  \end{equation*}
  i.e., it is second-order accurate with respect to the spatial discretization in $\bx$ and first-order accurate with respect to
  the discretization of the state variable.
\end{theorem}

Overall, integration over codimension-one level sets constitutes a central numerical challenge in adjoint-based inversion with state-dependent
constitutive relations. The methods outlined above provide a spectrum of trade-offs between geometric fidelity and algorithmic simplicity,
and their relative merits depend strongly on the problem dimension, grid structure, and required reconstruction accuracy.

\subsubsection{Verification and reproducibility}
\label{sec:verif_kappa}

As discussed in Sections~\ref{sec:barrier_slack}--\ref{sec:level_set_int}, recent developments in PDE-constrained optimization for inverse problems
reflect a gradual shift from purely analytical gradient derivations toward integrated computational frameworks that combine functional analysis,
numerical optimization, and modern data-driven techniques. In such settings, gradients of the cost functional with respect to state-dependent
parameters, denoted earlier by $\bnabla_p \mJ(p)$ and given in abstract form by~\eqref{eq:state_dep_gradient_general}, are often obtained
through complex adjoint-based procedures involving discretized forward and backward PDE solves. As a consequence, systematic verification of
computed gradients is an essential prerequisite before initiating any optimization or inversion process, both to ensure correctness of the
implementation and to promote reproducibility of computational results.

A widely adopted diagnostic tool for gradient verification in PDE-constrained optimization is the so-called \emph{$\kappa$-test}, which compares
directional derivatives of the cost functional computed in two independent ways: (i)~via finite-difference approximations of the cost
functional~$\mJ(p)$ itself, and (ii)~via inner products involving the adjoint-based gradient~$\bnabla_p \mJ$. This test has been extensively
employed in the context of state-dependent and time-dependent inverse problems~\cite{Bukshtynov2011,Bukshtynov2013,Homescu2002}, and is discussed
in a broader computational optimization context in~\cite{BukshtynovBook2023,ChunEdwardsBukshtynov2024}. Given a parameter $p$ and a perturbation
direction $\delta p$, the $\kappa$-test examines the quantity
\begin{equation}
  \kappa(\epsilon) := \dfrac{\epsilon^{-1}\bigl[\mJ(p + \epsilon \, \delta p)-\mJ(p)\bigr]}
  {\int_{\mD} \bnabla_p^{L_2}\mJ(s)\, \delta p(s)\,ds},
  \label{eq:kappa_test_review}
\end{equation}
evaluated over a range of perturbation magnitudes $\epsilon$. The numerator represents a finite-difference approximation of the directional
derivative of $\mJ$, while the denominator corresponds to the same quantity computed using the adjoint-based gradient through the appropriate
Riesz representation.

It is important to emphasize that, while the choice of inner product on the parameter space (e.g., $L_2$ versus Sobolev $H^1$) plays a crucial
algorithmic role in defining well-behaved and mesh-independent search directions, as discussed in Section~\ref{sec:Sobolev_gradients}, the
$\kappa$-test itself probes the consistency of the underlying directional derivative $\mJ'(p;\delta p)$ and is therefore metric-invariant at the
continuous level, provided the gradient and perturbation are represented consistently.

At the same time, excessively large values of the Sobolev smoothing parameter $\ell$ in~\eqref{eq:sobolev_filter_review} may lead to over-smoothing
of the gradient, which can artificially degrade the observed accuracy of the $\kappa$-test despite a formally correct adjoint implementation, as
noted in~\cite{AbdullaBukshtynovHagverdiyev2019}. In addition, transformations applied to the control variables (such as reduced-order or PCA-based
parameterizations introduced to decrease the effective dimensionality of the optimization problem) act directly on the gradient representation and
may likewise degrade the apparent quality of the $\kappa$-test, even when the underlying adjoint formulation is correct, as discussed
in~\cite{ChunEdwardsBukshtynov2024}.

If the adjoint-based gradient is implemented correctly and the discretization accurately approximates the underlying continuous problem, the
ratio $\kappa(\epsilon)$ should remain close to unity over an intermediate range of $\epsilon$. In practice, this plateau typically spans several
orders of magnitude and provides strong evidence of gradient correctness. Deviations of $\kappa(\epsilon)$ from unity for very small values of
$\epsilon$ are expected due to subtractive cancellation and round-off errors, whereas deviations for large $\epsilon$ arise from truncation
errors associated with the finite-difference approximation. Moreover, systematic refinement of the spatial or temporal discretization generally
extends the range over which $\kappa(\epsilon)\approx 1$, reflecting improved consistency between discrete and continuous gradients, particularly
within the ``optimize--then--discretize'' paradigm~\cite{Bukshtynov2011,Bukshtynov2013,BukshtynovBook2023,ArbicBukshtynov2024}.

Beyond serving as a correctness check, the $\kappa$-test also offers quantitative insight into the effective accuracy of gradient evaluations.
For example, plotting $\log_{10}|\kappa(\epsilon)-1|$ as a function of $\epsilon$ reveals the number of significant digits captured by the computed
gradient, providing a practical diagnostic for assessing the reliability of subsequent optimization iterations. As emphasized
in~\cite{BukshtynovBook2023}, incorporating such verification procedures as a standard component of computational workflows is essential for
ensuring robustness, transparency, and reproducibility in large-scale PDE-constrained inverse problems, particularly in modern settings involving
state-dependent parameters, multiphysics coupling, and hybrid physics--data-driven methodologies.

\subsubsection{Time-dependent and multiphysics inverse problems}
\label{sec:time_multiphysics}

While much of the early computational work in PDE-constrained inverse problems focused on steady (time-independent)
formulations~\cite{Bukshtynov2011,Egger2014,Grisvard1985}, many physically relevant applications involve \emph{time evolution} of state variables
and \emph{coupled multiphysics phenomena}. In such settings, the unknown constitutive relations or parameters may depend both on the evolving state
and on temporal dynamics, introducing significant analytical and computational challenges. It is important to distinguish between parameters that
depend explicitly on time and those that depend implicitly on time through the evolving state variables. In the latter case, state-dependent
constitutive relations embedded within time-dependent PDEs lead to additional nonlinearity and coupling in both the forward and adjoint systems,
complicating identifiability, sensitivity analysis, and numerical optimization.

A prototypical class of time-dependent inverse problems aims to reconstruct \emph{time-varying coefficients} or source terms in evolution equations,
such as parabolic, hyperbolic, or coupled systems. For instance, simultaneous reconstruction of time-dependent coefficients in parabolic heat equations
and related diffusion models has been studied using over-specified boundary or integral conditions, demonstrating uniqueness and providing algorithmic
strategies for combined time-space reconstruction~\cite{Huntul2021a,HuntulTekin2023,IbraheemHuntul2025}. Similarly, inverse source problems for
time-dependent heat equations with non-classical or nonlocal boundary conditions have been analyzed using integral equation methods and regularization
techniques~\cite{Hazanee2015,Hazanee2019,Kanca2013}. For hyperbolic dynamics, recovery of time-dependent potentials and displacement distributions
under periodic conditions has been pursued with tailored finite-difference and Tikhonov stabilization schemes~\cite{Huntul2021b}, as well as from
the viewpoint of stability theory~\cite{Krishnan2021}. These studies highlight the additional challenges posed by wave propagation, finite signal
speed, and the need for robust reconstruction strategies.

From an analytical perspective, time-dependent inverse problems are frequently only conditionally stable and require careful treatment of solution
regularity, particularly when coefficients depend on the evolving state. Techniques based on Carleman estimates, compactness arguments, and stability
inequalities have therefore played a central role in establishing uniqueness and stability for evolution equations, including problems with
time-dependent or state-dependent coefficients~\cite{Yamamoto2009}. From a PDE-constrained optimization viewpoint, many time-dependent inverse problems
are formulated as optimization over temporally evolving state variables and parameters. These formulations typically require the solution of forward
and adjoint \emph{evolution equations}, with adjoint variables integrated backward in time to compute sensitivities efficiently. Classical control
literature provides foundational treatments of time-dependent optimality systems and associated numerical methods in both continuous and discrete time
settings. For example, optimal control of parabolic and hyperbolic systems using adjoint-based gradient evaluation has been examined
extensively~\cite{TroltzschYousept2012}. More recently, space-time variational frameworks have been developed to improve well-posedness and computational
efficiency of time-dependent PDE-constrained optimization problems~\cite{BeranekReinholdUrban2023}.

To address the resulting computational challenges, algorithmic innovations such as domain decomposition in time and low-rank temporal
preconditioning have been developed. These techniques enable scalable solution of large forward-adjoint systems arising in time-dependent
PDE-constrained optimization~\cite{BarkerStoll2015,StollBreiten2015}. In parallel, iterative regularization methods tailored to the temporal structure
of parameter identification have been explored. Landweber--Kaczmarz schemes compare all-at-once and reduced formulations for identifying
time-dependent parameters, demonstrating trade-offs in convergence behavior and regularization efficacy over temporal intervals~\cite{Nguyen2019}.
Ensemble Kalman approaches and derivative-free ensemble methods provide alternative frameworks for spatiotemporal inverse problems, blending
regularization with data assimilation strategies~\cite{Iglesias2013,Iglesias2016}.

\emph{Multiphysics} inverse problems demand careful consideration not only of the dynamics of the governing PDEs but also of the intricate interactions
across different physical fields, as well as the significant computational costs associated with solving tightly coupled systems. In these settings,
unknown parameters often depend on multiple state variables that evolve simultaneously, giving rise to complex, highly nonlinear optimization
landscapes with expensive forward and adjoint solves. Consequently, robust parameter identification requires a synergistic integration of rigorous
theoretical analysis, algorithmic regularization, and scalable numerical optimization strategies. Representative examples of such systems include
(1)~thermo-fluid systems with state-dependent material properties~\cite{Bukshtynov2013}, (2)~parameter estimation in coupled fluid-structure
interaction models~\cite{Sunseri2020}, (3)~inverse modeling of thermo-hydro-mechanical processes in porous media~\cite{Amini2022},
(4)~multiphase poromechanics~\cite{Amini2023}, and (5)~geological and subsurface multiphysics inverse problems where multimodal data are exploited
for joint reconstructions~\cite{Yin2023}.

In thermo-fluid systems, optimization-based formulations have been developed to reconstruct state-dependent material properties, such as
temperature-dependent viscosity, within coupled momentum and energy equations. A notable example casts the inverse problem as a PDE-constrained
optimization problem, where gradients of the cost functional are efficiently computed via adjoint analysis and specialized integration along level
sets of the state fields~\cite{Bukshtynov2013}. In the context of fluid-structure interaction, joint state and parameter estimation approaches have
been employed to infer fluid or structural parameters from sparse measurements, highlighting the intricate interplay between fluid forces and
structural deformation within a fully coupled optimization framework~\cite{Sunseri2020}.

Beyond classical adjoint-based PDE-constrained optimization, recent years have witnessed the emergence of data-driven and hybrid approaches aimed at
mitigating the computational burden of tightly coupled multiphysics inversions. Sequential physics-informed neural networks~(PINNs) have been
proposed to extract material and transport parameters governing thermo-hydro-mechanical processes, leveraging the coupled PDE structure within the
optimization loop~\cite{Amini2022}. Similarly, PINN-based inversion strategies have been applied to nonisothermal multiphase poromechanics, enabling
the recovery of constitutive parameters in highly nonlinear coupled systems~\cite{Amini2023}. Learned surrogates and constraint-augmented inversion
frameworks have also been utilized for geological carbon storage monitoring, demonstrating that hybrid approaches combining surrogate models,
physics-informed constraints, and PDE solves can significantly enhance reconstruction accuracy while mitigating computational cost~\cite{Yin2023}.
Despite their flexibility, such approaches often face challenges related to training stability, interpretability of learned parameters, and the
enforcement of rigorous physical constraints, particularly in regimes with limited data or stiff multiphysics coupling.

Further insights into several of these emerging directions, including physics-informed learning, operator-learning frameworks, Bayesian and
uncertainty-quantification approaches, and hybrid optimization--data-driven methods, are discussed in the following section~\ref{sec:non_traditional_OPT}.
While not intended as a comprehensive review, this discussion highlights representative ideas that transcend traditional deterministic optimization
paradigms. It should also be noted that not all works cited in this section explicitly address state-dependent parameter reconstruction; many focus
on parameters that vary in time or space. Nonetheless, the approaches surveyed (particularly PDE-constrained optimization with adjoint-based sensitivities
and hybrid physics-informed strategies) are naturally suited for state-dependent parameter identification in time-dependent and multiphysics systems,
and their systematic development therefore represents a compelling direction for future research in high-fidelity multiphysics modeling.

\subsubsection{Transcending traditional deterministic optimization}
\label{sec:non_traditional_OPT}

The increasing complexity of time-dependent and multiphysics inverse problems, together with limited and noisy observational data, has motivated
approaches that go beyond classical deterministic PDE-constrained optimization. In particular, modern formulations increasingly emphasize
uncertainty quantification, data-driven representation of solution operators or constitutive relations, and hybrid strategies that combine
\emph{physics-based} optimization with \emph{learning-based} components. While these approaches do not replace adjoint-based methods, they offer
complementary perspectives that are especially valuable for state-dependent parameter reconstruction in highly nonlinear, high-dimensional, and
computationally expensive multiphysics systems.

\textbf{Bayesian and uncertainty-quantification approaches.}
Bayesian inverse problem formulations provide a probabilistic framework for parameter identification by combining observational data with prior
information and explicitly quantifying uncertainty in the reconstructed parameters. In the context of time-dependent and multiphysics PDEs,
Bayesian approaches have been developed using Gaussian process priors, hierarchical models, and ensemble-based sampling or approximation techniques,
enabling posterior inference over spatiotemporal parameters~\cite{Stuart2010,DashtiStuart2017}. For large-scale evolution equations, computationally
tractable approximations such as ensemble Kalman methods and their variants have been widely adopted, particularly in data assimilation and
parameter estimation for dynamical systems~\cite{Evensen2009,Iglesias2013,Iglesias2016}. Recent work has extended Bayesian formulations to
state-dependent and nonlinear inverse problems, including PDEs with heterogeneous and state-dependent coefficients, emphasizing scalable approximate
inference and uncertainty quantification for complex multiphysics models~\cite{BarajasTartakovsky2019,Lan2023}. A representative example is the work
of Zhao et al.~\cite{Zhao2021}, which combines adjoint-based PDE-constrained optimization with Bayesian inference to reconstruct multiple state-dependent
constitutive laws from spatiotemporal observations in coupled reaction-diffusion systems. A closely related example in the electrochemical context is
provided by~\cite{Sethurajan2019}, which combines adjoint-based PDE-constrained optimization with Bayesian inference to quantify uncertainty in
reconstructed \emph{state-dependent transport laws}. Building on deterministic function-valued reconstructions of electrolyte transport properties,
the Bayesian formulation yields posterior distributions and confidence bands for constitutive relations such as concentration-dependent diffusivities
and transference numbers. This work illustrates how uncertainty quantification can be layered on top of classical variational optimization without
abandoning adjoint-based structure, and directly complements the deterministic electrochemical applications discussed in
Section~\ref{sec:apps_echem_2015_2026}, where the focus is on reconstruction rather than uncertainty assessment. These methods provide principled
uncertainty estimates but often face challenges related to computational cost, prior specification, and scalability in tightly coupled multiphysics
settings. For broader and pedagogical perspectives on Bayesian inverse problems and uncertainty quantification in physics-based models, including
probabilistic formulation, posterior analysis, and computational considerations, we refer the reader to recent review and survey
articles~\cite{Bingham2024,Reiser2025,Bai2024}.

\textbf{Physics-informed and operator-learning methods.}
Physics-informed learning approaches, including physics-informed neural networks, aim to incorporate governing PDEs directly into learning
objectives, enabling simultaneous inference of state variables, parameters, and constitutive relationships from sparse or noisy
data~\cite{Raissi2019,Raissi2020,Tartakovsky2020}. More recently, operator-learning frameworks such as DeepONets and neural operator architectures
have been introduced to approximate solution operators mapping inputs (e.g., coefficients, sources, or boundary data) to PDE solutions, offering
mesh-independent and data-efficient representations for parametric and time-dependent PDEs~\cite{Lu2021,Kovachki2023}. These methods have been applied
to inverse problems and parameter identification in time-dependent and multiphysics systems, where they can capture complex nonlinear dependencies
between evolving states and constitutive relations. Learning-informed inversion strategies that explicitly combine data-driven components with PDE
structure have also been proposed for nonlinear, time-dependent parameter identification, illustrating how physics-informed learning can augment
classical inversion frameworks while preserving interpretability and physical consistency~\cite{Aarset2023}. While physics-informed and operator-learning
methods offer flexibility and reduced reliance on adjoint solvers, they raise important questions regarding training stability, generalization outside
the training regime, and the enforcement of physical constraints, particularly for stiff or strongly coupled PDE systems. For comprehensive overviews
of physics-informed learning and operator-learning frameworks in scientific computing, including theoretical foundations, algorithmic developments,
and applications to inverse problems, the reader is referred to recent survey articles~\cite{Karniadakis2021,Luo2025,Kovachki2024}.

\textbf{Hybrid optimization--data-driven methods.}
Hybrid approaches seek to combine the rigor and interpretability of PDE-constrained optimization with the flexibility of data-driven models. Typical
strategies include the use of learned surrogates or reduced-order models within optimization loops, machine-learning--assisted regularization,
and differentiable surrogate models embedded in adjoint-based or gradient-based frameworks~\cite{Xu2022,Sun2023,MowlaviNabi2023}. In large-scale
multiphysics inverse problems, such hybrid methods have been shown to significantly reduce computational cost while preserving physical consistency,
for example in subsurface flow, seismic monitoring, and geological carbon storage applications~\cite{Yin2023,Louboutin2023}. Recent contributions
also emphasize adaptive and ``one-shot'' surrogate construction strategies, in which surrogate models are learned online or jointly with the
optimization variables to accelerate repeated PDE solves and reduce simulation burden~\cite{Keil2022,Guth2024}. A closely related but conceptually
distinct example is provided by~\cite{Akerson2025}, which employs neural-network parameterizations \emph{within} a classical adjoint-based
PDE-constrained optimization framework to represent high-dimensional constitutive relations. In contrast to physics-informed or operator-learning
approaches, the neural networks in~\cite{Akerson2025} are not trained in a data-driven sense; rather, they serve as flexible function spaces for
state-dependent coefficients, with all network weights optimized deterministically through reduced-space gradient-based methods. The governing PDE
remains enforced exactly, gradients are computed via adjoint equations, and the optimization objective retains a standard output least-squares form.
This formulation preserves the mathematical structure, interpretability, and convergence properties of classical PDE-constrained optimization while
enabling richer representations of constitutive laws than traditional low-dimensional parametric ansatzes. As such, it exemplifies a principled
extension of deterministic optimization rather than a departure from it. Overall, hybrid strategies illustrate how data-driven components can be
leveraged to accelerate forward solves or enrich constitutive representations, while careful coupling with physics-based constraints remains
essential to ensure robustness, interpretability, and reliability of reconstructed state-dependent parameters.

Taken together, the works discussed in this section illustrate how contemporary inverse problem methodologies increasingly operate at the
interface between physics-based modeling and data-driven inference. Although many of the cited contributions are not formulated explicitly in terms
of state-dependent parameter reconstruction, they introduce probabilistic, operator-level, and surrogate-based perspectives that fundamentally
reshape how inverse problems in time-dependent and multiphysics PDE systems can be posed and solved. In this sense, Bayesian uncertainty-aware
formulations, physics-informed and operator-learning techniques, and hybrid optimization--data-driven strategies should be viewed not as isolated
alternatives to classical PDE-constrained optimization, but as complementary components of a broader inversion paradigm capable of addressing the
nonlinearity, scale, and complexity inherent in state-dependent multiphysics models.

\section{Applications and Case Studies}
\label{sec:apps}

The mathematical and computational developments surveyed in the preceding sections have enabled a transition from proof-of-concept studies to
\emph{end-to-end computational frameworks} capable of addressing realistic inverse problems with state-dependent constitutive relations. Since
roughly 2011, advances in adjoint-based PDE-constrained optimization, scalable numerical solvers, regularization strategies, and high-performance
computing have made it possible to deploy these methods in practical settings arising in physics, engineering, and materials science. In contrast
to early studies focused on simplified model problems or low-dimensional parameterizations, contemporary applications increasingly involve
time-dependent, multiphysics systems, large-scale discretizations, and nonparametric representations of constitutive laws, often integrated within
robust optimization pipelines.

The purpose of this section is to survey \emph{representative applications} and \emph{case studies} in which state-dependent parameter identification
has been demonstrated as a viable and effective \emph{computational tool}. Rather than providing a tutorial treatment, we emphasize complete or
near-complete computational workflows, highlighting how theoretical concepts (identifiability intervals, adjoint-based gradients, Sobolev
regularization, admissible set enforcement, etc.) are realized in practice. Each case study is organized around three guiding elements: (i)~a concise
description of the underlying physical or engineering problem and the associated inverse formulation; (ii)~the principal computational components
employed, including discretization strategies, adjoint solvers, regularization mechanisms, and optimization algorithms; and (iii)~the main achievements,
limitations, and lessons learned, with an emphasis on reproducibility, scalability, and robustness.

Collectively, the applications discussed below illustrate both the \emph{maturity} and the \emph{diversity} of modern state-dependent inverse problem
methodologies. They span stationary and time-dependent regimes, single-physics and coupled multiphysics models, and deterministic as well as hybrid
data-driven formulations. Beyond documenting successful deployments, this section aims to guide readers toward existing solution strategies that can
be adapted to related problems, and to identify gaps where current frameworks remain incomplete or where new methodological developments are needed.
In this sense, the case studies serve not only as validation of the theory and algorithms reviewed earlier, but also as a roadmap for future research
and application-driven innovation in PDE-constrained optimization with state-dependent parameters.

\subsection{Heat Conduction with Temperature-Dependent Conductivity}
\label{sec:apps_heat_2011}

We begin the applications section with the stationary inverse heat conduction problem with \emph{temperature-dependent thermal conductivity}, which
serves as a canonical and computationally complete example of state-dependent parameter identification. This problem was studied in detail
in~\cite{Bukshtynov2011} and represents one of the earliest end-to-end PDE-constrained optimization frameworks for reconstructing a constitutive law
of the form $k = k(T)$ without restricting it to a finite-dimensional parametric ansatz. Because of its analytical transparency and clear physical
interpretation, this example plays a foundational role in demonstrating how the abstract concepts developed in Sections~\ref{sec:math} and~\ref{sec:comp}
translate into a practical, verifiable, and reproducible computational methodology.

\textbf{Problem setting and inverse formulation.}
The forward model is the stationary heat conduction equation~\eqref{eq:heat}, posed on a bounded domain $\Omega \subset \RR^d$, $d=1,2,3$, with prescribed
boundary temperature $T_0$ and known volumetric heat source $f$. The unknown quantity is the thermal conductivity $k(T)$, assumed to depend solely on
the temperature field $T$ and to satisfy uniform positivity bounds ensuring ellipticity and thermodynamic admissibility. Physically, this setting models
heat transport in materials whose conductive properties vary with temperature, a situation commonly encountered in high-temperature processes, liquid
metals, and engineered materials.

The inverse problem consists in reconstructing the constitutive relation $k(T)$ from temperature measurements $\tilde T(\bx)$ available on a sensing
subdomain $\Sigma \subseteq \Omega$. As discussed in Sections~\ref{sec:problem_formulation} and~\ref{sec:var_adj_framework}, this problem differs
fundamentally from identifying a spatially varying coefficient $k(\bx)$: the unknown is a \emph{function of the state} $T$, and only values of $k(T)$
corresponding to temperatures actually attained by the solution are identifiable. The reconstruction is therefore meaningful only on the associated
\emph{identifiability interval} $\mI = [T_\alpha,T_\beta]$, while outside this interval the inverse problem is underdetermined and must rely on prior
information or regularization. The inverse problem is formulated as a reduced-space PDE-constrained optimization problem of the form~\eqref{eq:final_opt},
with a least-squares misfit functional and optional Tikhonov regularization acting directly on the constitutive law over a prescribed interval
$\mD \supset \mI$, cf.~\eqref{eq:intI}--\eqref{eq:intD}.

\textbf{Computational framework.}
The computational approach adopted in~\cite{Bukshtynov2011} follows the \emph{optimize--then--discretize} paradigm described in
Section~\ref{sec:var_adj_framework} and emphasized throughout this review. A key modeling device enabling a computationally tractable adjoint formulation
is the \emph{Kirchhoff transformation} (see, e.g., classical treatments of nonlinear diffusion in~\cite{Kirchhoff1894}), which introduces an auxiliary
variable $V(T)=\int_{T_\alpha}^T k(s)\,ds$ and converts the nonlinear diffusion operator into a linear Poisson problem in the transformed variable.
This transformation plays a central role in deriving adjoint equations and gradient expressions that avoid explicit Dirac delta distributions and
are robust at the discrete level.

The forward and adjoint problems are discretized using standard numerical methods (finite differences in the one-dimensional test case considered
in~\cite{Bukshtynov2011}), while the constitutive law $k(T)$ is represented independently on the temperature interval using smooth interpolation.
This separation between the discretization of the physical domain and that of the constitutive relation is essential: it allows refinement of
the PDE solver without introducing artificial structure into the reconstruction of $k(T)$ and underscores the method's independence from any
particular spatial discretization technology.

Gradients of the reduced cost functional are computed using a continuous adjoint formulation derived in Section~\ref{sec:var_adj_grads}. As shown
in Theorem~\ref{thm:Bukshtynov2011_gradient}, the resulting gradient admits a level-set-based representation that aggregates sensitivity information
over spatial regions where the temperature field attains a given value. For numerical purposes, this representation is recast into equivalent volume
and boundary integral expressions involving the adjoint solution, yielding a gradient formula that is straightforward to implement and evaluate.

To ensure stability and mesh-independent convergence, \emph{Sobolev gradients} on the temperature interval are employed, as discussed in
Section~\ref{sec:Sobolev_gradients}. Raw $L_2$ gradients (which are generally discontinuous and confined to the identifiability interval) are
filtered through an $H^1(\mD)$ Riesz map, producing smooth descent directions defined on the entire reconstruction interval. This Sobolev
\emph{preconditioning} acts simultaneously as a regularization mechanism and as gradient smoothing controlled by an intrinsic length-scale parameter.
Additional stabilization is achieved through explicit Tikhonov regularization, with both $L_2$- and $H^1$-type penalties examined to control noise
amplification.

The optimization problem is solved using gradient-based algorithms in function space, including steepest descent and quasi-Newton methods,
cf.~Sections~\ref{sec:opt_problem} and \ref{sec:gradient_based_OPT}. Each optimization iteration requires one forward and one adjoint PDE solve,
resulting in a computational cost comparable to that of classical PDE-constrained optimization problems with spatially distributed controls.

\textbf{Achievements and demonstrated capabilities.}
The numerical experiments reported in~\cite{Bukshtynov2011} demonstrate that the proposed framework can accurately reconstruct nontrivial
temperature-dependent conductivity laws from partial and noisy temperature data. In particular, the study illustrates: (i)~the decisive role
of the identifiability interval in determining which portions of $k(T)$ can be recovered from a given experimental configuration;
(ii)~the effectiveness of Sobolev gradients and $H^1$-based regularization in producing stable, physically meaningful reconstructions; and
(iii)~the necessity of adjoint-based gradient evaluation for computational efficiency and scalability.

A notable feature of the work is the systematic verification of adjoint-based gradients using $\kappa$-tests (cf.~Section~\ref{sec:verif_kappa}),
carried out over many orders of magnitude in the perturbation size and for multiple discretization resolutions. This level of verification provides
strong evidence for the correctness of the continuous adjoint derivation and for the consistency of the discrete implementation. The computational
results further demonstrate that the framework can be extended beyond a single identifiability interval by actively modifying boundary conditions
and source terms, thereby shifting the temperature range explored by the system. This procedure effectively casts the inverse problem as
an \emph{adaptive experimental design} strategy, in which system inputs are tuned to generate sensitivity information over previously inaccessible
portions of the state space.

\textbf{Limitations and lessons learned.}
At the same time, the study highlights several intrinsic limitations characteristic of state-dependent inverse problems. First, reconstruction
is fundamentally restricted to the identifiability interval associated with a given experiment; outside this interval, sensitivity information is
absent and the inverse problem is underdetermined. Second, even in this comparatively simple elliptic setting, the optimization problem is nonconvex,
and numerical experiments demonstrate a tangible dependence of the reconstruction on the choice of the initial guess, indicating convergence to
different local minimizers. Third, the problem is severely ill-posed with respect to measurement noise: increasing noise levels lead to oscillatory
instabilities in the reconstructed constitutive laws unless regularization is employed.

The results show that Sobolev gradients provide a degree of implicit regularization, while explicit Tikhonov regularization is essential for robust
performance at moderate to high noise levels, with $H^1$-type penalties offering improved control of oscillatory modes. These findings underscore
the importance of combining appropriate function-space metrics, regularization strategies, and experimental design considerations.

From a broader perspective, this application serves as a \emph{reference case study} for state-dependent parameter identification. It establishes
a complete and reproducible computational pipeline (from variational formulation and adjoint analysis to discretization, gradient verification,
regularization, and optimization) that can be systematically extended to more complex models. Many of the algorithmic ideas introduced here,
including Kirchhoff transformations, Sobolev-gradient preconditioning, identifiability-aware reconstruction, and adaptive experimental design,
reappear in more advanced applications discussed in the subsequent subsections, where additional challenges arise from time dependence, coupling,
and increased model complexity.

\subsection{Temperature-Dependent Viscosity in Coupled Thermo-Fluid Flow}
\label{sec:apps_NS_2013}

Here, we summarize the second application, developed in~\cite{Bukshtynov2013}, which extends the state-dependent identification paradigm from
Section~\ref{sec:apps_heat_2011} to a genuinely \emph{multiphysics} setting. The goal is to reconstruct a constitutive law of the form $\mu=\mu(T)$
appearing in the momentum balance~\eqref{eq:NS_heat_a}--\eqref{eq:NS_heat_b}, where the temperature field $T$ is governed by a \emph{separate}
energy equation~\eqref{eq:NS_heat_c}. In contrast to the previous heat-conduction example, the defining computational feature here is that the
reduced gradient involves \emph{integrals over level sets of the temperature field}, cf.~\eqref{eq:B13_grad_L2}, which necessitates specialized
numerical quadrature strategies.

\textbf{Problem setting and inverse formulation.}
The forward model consists of the incompressible Navier--Stokes equations coupled with an energy equation for heat transfer (cf.~the general
multiphysics PDE-constrained formulation~\eqref{eq:NS_heat} in Section~\ref{sec:gov_eqn}). The unknown control is the temperature-dependent viscosity
$\mu(T)$ used in the momentum equation, while the state comprises the velocity-pressure pair and the temperature. The inverse problem is posed as
a reduced-space optimization problem of the type~\eqref{eq:final_opt}, with a least-squares misfit between predicted and observed temperatures at a
finite set of sensors and over a time window, as done in~\cite{Bukshtynov2013}. As in Section~\ref{sec:anal_properties}, identifiability is
inherently \emph{state-dependent}: only values of $\mu(T)$ corresponding to temperatures attained by the solution are informed by the data, so the
reconstruction is effectively meaningful on an identifiability interval (and associated measurement span) determined by the experiment.

\textbf{Computational framework.}
The framework follows the \emph{optimize--then--discretize} viewpoint adopted throughout the review, cf.~Section~\ref{sec:var_adj_framework}.
A continuous adjoint system is derived for the coupled thermo-fluid PDE, and the reduced gradient is expressed in the level-set form described in
Section~\ref{sec:var_adj_grads}, cf.~Theorem~\ref{thm:Bukshtynov2013_gradient}: the sensitivity of the objective with respect to $\mu(\cdot)$
at a temperature value $T$ aggregates contributions from spatial regions where the state attains that temperature. In~\cite{Bukshtynov2013},
this yields gradient expressions involving integrals over temperature level sets, i.e., over manifolds defined implicitly by $T(\bx,t)=\mathrm{const}$,
which is structurally different from standard distributed-control problems and requires dedicated numerical evaluation.

A central algorithmic component is therefore the numerical approximation of these level-set integrals. The paper systematically compares three
representative strategies: (i)~geometric line integration over approximate level sets, (ii)~regularized Dirac-delta/area-integration approaches,
and (iii)~hybrid contour-to-area formulations, cf.~Section~\ref{sec:level_set_int}. Importantly, since the control lives on a temperature interval,
the paper highlights the interplay between discretization in physical space (mesh size $h$) and discretization of the temperature interval (step size
$h_T$), and how this interaction affects overall reconstruction accuracy, cf.~Theorem~\ref{thm:lsm_int_acc}.

Regularization and robustness mechanisms follow the tools introduced earlier. First, Sobolev-gradient smoothing on the temperature interval
(Section~\ref{sec:Sobolev_gradients}) is used to obtain stable descent directions for $\mu(\cdot)$. Second, explicit Tikhonov regularization is added
to mitigate noise amplification, with the regularization strength acting as a tunable robustness parameter. The resulting infinite-dimensional
optimization problem is solved by gradient-based methods, including nonlinear conjugate-gradient and BFGS-type updates with line search, as in
Sections~\ref{sec:opt_problem} and~\ref{sec:gradient_based_OPT}, where each iteration requires one forward unsteady thermo-fluid solve and one
backward-in-time adjoint solve.

\textbf{Achievements and demonstrated capabilities.}
The study~\cite{Bukshtynov2013} demonstrates an end-to-end reconstruction pipeline for $\mu=\mu(T)$ in a coupled unsteady thermo-fluid system
(a two-dimensional lid-driven cavity with heat transfer). The main technical achievement is a practically usable adjoint-based gradient whose
level-set structure is respected at the discrete level via level-set quadrature (integration on isotherms), together with a quantitative comparison
of three numerical gradient-evaluation strategies. In addition, the paper explicitly shows that the \emph{recoverable} portion of the constitutive
law is governed by the temperature range explored by the experiment, and it demonstrates that one can \emph{shift/enlarge} the identifiability
region by adjusting temperature boundary conditions, thereby obtaining accurate reconstructions on different target intervals (an experimental-design
interpretation consistent with Section~\ref{sec:anal_properties}).

Robustness to measurement noise is assessed in a controlled manner: without Tikhonov regularization, increasing noise levels induce oscillatory
instabilities in the reconstructed constitutive law, whereas increasing the regularization parameter suppresses these instabilities and yields
smoother reconstructions. The experiments also show that noise can perturb the optimization landscape, leading to convergence to different local
minimizers while still preserving the qualitative features of the true constitutive law.

\textbf{Limitations and lessons learned.}
The work underscores several limitations characteristic of state-dependent inverse problems in coupled unsteady PDEs. First, identifiability remains
experiment-limited: outside the temperature range generated by the system, the inverse problem is underdetermined. Second, the reduced functional is
nonconvex and may admit multiple local minimizers, an effect exacerbated by noise; consequently, initialization and regularization choices can
materially influence the final reconstruction. Third, the computational burden is substantially higher than in elliptic single-physics examples,
since each optimization iteration requires solving both a forward unsteady thermo-fluid problem and a backward adjoint problem, and the gradient
evaluation requires additional geometric and numerical care due to level-set integrals. While the adjoint formulation provides algorithmic scalability
in the sense of one adjoint solve per iteration independent of the discretized control dimension, the study does not address parallel performance;
extending the approach to large-scale three-dimensional configurations would require additional HPC-oriented solver and preconditioning strategies.

From a methodological standpoint, the main lesson is that in multiphysics state-dependent identification problems, \emph{gradient representation}
and \emph{its numerical realization} can become the dominant issue: reliable reconstruction hinges not only on adjoint correctness
(Section~\ref{sec:verif_kappa}) and regularization (Section~\ref{sec:noisy_data}), but also on accurately and efficiently computing level-set
sensitivities, and on aligning discretizations in physical and state spaces. These insights motivate subsequent applications where coupling,
time dependence, and state-space geometry play an even more pronounced role.

\subsection{State-Dependent Vorticity Reconstruction in Inviscid Flows}
\label{sec:apps_vorticity_2015}

We next consider an application drawn from inviscid fluid mechanics, developed in~\cite{DanailaProtas2015}, which provides a particularly clean and
transparent illustration of state-dependent parameter identification within a classical deterministic PDE-constrained optimization framework. In
contrast to the diffusion- and viscosity-driven examples discussed in Sections~\ref{sec:apps_heat_2011} and~\ref{sec:apps_NS_2013}, the present
setting involves the reconstruction of a nonlinear \emph{vorticity-streamfunction constitutive law} governing steady vortex dynamics. Despite the
apparent simplicity of the governing equations, this problem exhibits many of the structural features emphasized throughout this review, including
infinite-dimensional controls, adjoint-based gradient evaluation, Sobolev-type preconditioning, and identifiability governed by the range of the
state variable.

\textbf{Problem setting and inverse formulation.}
The forward model is based on the steady, incompressible, inviscid Euler equations in axisymmetric form without swirl. Introducing the Stokes
streamfunction $\psi$ in cylindrical coordinates $(z,r,\theta)$, the velocity field is given by
\begin{equation*}
  \bv = \frac{1}{r}\left( \partial_r \psi,\,-\partial_z \psi , 0 \right),
\end{equation*}
and the azimuthal vorticity $\omega = \partial_z \bv_r - \partial_r \bv_z$ satisfies
\begin{equation}
  \Delta^* \psi = - \omega, \qquad
  \Delta^* \psi := \partial_z\!\left( \frac{1}{r}\,\partial_z \psi \right) +\partial_r\!\left( \frac{1}{r}\,\partial_r \psi \right),
  \label{eq:danaila_stream}
\end{equation}
cf.~\cite{Batchelor1988,Saffman1992}. In steady inviscid flow, the vorticity is transported along streamlines, implying that $\omega$ can be expressed
as a function of the streamfunction (a rescaling commonly used in axisymmetric formulations), $\omega = r \, f(\psi)$, where $f$ is an \emph{a priori
unknown constitutive function}. Substituting this relation into the vorticity-streamfunction formulation~\eqref{eq:danaila_stream} yields the nonlinear
elliptic boundary value problem
\begin{equation}
  \Delta^* \psi = - r\,f(\psi) \qquad \text{in } \Omega,
  \label{eq:danaila_forward}
\end{equation}
supplemented with appropriate boundary conditions on $\partial\Omega$.

The inverse problem considered in~\cite{DanailaProtas2015} is to reconstruct the vorticity function $f(\psi)$ from measurements of the flow, typically
available in the form of velocity or streamfunction data on (parts of) the domain boundary, serving as the observation operator in the cost functional.
As in the previous applications, the unknown parameter is not a spatial field but a \emph{function of the state variable} $\psi$ itself. Consequently,
reconstruction is meaningful only on the identifiability interval $\mI = [0,\psi_{\max}]$, where $\psi_{\max}$ denotes the maximum streamfunction value
attained by the solution. Following the general framework outlined in Sections~\ref{sec:problem_formulation} and~\ref{sec:var_adj_framework}, the inverse
problem is formulated as a reduced-space PDE-constrained optimization problem in which a least-squares misfit functional measuring discrepancy between
predicted and observed data is minimized with respect to $f(\cdot)$, subject to the nonlinear constraint~\eqref{eq:danaila_forward}. Regularization is
introduced by penalizing the roughness of $f$ in a suitable Sobolev space defined on $\mI$.

\textbf{Computational framework.}
The computational approach in~\cite{DanailaProtas2015} follows the \emph{optimize--then--discretize} paradigm emphasized throughout this review.
A continuous adjoint problem is derived for the nonlinear elliptic constraint~\eqref{eq:danaila_forward}, leading to an explicit expression for the
gradient of the reduced cost functional with respect to the control $f(\psi)$. As in the heat conduction and thermo-fluid examples, the resulting
gradient exhibits a \emph{level-set structure}, cf.~\eqref{eq:state_dep_gradient_general}: sensitivity with respect to the value of $f$ at a given
$\psi = \psi_0$ accumulates contributions from spatial locations where the streamfunction attains that value. This structure can be interpreted
through coarea-type arguments and provides a direct analogue of the temperature-level-set gradients discussed in Section~\ref{sec:var_adj_grads}.

For numerical implementation, the forward and adjoint problems are discretized using finite-element methods, while the control $f(\psi)$ is represented
independently on the interval $\mI$. To obtain smooth and mesh-independent descent directions, the raw $L_2(\mI)$ gradient is preconditioned using
an $H^1(\mI)$ Riesz map, yielding a Sobolev gradient with an intrinsic smoothing length scale, cf.~Section~\ref{sec:Sobolev_gradients}. This
Sobolev-gradient framework serves simultaneously as a regularization mechanism and as a means of improving convergence behavior. The resulting
infinite-dimensional optimization problem is solved using classical gradient-based algorithms, including steepest descent and nonlinear conjugate-gradient
methods with line search. As in the earlier applications, each optimization iteration requires one forward and one adjoint PDE solve, resulting in
a computational cost comparable to that of standard PDE-constrained optimization problems.

\textbf{Achievements and demonstrated capabilities.}
The numerical experiments reported in~\cite{DanailaProtas2015} demonstrate that the proposed framework can accurately reconstruct nontrivial
vorticity-streamfunction relations from limited flow measurements, thereby recovering the underlying constitutive law governing vortex structure
rather than merely fitting a flow field. In particular, the study shows that: (i)~the identifiability of $f(\psi)$ is strictly limited to the range
of $\psi$ explored by the flow, reinforcing the central role of the identifiability interval; (ii)~Sobolev-gradient preconditioning is essential
for suppressing spurious oscillations and for obtaining physically meaningful reconstructions; and (iii)~adjoint-based gradient evaluation enables
efficient treatment of an infinite-dimensional control without resorting to parametric approximations.

The work further illustrates the robustness of the method through reconstructions of both analytically defined vortices (such as Hill-type vortex
solutions) and numerically generated vortex-ring data. By relying solely on boundary measurements, the study highlights the ability of state-dependent
inversion frameworks to extract interior constitutive information from indirect and incomplete observations, a feature of practical relevance in
experimental fluid dynamics. The methodological relevance of this approach is further underscored by its subsequent discussion within a broader modeling
framework for steady inviscid vortex rings, where similar reconstruction ideas are presented from a modeling and interpretation perspective,
see Section~2.5 of~\cite{Danaila2021}.

\textbf{Limitations and lessons learned.}
At the same time, the study highlights limitations characteristic of state-dependent inverse problems governed by nonlinear elliptic equations. First,
reconstruction remains confined to the identifiability interval associated with the realized flow configuration; outside this interval, the inverse
problem is underdetermined and regularization dominates. Second, the reduced optimization problem is nonconvex, and numerical experiments demonstrate
sensitivity to initialization, with different initial guesses potentially leading to distinct local minimizers. Third, while the inviscid and steady
setting provides a clean testbed for methodological development, extensions to unsteady or viscous flows would introduce additional analytical and
computational challenges, including time-dependent adjoints and increased ill-posedness. In addition, the reconstructed constitutive law is inherently
tied to the validity of the steady inviscid model, so modeling errors (e.g., neglected viscosity or unsteadiness) may be indistinguishable from
inversion error in practical applications.

From a broader perspective, this application provides a conceptually important complement to the diffusion- and viscosity-based examples discussed
earlier. It demonstrates that the classical PDE-constrained optimization framework for state-dependent parameter identification applies equally to
transport-dominated systems with implicit constitutive structure, and that adjoint-based, Sobolev-regularized optimization can recover nonlinear
state-dependent laws even when the governing physics differs substantially from parabolic or viscous models.

\subsection{State-Dependent Closure Identification in Reduced-Order Models}
\label{sec:apps_ROM_closure_2015}

We next consider an application focused on the identification of \emph{state-dependent closure models} in reduced-order representations of fluid
flows, developed in~\cite{ProtasNoack2014,ProtasNoack2015}. This line of work extends the classical PDE-constrained optimization framework for
state-dependent parameter identification discussed in the preceding full-field applications to the level of \emph{reduced-order models} (ROMs),
where unresolved dynamics must be represented through effective constitutive relations. In contrast to the heat conduction, thermo-fluid, and
vorticity-based examples (cf.~Sections~\ref{sec:apps_heat_2011}--\ref{sec:apps_vorticity_2015}), the present setting operates on dynamical systems
obtained by Galerkin projection, while retaining the defining features of deterministic variational optimization, adjoint-based gradient evaluation,
and nonparametric state-dependent reconstruction.

\textbf{Problem setting and inverse formulation.}
The forward model is a Galerkin reduced-order approximation of the incompressible Navier--Stokes equations, obtained by Galerkin projection of the
governing equations onto a finite set of spatial basis functions (e.g., proper orthogonal decomposition modes), cf.~\cite{ProtasNoack2015}. Denoting
by $\ba(t) \in \RR^N$ the vector of modal amplitudes, the ROM dynamics take the generic form
\begin{equation}
  \dot{\ba}(t) = \bF(\ba(t)) + \bC(\ba(t)),
  \label{eq:rom_general}
\end{equation}
where $\bF(\ba)$ represents the resolved (projected) dynamics and $\bC(\ba)$ is a closure term accounting for the effect of unresolved scales,
see~\cite{ProtasNoack2014,ProtasNoack2015}. The Galerkin projection framework underlying such reduced-order models is classical; for comprehensive
treatments of POD--Galerkin modeling in fluid mechanics, see, e.g.,~\cite{Holmes2012,Noack2011}.

A central modeling choice, introduced in~\cite{ProtasNoack2014} and further developed in~\cite{ProtasNoack2015}, is to restrict the closure to
depend on a scalar global state quantity characterizing the system, most notably the resolved fluctuation energy
\begin{equation}
  E(t) = \tfrac{1}{2}\,\|\ba(t)\|_2^2,
  \label{eq:rom_energy}
\end{equation}
leading to closure models in~\eqref{eq:rom_general} of the form
\begin{equation}
  \bC(\ba) = -\,\nu_T(E)\,\bL\,\ba,
  \label{eq:eddy_visc}
\end{equation}
where $\nu_T(E)$ is an unknown \emph{state-dependent eddy-viscosity function} and $\bL$ is a known linear operator (negative semidefinite for
the boundary conditions considered)\footnote{The representation~\eqref{eq:eddy_visc} makes explicit that the state-dependent eddy-viscosity
$\nu_T(E)\ge 0$ contributes additional dissipation.}, representing the linear viscous contribution and arising from the Galerkin projection,
cf.~\cite{ProtasNoack2015}. The inverse problem consists in reconstructing the function $\nu_T(\cdot)$ from time-resolved measurements of the ROM
state~$\ba(t)$, typically obtained by projecting high-fidelity simulation or experimental data onto the reduced basis, see~\cite{ProtasNoack2014}
for a minimal validation setting and~\cite{ProtasNoack2015} for realistic flow configurations. The reconstruction is carried out under physical
admissibility constraints, most notably non-negativity of the eddy-viscosity function, ensuring consistency with energy dissipation and the
underlying physical interpretation of the closure~\cite{ProtasNoack2015}.

As in the preceding applications, the unknown is not a finite set of parameters but an \emph{infinite-dimensional function of the state}.
Consequently, identifiability is restricted to the interval of energies explored by the observed trajectories, $\mI = [E_{\min},E_{\max}]$,
and the reconstruction problem is formulated as a reduced-space optimization problem of the form~\eqref{eq:final_opt}, with a least-squares misfit
functional measuring discrepancy between observed and predicted modal dynamics, cf.~\cite{ProtasNoack2014,ProtasNoack2015}. Regularization is
introduced by penalizing the roughness of $\nu_T(E)$ in an appropriate Sobolev space defined on $\mI$.

\textbf{Computational framework.}
The computational strategy adopted in~\cite{ProtasNoack2014,ProtasNoack2015} is firmly rooted in the classical variational framework for deterministic
inverse problems and adheres to the \emph{optimize--then--discretize} philosophy outlined earlier in this review. Starting from the continuous-time
reduced-order dynamical system \eqref{eq:rom_general}--\eqref{eq:eddy_visc}, adjoint equations are derived at the level of the ROM, yielding analytic
expressions for the sensitivity of the reduced cost functional with respect to variations in the closure function $\nu_T(E)$.

A distinctive feature of this setting is that the control depends on a scalar state variable \emph{evolving in time}. As a consequence, the gradient
of the reduced functional with respect to $\nu_T(E)$ admits a \emph{trajectory-based state aggregation} structure: the sensitivity at a given energy
level $E=E_0$ is obtained by accumulating contributions from all time instances at which the ROM trajectory attains that energy. This time-aggregation
mechanism plays a role analogous to the level-set integrals encountered in the PDE-based applications, but here it is realized along trajectories in
state space rather than over spatial manifolds.

For numerical realization, the forward and adjoint ROM equations are discretized in time, while the closure function $\nu_T(E)$ is discretized
independently on the admissible energy interval $\mI$. This separation of discretizations mirrors the strategy employed in the full-PDE examples and
allows the resolution of the reduced dynamics to be refined without imposing artificial structure on the reconstructed closure. The resulting gradients
are used within standard gradient-based optimization algorithms, enabling efficient solution of the inverse problem at a computational cost that scales
with the ROM dimension rather than with the underlying high-fidelity model.

\textbf{Achievements and demonstrated capabilities.}
The numerical studies in~\cite{ProtasNoack2014} demonstrate the feasibility of nonparametric state-dependent closure identification in a minimal
three-state model describing the transient approach to a stable limit cycle, thereby providing a controlled validation of the methodology. Building
on this foundation,~\cite{ProtasNoack2015} applies the framework to realistic POD--Galerkin ROMs of shear flows and bluff-body wakes, with reduced
dimensions ranging from tens to hundreds of modes. In these settings, the reconstructed state-dependent eddy-viscosity functions significantly
improve ROM stability and predictive accuracy, enabling autonomous reduced-order simulations that faithfully reproduce transient dynamics and
long-time statistics. In particular, the identified closures act as stabilization mechanisms for POD--Galerkin models that are otherwise prone
to spurious energy growth and loss of robustness, a central motivation emphasized in~\cite{ProtasNoack2015}.

A key outcome of this work is the demonstration that closure models can be reconstructed directly as functions of physically meaningful state
quantities, rather than prescribed a priori through ad hoc parametric forms. The results further illustrate that adjoint-based optimization provides
a scalable and systematic route to closure identification, even in moderate-dimensional dynamical systems, and that Sobolev-regularized reconstruction
is essential for suppressing spurious oscillations induced by noise or limited data.

\textbf{Limitations and lessons learned.}
The studies ~\cite{ProtasNoack2014,ProtasNoack2015} also highlight limitations intrinsic to state-dependent closure identification in reduced-order
settings. First, identifiability is fundamentally restricted to the range of energies explored by the training trajectories; outside this range, the
closure function is underdetermined and dominated by regularization. Second, the resulting optimization problems are nonconvex, and numerical
experiments demonstrate sensitivity to initialization and to the choice of regularization parameters, particularly in higher-dimensional ROMs. Third,
the reconstructed closures are inherently model-dependent: inaccuracies in the reduced basis or neglected physical mechanisms may be absorbed into
the inferred closure function, making it difficult to disentangle modeling error from closure error.

From a broader perspective, this application demonstrates that the classical PDE-constrained optimization paradigm for state-dependent parameter
identification extends naturally from full-field PDE models to reduced-order dynamical systems. It provides a systematic bridge between the full-PDE
inverse problems discussed earlier and data-driven closure modeling, showing that adjoint-based, Sobolev-regularized optimization can recover
interpretable state-dependent closures while preserving the determinism, transparency, and computational efficiency characteristic of classical
variational approaches.

\subsection{State-Dependent Eddy-Viscosity Reconstruction and Fundamental Limits in Turbulent Flows}
\label{sec:apps_MP_2020_2024}

We next consider a sequence of applications developed in~\cite{MatharuProtas2020,MatharuProtas2022,MatharuProtas2024}, which extend the classical
PDE-constrained optimization framework for state-dependent parameter identification to the context of turbulent flows and large-eddy simulation~(LES)
modeling. Building on the closure-identification ideas introduced earlier for reduced-order models (Section~\ref{sec:apps_ROM_closure_2015}), this
line of work addresses the reconstruction of \emph{state-dependent eddy-viscosity closures} at the level of \emph{full PDEs}, and, importantly,
investigates the fundamental limitations of such reconstructions. In contrast to the preceding applications, which primarily demonstrate feasibility,
the present studies systematically examine questions of \emph{well-posedness}, \emph{optimality}, and \emph{nonexistence} in state-dependent inverse
problems.

\textbf{Problem setting and inverse formulation.}
The forward model is based on the incompressible Navier--Stokes equations (or their large-eddy simulation (LES) filtered counterparts), in which the
effect of unresolved scales is modeled through an eddy-viscosity closure. In its simplest form, the governing equations read
\begin{equation}
  \label{eq:MP_NS}
  \begin{aligned}
    \partial_t \bv + (\bv \cdot \nabla)\bv + \nabla p &= \nabla \cdot \bigl[ (\nu + \nu_T(q)) \nabla \bv \bigr] + \boldsymbol{f}, \\
    \nabla \cdot \bv &= 0,
  \end{aligned}
\end{equation}
where $\bv(\bx,t)$ denotes the velocity field, $p$ is the pressure, $\nu$ is the molecular viscosity, $\nu_T(q)$ is an \emph{eddy-viscosity
closure}\footnote{The subscript $T$ emphasizes the turbulent (subgrid-scale) contribution to momentum transport.}, and $\boldsymbol{f}$
represents a prescribed body-force or external driving term used to sustain statistically stationary flow regimes.

The scalar variable $q$ characterizes the local or global flow state and serves as the argument of the unknown constitutive function $\nu_T(\cdot)$.
Its precise definition depends on the modeling context and the specific inverse problem under consideration. For example, in~\cite{MatharuProtas2020}
the state variable is chosen as a scalar measure of local velocity gradients in a simplified one-dimensional setting, while in~\cite{MatharuProtas2022}
and~\cite{MatharuProtas2024} it represents filtered kinetic-energy--like quantities arising naturally in LES formulations of two-dimensional turbulent
flows. In all cases, $q$ is constructed from the resolved velocity field and is intended to capture the intensity of the local or global flow activity.

The inverse problem consists in reconstructing the function $\nu_T(q)$ from observational data, typically in the form of instantaneous flow fields or
target statistical quantities obtained from high-fidelity simulations. As in the earlier applications, the control variable is not a finite set of
parameters but an \emph{infinite-dimensional function of the state}. Consequently, identifiability is restricted to the range of state values explored
by the flow, and the optimization problem is formulated in reduced space as in~\eqref{eq:final_opt}, with a least-squares objective measuring discrepancy
between model predictions and reference data, subject to the nonlinear PDE constraint~\eqref{eq:MP_NS}.

\textbf{Computational framework.}
All three studies adopt a classical \emph{optimize--then--discretize} approach and derive continuous adjoint equations associated with the
state-dependent closure. The resulting gradients of the reduced cost functional with respect to $\nu_T(q)$ exhibit the familiar \emph{state-aggregation
structure} discussed in earlier sections: sensitivity with respect to the value of $\nu_T$ at a given state $q=q_0$ accumulates contributions from
space-time regions where the flow attains that state.

In~\cite{MatharuProtas2020}, this framework is first examined in a simplified one-dimensional setting, serving as a controlled environment for
validating adjoint-based gradient computation and Sobolev-regularized functional reconstruction. The analysis is then extended
in~\cite{MatharuProtas2022} to two-dimensional turbulent flows and LES-type closures, where the cost functional targets statistical quantities rather
than pointwise flow fields. Typical choices include long-time--averaged or distributional observables (e.g., spectral or structure-function--type
summaries and related turbulence statistics), which makes the inverse problem inherently ``operator/statistics matching'' rather than snapshot
matching, cf.~\cite{MatharuProtas2022,MatharuProtas2024}. Finally,~\cite{MatharuProtas2024} generalizes the formulation by incorporating explicit
state constraints on $\nu_T(q)$, leading to optimization problems posed on constrained manifolds in function space.

Across all cases, Sobolev regularization plays a central role in stabilizing the inversion and enforcing smoothness of the reconstructed closure.
In particular, \cite{MatharuProtas2022} advocates an $H^2$ Sobolev metric on the state interval for the closure function, introducing an $H^2$ inner
product with tunable length-scale parameters and computing the corresponding Sobolev gradient via the associated Riesz map. This higher-order~($H^2$)
preconditioning provides enhanced suppression of spurious oscillations compared with $H^1$ smoothing and proves useful when the targeted statistics
induce rough $L_2$ sensitivities. To solve the resulting infinite-dimensional problems, gradient-based optimization algorithms are employed
with one forward and one adjoint PDE solve required for each optimization iteration. This results in computational costs comparable to those of other
large-scale PDE-constrained optimization problems discussed, e.g., in~Section~\ref{sec:apps_NS_2013}.

\textbf{Achievements and demonstrated capabilities.}
A key achievement of this sequence of works is the systematic characterization of what state-dependent eddy-viscosity closures \emph{can} and
\emph{cannot} achieve when inferred through deterministic variational optimization. The results in~\cite{MatharuProtas2020} demonstrate that,
in simplified settings, adjoint-based reconstruction can successfully recover smooth, physically interpretable state-dependent closures. However,
the more realistic turbulent-flow studies in~\cite{MatharuProtas2022} reveal a fundamental limitation: for certain natural pointwise/instantaneous
matching objectives, the inverse problem may fail to admit a minimizer, even in the absence of noise. This nonexistence result highlights an intrinsic
incompatibility between certain state-dependent closure ansatz and the chosen optimization objectives.

The subsequent analysis in~\cite{MatharuProtas2024} shows that imposing appropriate state constraints on the admissible class of closure functions can
restore existence (attainability) and yield meaningful optimal solutions. Moreover, within the constrained admissible set, the optimizer can be interpreted
as a \emph{variationally optimal} state-dependent closure for the chosen objective. Together, works~\cite{MatharuProtas2020,MatharuProtas2022,MatharuProtas2024}
provide one of the first rigorous demonstrations that state-dependent closure identification in turbulent PDEs is subject to \emph{fundamental performance
limits} that cannot be overcome by algorithmic or numerical refinements alone.

\textbf{Limitations and lessons learned.}
Beyond their specific model problems, these studies deliver several broad lessons for state-dependent inverse problems in complex flows. First,
identifiability alone does not guarantee the existence of an optimal state-dependent closure: the choice of objective functional plays a decisive role
in determining well-posedness. Second, unconstrained functional optimization may lead to nonphysical or nonattainable optima, underscoring the necessity
of incorporating physical admissibility and state constraints directly into the inverse formulation. Third, the results demonstrate that increasing model
fidelity or data availability does not automatically resolve ill-posedness when the inverse problem is structurally incompatible with the underlying
dynamics. In particular, these results highlight that restricting the closure to a single scalar argument $q$ and a purely dissipative eddy-viscosity
form may be insufficient to attain certain optimization targets, even with perfect data, cf.~\cite{MatharuProtas2022}.

From a broader perspective, this application complements the earlier feasibility-focused examples by exposing the analytical and variational limits
of state-dependent parameter identification in turbulent PDEs. It underscores the need for careful problem formulation, objective selection, and
constraint design when extending adjoint-based optimization frameworks to realistic multiphysics and turbulence models, and it delineates a clear
frontier for future research in this area.

\subsection{State-Dependent Transport and Kinetic Law Identification in Battery Models}
\label{sec:apps_echem_2015_2026}

We next consider a family of applications concerned with the identification of \emph{state-dependent transport} and \emph{kinetic laws} in
electrochemical energy-storage models, developed in~\cite{Sethurajan2015,Escalante2020,Daniels2023,Ahmadi2026}. These works extend the classical
deterministic PDE- and ODE-constrained optimization framework for state-dependent parameter identification to battery systems governed by
electrochemical transport, interfacial kinetics, and degradation dynamics. Despite differences in physical setting and modeling fidelity, all
four studies share a common methodological core: nonparametric reconstruction of constitutive relations expressed as functions of evolving state
variables, using adjoint-based gradients and variational optimization. Together, they form a coherent application class within the scope of this
review. Related work~\cite{Richardson2018}, while not itself an inverse formulation, plays an important interpretive role by analyzing the
physical implications and limitations of reconstructed transport laws and is therefore discussed in the context of lessons learned.

\textbf{Problem setting and inverse formulation.}
The forward models considered in~\cite{Sethurajan2015,Escalante2020,Daniels2023,Ahmadi2026} are based on continuum electrochemical descriptions
of battery operation, coupling transport of ionic species, charge conservation, interfacial reaction kinetics, and (in some cases) degradation
mechanisms. While differing in dimensionality and level of model reduction, all formulations can be viewed as instances of nonlinear PDE- or
ODE-constrained systems in which constitutive relations \emph{depend explicitly on evolving state variables}, placing them squarely within the
state-dependent inverse-problem framework discussed at length in Sections~\ref{sec:math}--\ref{sec:comp}.

A representative example, drawn from electrolyte transport modeling~\cite{Sethurajan2015,Escalante2020}, is the ionic species conservation law
(subject to electrochemical coupling constraints)
\begin{subequations}
  \label{eq:echem_transport}
  \begin{align}
    \partial_t c &= \nabla \cdot \bigl( D(c)\,\nabla c \bigr) + \bS,
    && \text{in } \Omega \times (0,T], \label{eq:echem_transport_pde}\\
    - D(c)\,\nabla c \cdot \bn &= \bB,
    && \text{on } \partial\Omega \times (0,T], \label{eq:echem_transport_bc}\\
    c(\bx,0) &= c_0(\bx),
    && \text{in } \Omega, \label{eq:echem_transport_ic}
  \end{align}
\end{subequations}
where $c(\bx,t)$ denotes the electrolyte concentration, $\phi$ is the electric potential, $D(c)$ is a concentration-dependent electrolyte diffusivity,
and $\bS(\cdot)$ and $\bB(\cdot)$ represent volumetric and interfacial source terms associated with electrochemical reactions (their detailed form
being model-dependent and potentially involving additional constitutive functions). Importantly, in~\cite{Sethurajan2015} the \emph{transference number}
$t^+(c)$ appears explicitly in both the source term $\bS(\cdot;\,t^+(c))$ in~\eqref{eq:echem_transport_pde} and in the boundary flux $\bB(\cdot;\,t^+(c))$
in~\eqref{eq:echem_transport_bc}, so that the unknown constitutive function enters both the bulk transport equation and the boundary conditions.

Other state-dependent constitutive relations reconstructed in this application family include:
(i)~transport properties governing lithium motion \emph{inside the graphite insertion material}, inferred from in-situ MRI data via inverse modelling:
in a diffusion-based description, a solid-phase diffusivity $D(c)$ is reconstructed, whereas in a Cahn--Hilliard description the inferred quantities
are coefficients entering the chemical-potential model (and hence the associated Cahn--Hilliard flux), cf.~\cite{Escalante2020}; (ii)~state-dependent
exchange-current densities $i_0(c)$ in Butler--Volmer reaction kinetics, reconstructed from voltage data in~\cite{Daniels2023}; and (iii)~effective
degradation or side-reaction rates $\omega(\bz(t))$, expressed as functions of reduced electrochemical state variables $\bz(t)$ in lumped or
reduced-order battery models~\cite{Ahmadi2026}.

The inverse problems addressed in these studies consist in reconstructing such constitutive laws from experimental measurements or high-fidelity
simulation data, including concentration profiles, terminal voltage curves, and time-averaged or derived observables. In all cases, the control
variable is not a finite set of parameters but an \emph{infinite-dimensional function of the state}, directly paralleling the temperature-dependent
conductivity reconstruction in Section~\ref{sec:apps_heat_2011}, the viscosity reconstruction in Section~\ref{sec:apps_NS_2013}, and the
closure-identification problems discussed in Sections~\ref{sec:apps_ROM_closure_2015} and~\ref{sec:apps_MP_2020_2024}.

As in the preceding applications, identifiability is intrinsically state-dependent: only values of the constitutive functions corresponding to state
regimes explored by the system under the given operating conditions can be informed by the data. The inverse problems are therefore formulated in
reduced space as PDE- or ODE-constrained optimization problems of the form~\eqref{eq:final_opt}, with least-squares misfit functionals subject to
the governing electrochemical model equations~\eqref{eq:echem_transport} and their analogues in the respective studies.

\textbf{Computational framework.}
All four studies employ a classical \emph{optimize--then--discretize} strategy and derive continuous adjoint equations associated with the governing
electrochemical models. Gradients of the reduced cost functional with respect to the constitutive functions are obtained in closed form and exhibit
the characteristic \emph{state-aggregation structure} discussed throughout this review: sensitivity with respect to the value of a constitutive law
at a given state accumulates contributions from space-time regions (or time instances) at which the system attains that state.

In~\cite{Sethurajan2015} and~\cite{Escalante2020}, this structure arises from PDE-based transport models, while~\cite{Daniels2023} and
~\cite{Ahmadi2026} demonstrate analogous behavior in reduced electrochemical and degradation models governed by ODEs. To stabilize the inherently
ill-posed reconstructions, all studies employ Sobolev-type regularization in the space of admissible constitutive functions, effectively replacing
raw $L_2$ gradients with Sobolev gradients obtained via Riesz maps, cf.~Section~\ref{sec:Sobolev_gradients}. Gradient-based optimization algorithms
are then used to solve the resulting infinite-dimensional problems, with each iteration requiring one forward and one adjoint model evaluation.

\textbf{Achievements and demonstrated capabilities.}
A central achievement of this body of work is the demonstration that key electrochemical transport and kinetic laws can be reconstructed directly
as functions of physically meaningful state variables, without imposing restrictive parametric ansatzes. This philosophy mirrors earlier applications
involving temperature-dependent conductivity (Section~\ref{sec:apps_heat_2011}), state-dependent viscosity and constitutive laws in fluid mechanics
(Sections~\ref{sec:apps_NS_2013} and~\ref{sec:apps_vorticity_2015}), as well as closure identification in reduced-order and turbulent-flow models
(Sections~\ref{sec:apps_ROM_closure_2015} and~\ref{sec:apps_MP_2020_2024}). The study~\cite{Sethurajan2015} shows that concentration-dependent
electrolyte transport properties can be inferred from combined electrochemical and imaging data, with in-situ concentration measurements playing
a decisive role as an observation operator enabling functional reconstruction that would be impossible from voltage data alone. The
work~\cite{Escalante2020} demonstrates that voltage data alone may be insufficient to uniquely identify nonlinear diffusion laws, and further shows
how inverse modelling can be used to calibrate and assess competing transport descriptions (diffusion versus Cahn--Hilliard-type models) against
high-resolution MRI measurements, thereby linking constitutive identification with model-class selection.

The work~\cite{Daniels2023} extends state-dependent reconstruction to interfacial kinetics, showing that exchange-current functions can be inferred
variationally from voltage curves and that the resulting reconstructions differ qualitatively from standard parametric models, demonstrating that
even macroscale terminal-voltage observables can carry sufficient information to identify nonparametric kinetic laws when combined with appropriate
modeling assumptions and admissibility constraints. Finally, \cite{Ahmadi2026} generalizes the framework to multi-argument constitutive relations
in degradation modeling, using reduced electrochemical models derived from higher-fidelity descriptions to enable tractable inverse identification
of state-dependent plating, stripping, or side-reaction laws from experimental data.

\textbf{Limitations and lessons learned.}
These studies also reveal limitations characteristic of state-dependent inverse problems in electrochemical systems. First, reconstructions are
fundamentally limited to the range of state variables explored by the experiment or simulation, and extrapolation beyond this range is dominated by
regularization. Second, different choices of observables (e.g., voltage versus concentration data) can lead to markedly different reconstructed
constitutive laws, underscoring the role of experimental design and data richness in ensuring identifiability~\cite{Escalante2020}. In practice,
stable and physically meaningful reconstructions often require the explicit enforcement of admissibility constraints (e.g., non-negativity, boundedness,
or thermodynamic consistency) in addition to Sobolev regularization, particularly when the observation operator is information-poor.

Importantly, independent physical analysis~\cite{Richardson2018} shows that some features observed in reconstructed transport properties (such as
apparent anomalies in concentration-dependent transference numbers) cannot be resolved by minor model refinements (e.g., ion pairing) and instead
signal missing physics or structural model inadequacy. This analysis clarifies that reconstructed constitutive laws should be interpreted as
\emph{effective closures} within a given model class, and that inverse identification may faithfully reproduce data while remaining physically
misleading if the underlying model structure is incomplete. This highlights a broader lesson consistent with other applications in this review:
inverse reconstruction can faithfully identify effective state-dependent laws \emph{within a given model class}, but it \emph{cannot compensate}
for \emph{fundamental modeling deficiencies}.

From a broader perspective, this application demonstrates that classical PDE- and ODE-constrained optimization provides a powerful and unifying
framework for state-dependent constitutive identification in battery models, while simultaneously clarifying the limits of what such reconstructions
can achieve without enriched physics or additional data.

\section{Open Problems, Challenges, and Future Directions}
\label{sec:challenges}

Despite substantial progress in the formulation, analysis, and numerical treatment of inverse problems with state-dependent parameters, as surveyed
in Sections~\ref{sec:math}--\ref{sec:apps}, a number of fundamental challenges continue to limit the robustness, scalability, and predictive
reliability of current computational approaches. Importantly, these challenges do not arise from a lack of variational or adjoint-based methodology,
which is now well established, but rather from intrinsic structural features of inverse problems in which the unknown is an infinite-dimensional
constitutive law depending on the solution of nonlinear PDE and/or ODE systems. The following discussion synthesizes open issues that recur across
analytical theory, numerical algorithms, and applications, and highlights how they manifest in the concrete case studies presented throughout
Sections~\ref{sec:apps_heat_2011}--\ref{sec:apps_echem_2015_2026}.

\subsection{Nonlinearity and Nonconvexity Informed by Applications}
\label{sec:challenges_nonlin_nonconv_apps}

The applications surveyed in Section~\ref{sec:apps} demonstrate that nonlinearity and nonconvexity are not merely theoretical complications of
state-dependent inverse problems, but persistent, structurally induced challenges that shape what can be achieved in practice. While their mathematical
origins were established in Section~\ref{sec:math} and their algorithmic implications discussed in Section~\ref{sec:comp}, the case studies reveal how
these phenomena arise concretely in realistic physical, engineering, and materials-science settings. This subsection synthesizes these observations and
frames them as open challenges that remain only partially addressed by current methodologies, emphasizing structural limitations rather than
algorithm-specific difficulties.

\textbf{Identifiability-induced nonconvexity.}
A central challenge highlighted by essentially all applications is that nonconvexity is inseparable from \emph{state-dependent identifiability}. As shown
analytically in Section~\ref{sec:math}, the Fr\'echet derivative of the parameter--to--observation map is supported only on the identifiability interval
associated with the attained state range. In practice, this leads to reduced cost functionals with flat or weakly curved directions outside this interval,
giving rise to nonuniqueness and multiple local minimizers. This behavior is consistently observed across applications, including temperature-dependent
heat conduction~\cite{Bukshtynov2011}, thermo-fluid viscosity reconstruction~\cite{Bukshtynov2013}, vorticity reconstruction in inviscid
flows~\cite{DanailaProtas2015}, reduced-order closure identification~\cite{ProtasNoack2014,ProtasNoack2015}, turbulent-flow
closures~\cite{MatharuProtas2020,MatharuProtas2022,MatharuProtas2024}, and electrochemical transport and kinetic-law
identification~\cite{Sethurajan2015,Escalante2020,Daniels2023,Ahmadi2026}. Despite substantial differences in governing equations and data modalities,
all these problems admit families of constitutive laws that are indistinguishable on the realized state range. A fundamental open challenge is to
characterize when such nonconvexity can be reduced through experimental design or admissible-set restrictions, and when it is intrinsic and unavoidable.

\textbf{Geometric nonlinearity from state aggregation.}
A second challenge arises from the geometric structure of adjoint-based sensitivities discussed in Section~\ref{sec:var_adj_grads}. In PDE-based applications,
gradients with respect to state-dependent constitutive laws involve integration over codimension-one level sets of the state variable, as shown explicitly
for thermo-fluid systems~\cite{Bukshtynov2013} and inviscid vortex dynamics~\cite{DanailaProtas2015}. Small perturbations of the constitutive law can alter
the geometry, connectivity, or distribution of these level sets, producing highly nonlinear responses in the reduced gradient and Hessian. Analogous
aggregation mechanisms occur in reduced-order models and electrochemical degradation models, where sensitivities accumulate along trajectories in state space
rather than over spatial manifolds~\cite{ProtasNoack2015,Ahmadi2026}. In such settings, the optimization landscape depends sensitively on how frequently and
how uniformly the system visits different regions of the state space. A persistent challenge is to develop optimization and continuation strategies that
explicitly account for this geometric coupling between constitutive laws and state evolution, rather than treating it as a numerical side effect.

\textbf{Objective-dependent nonexistence and structural limits.}
The turbulent-flow studies reviewed in Section~\ref{sec:apps_MP_2020_2024}~\cite{MatharuProtas2020,MatharuProtas2022,MatharuProtas2024} expose a deeper and
less explored challenge: nonconvexity may manifest as \emph{nonexistence of minimizers}. In these works, certain natural matching objectives lead to inverse
problems that admit no optimal state-dependent eddy-viscosity closure within unconstrained admissible classes, even in the absence of noise and with exact
adjoints. These results demonstrate that some inverse formulations are variationally incompatible with the chosen constitutive ansatz, revealing fundamental
limits that cannot be overcome by regularization tuning or improved solvers.While the imposition of explicit state constraints can restore
existence~\cite{MatharuProtas2024}, a general theoretical framework for predicting objective-dependent nonattainability in infinite-dimensional state-dependent
inverse problems remains an open problem.

\textbf{Dependence on data modality and excitation.}
Applications further demonstrate that the severity of nonconvexity depends strongly on the choice of observables and experimental or computational excitation.
In electrochemical modeling, reconstructions based solely on terminal-voltage data~\cite{Escalante2020,Daniels2023} exhibit markedly flatter optimization
landscapes than those incorporating internal concentration measurements~\cite{Sethurajan2015}. Similarly, in reduced-order closure
identification~\cite{ProtasNoack2015}, limited or weakly exciting training trajectories lead to closures that perform well locally but fail to generalize.
These observations point to an open challenge at the interface of inverse problems and experimental design: how to systematically design experiments or
data-acquisition strategies that enrich state-space coverage and thereby mitigate nonconvexity, without relying on trial-and-error adjustments of boundary
conditions or forcing terms.

\textbf{Model-form error and effective optimality.}
Several applications underscore that reconstructed state-dependent constitutive laws should be interpreted as \emph{effective closures} within a prescribed
model class. In battery applications, independent physical analysis~\cite{Richardson2018} shows that some features of reconstructed transport properties
cannot be reconciled with established physics, indicating compensation for missing model ingredients. Similar conclusions arise in turbulent-flow
modeling~\cite{MatharuProtas2022}, where optimal closures may fail to exist for certain objectives because the closure ansatz itself is too restrictive.
A key open challenge is therefore to distinguish nonconvexity arising from limited data or identifiability from that induced by model-form inadequacy.
Without such distinction, optimization may converge reliably while producing constitutive laws that are mathematically optimal yet physically questionable.

\textbf{Outlook.}
Taken together, the applications reviewed in Section~\ref{sec:apps} demonstrate that nonlinearity and nonconvexity in state-dependent inverse problems
are driven by intertwined effects of identifiability, geometry, objective selection, data modality, and model adequacy. While current computational
frameworks can manage these effects to a degree, they do not yet provide systematic criteria for when nonconvexity can be mitigated and when it reflects
fundamental limits. Addressing these questions remains a central challenge for future research at the intersection of PDE-constrained optimization,
experimental design, and physics-based modeling.

\subsection{Identifiability and Ill-Posedness}
\label{sec:challenges_identif_ill_posed}

A recurring and unifying challenge across all state-dependent inverse problems surveyed in Sections~\ref{sec:apps_heat_2011}--\ref{sec:apps_echem_2015_2026}
is that \emph{identifiability of the unknown constitutive relation is intrinsically restricted to the range of state values actually explored by the system}.
This limitation is not an artifact of discretization, regularization choice, or optimization algorithm, but a structural property of the
parameter--to--observation map, rooted in the mathematical formulation of the inverse problem itself.

\textbf{Identifiability intervals and loss of injectivity.}
From a mathematical standpoint, Section~\ref{sec:math} establishes that when a constitutive relation is expressed as a function of a state variable,
the Fr\'echet derivative of the parameter--to--observation operator is supported only on the \emph{identifiability interval} $\mI \subset \mD$ associated
with the attained state range. Variations of the constitutive law outside this interval do not induce first-order changes in the state or in the observations.
As a result, the inverse map is non-injective on the interval $\mD$ of desired reconstruction, and the inverse problem is underdetermined even in the absence
of noise. This phenomenon was demonstrated explicitly in the stationary heat-conduction problem with temperature-dependent conductivity~\cite{Bukshtynov2011},
where reconstructions are provably meaningful only on the temperature interval realized by the forward solution. Outside this interval, the reconstructed
conductivity is determined entirely by regularization and admissible-set assumptions, rather than by data.

\textbf{Ill-posedness beyond noise amplification.}
While classical inverse-problem theory often associates ill-posedness with sensitivity to measurement noise, the applications reviewed here show that
state-dependent inverse problems are ill-posed in a fundamentally stronger sense. Even with exact data, perfect forward models, and adjoint-consistent gradients,
the inverse problem may admit infinitely many solutions or fail to possess a minimizer. This distinction is already implicit in the analytical discussion of
Section~\ref{sec:math}, but becomes concrete in applications discussed in Section~\ref{sec:apps}. In the thermo-fluid viscosity reconstruction
problem~\cite{Bukshtynov2013}, distinct viscosity functions that coincide on the realized temperature range produce identical observations, leading to flat
directions in the reduced cost functional. In reduced-order closure identification~\cite{ProtasNoack2014,ProtasNoack2015}, closures inferred from limited
trajectories are non-unique outside the explored energy range, even when the data are noise-free.

\textbf{Nonexistence as an extreme form of ill-posedness.}
The turbulent-flow closure studies of Section~\ref{sec:apps_MP_2020_2024}~\cite{MatharuProtas2020,MatharuProtas2022,MatharuProtas2024} reveal an even more
severe manifestation of ill-posedness: for certain natural objective functionals, the inverse problem does not admit a minimizer within the chosen admissible
class of state-dependent closures. Importantly, this nonexistence is not caused by insufficient data or numerical instability, but by a structural
incompatibility between the objective and the closure ansatz. These results elevate ill-posedness from a numerical concern to a \emph{variational obstruction},
showing that identifiability failure can preclude not only uniqueness but also attainability of optimal solutions. The restoration of existence through
explicit state constraints~\cite{MatharuProtas2024} underscores that admissible-set design is a central component of well-posedness in state-dependent inverse
problems.

\textbf{Algorithmic mitigation versus structural limits.}
Section~\ref{sec:comp} discusses a range of computational strategies (e.g., Sobolev gradients, Tikhonov regularization, admissible-set enforcement, and
experimental design) that can stabilize reconstructions and mitigate the practical consequences of ill-posedness. The heat conduction
example~\cite{Bukshtynov2011} demonstrates that Sobolev-gradient preconditioning and $H^1$-type regularization yield stable and mesh-independent reconstructions
on the identifiability interval. Similarly, in battery-model inversions (Section~\ref{sec:apps_echem_2015_2026}), regularization and admissibility constraints
enable robust recovery of transport and kinetic laws within the explored concentration or state ranges. However, these techniques do not eliminate ill-posedness;
rather, they select a \emph{preferred representative} from an equivalence class of admissible solutions. As such, they address stability but not identifiability
itself.

\textbf{Identifiability as an experimental-design problem.}
Several applications illustrate that identifiability can be partially improved by modifying the experimental or computational setup to expand the explored
state range. The identifiability shifting strategy introduced in~\cite{Bukshtynov2011} and further developed in Section~\ref{sec:shifting_I} exemplifies this
idea: by altering boundary conditions or forcing, one can generate a sequence of experiments whose combined identifiability intervals cover a larger target
domain. Analogous observations arise in electrochemical modeling, where the use of internal concentration measurements~\cite{Sethurajan2015} dramatically
improves identifiability compared with voltage-only data~\cite{Escalante2020}. These examples highlight that identifiability is not solely a property of
the inverse problem formulation, but also of the data-acquisition strategy.

\textbf{Open challenges.}
Despite these advances, a general theory connecting identifiability, admissible-set design, objective selection, and experimental excitation remains elusive.
In particular, there is no systematic criterion for determining when identifiability limitations can be overcome through experimental design, and when they
reflect intrinsic information barriers imposed by the model structure. Developing such criteria (and embedding them into optimization and design frameworks)
remains a central open challenge for state-dependent inverse problems.

\subsection{Regularization and Stabilization Strategies}
\label{sec:challenges_reg_stabiliz}

All state-dependent inverse problems reviewed in Section~\ref{sec:apps} rely critically on regularization to obtain stable and interpretable reconstructions
of constitutive laws. This dependence is not incidental: as established in Section~\ref{sec:math}, the inverse maps under consideration are typically
non-injective, and as discussed in Section~\ref{sec:comp}, adjoint-based optimization without regularization leads to severe instability, mesh dependence,
and amplification of data errors. The applications surveyed in Sections~\ref{sec:apps_heat_2011}--\ref{sec:apps_echem_2015_2026} demonstrate both the
effectiveness and the limitations of current regularization strategies, and highlight several unresolved challenges.

\textbf{Regularization as an intrinsic component of well-posedness.}
From the perspective of classical inverse-problem theory, regularization is required whenever the inverse map fails to be continuous or injective; see,
e.g., \cite{EnglHankeNeubauer1996}. State-dependent inverse problems fall squarely into this category: even in noise-free settings, the restriction of
identifiability to state-dependent intervals renders the inverse problem underdetermined outside the realized state range. In applications such as
temperature-dependent heat conduction~\cite{Bukshtynov2011}, thermo-fluid viscosity reconstruction~\cite{Bukshtynov2013}, and vorticity reconstruction in
inviscid flows~\cite{DanailaProtas2015}, regularization is therefore not merely a noise-mitigation tool, but an essential mechanism for selecting
a physically meaningful representative from a family of admissible solutions. This interpretation aligns closely with the classical role of Tikhonov
regularization as a selection principle rather than a cure for ill-posedness.

\textbf{Sobolev regularization and metric selection in function space.}
Across essentially all applications reviewed, Sobolev-type regularization plays a central role. In the heat-conduction and thermo-fluid
examples~\cite{Bukshtynov2011,Bukshtynov2013}, raw $L_2$ gradients with respect to the constitutive law are discontinuous and confined to the
identifiability interval; Sobolev gradients, obtained via Riesz maps on $H^1(\mD)$ or higher-order spaces, provide both smoothing and mesh-independent
descent directions. Similar mechanisms are essential in vorticity reconstruction~\cite{DanailaProtas2015}, reduced-order closure
identification~\cite{ProtasNoack2014,ProtasNoack2015}, and turbulent-flow inversions (Section~\ref{sec:apps_MP_2020_2024}), where $H^1$- and $H^2$-type
metrics suppress spurious oscillations induced by state aggregation. Despite their empirical success, the choice of Sobolev order, length-scale parameters,
and weighting remains largely heuristic and problem-specific, reflecting the absence of general principles for metric selection in state-dependent inverse
problems.

\textbf{Admissible-set constraints and constraint-induced stabilization.}
Several applications demonstrate that regularization through admissible-set restrictions is not merely beneficial but sometimes indispensable. In
turbulent-flow closure reconstruction, unconstrained optimization may fail to admit a minimizer altogether~\cite{MatharuProtas2022}; imposing explicit state
constraints restores existence and yields variationally meaningful solutions~\cite{MatharuProtas2024}. From a broader perspective, such constraints act as
a form of nonlinear regularization, shaping the feasible set to enforce physical admissibility (e.g., non-negativity, boundedness, or monotonicity) that cannot
be guaranteed by quadratic penalties alone. Designing admissible sets that balance physical realism with mathematical tractability remains an open challenge,
particularly in multiphysics and turbulence applications.

\textbf{Penalty and constraint-relaxation strategies.}
An alternative stabilization paradigm, complementary to classical Tikhonov regularization, is provided by penalty-based formulations of PDE-constrained
optimization. The work of van~Leeuwen and Herrmann~\cite{Leeuwen2016} demonstrates that relaxing the PDE constraint via a quadratic penalty can significantly
reduce the effective nonlinearity of the optimization landscape and mitigate sensitivity to initialization, while retaining computational complexity comparable
to reduced formulations. From the standpoint of state-dependent inverse problems, such penalty methods suggest a promising route for balancing data fidelity,
model enforcement, and regularization. At the same time, they introduce new challenges: the choice of penalty parameters, continuation strategies, and their
interaction with function-space regularization remain poorly understood in infinite-dimensional, state-dependent settings. A systematic analysis of how penalty
relaxation interacts with identifiability intervals and Sobolev metrics is still lacking.

\textbf{Regularization versus modeling error.}
Several applications, particularly in electrochemical modeling (Section~\ref{sec:apps_echem_2015_2026}) and turbulent-flow closure reconstruction
(Section~\ref{sec:apps_MP_2020_2024}), underscore that regularization may compensate not only for ill-posedness but also for deficiencies in the underlying
model class itself. In battery applications~\cite{Escalante2020,Daniels2023}, strong regularization stabilizes the inverse problem even when available observables
are indirect, but the resulting transport or kinetic laws may reflect regularized surrogates rather than physically complete descriptions. Analogous behavior is
observed in turbulent-flow inversions~\cite{MatharuProtas2022}, where regularization and admissible constraints yield stable optimizers despite structural
limitations of the eddy-viscosity ansatz. In such cases, increasingly strong regularization yields stable reconstructions that fit the data yet encode effective
rather than physically faithful constitutive laws. Distinguishing regularization-induced bias from genuine constitutive features remains a major open problem,
especially when multiple regularization mechanisms (Sobolev smoothing, admissible-set constraints, and penalty relaxation) are combined.

\textbf{Open challenges.}
Despite significant progress, regularization in state-dependent inverse problems remains largely artisanal. There is no general theory linking identifiability
structure, noise level, objective functional, admissible-set design, and regularization strength in a predictive way. Developing such a theory (and translating
it into practical guidelines for large-scale, multiphysics inverse problems) is a central challenge for future research.

\subsection{Adjoint and Gradient Computation for State-Dependent Parameters}
\label{sec:challenges_adj_grad}

The derivation and numerical realization of adjoint equations for state-dependent inverse problems are substantially more involved than in classical
parameter-identification settings, where the control enters the governing equations explicitly and linearly. In state-dependent problems, the unknown
constitutive law influences the forward model only indirectly through the state, leading to gradient expressions that are nonlinear, nonlocal, and
tightly coupled to the solution of the forward problem. As established in Section~\ref{sec:math} and illustrated throughout Section~\ref{sec:apps},
this structural complexity shapes both the analytical form of the adjoint system and the numerical strategies required for its robust evaluation.

\textbf{Indirect dependence and variational complexity.}
From a variational standpoint, the key difficulty is that perturbations of the constitutive law affect the objective only through induced changes in
the state. Section~\ref{sec:var_adj_framework} shows that the resulting Fr\'echet derivatives involve compositions of the linearized forward operator
with state-dependent mappings, rather than simple pointwise sensitivities. Even in the scalar stationary heat-conduction problem~\cite{Bukshtynov2011},
this indirect dependence necessitates special analytical devices (e.g., Kirchhoff-type transformations) to obtain adjoint equations that are well defined
in function space. In coupled multiphysics systems, such as thermo-fluid viscosity reconstruction~\cite{Bukshtynov2013} or electrochemical transport models
(Section~\ref{sec:apps_echem_2015_2026}), the same issue manifests through adjoint equations that couple multiple fields and inherit the full nonlinearity
of the forward model.

\textbf{State aggregation and level-set sensitivity structures.}
A defining feature of adjoint-based gradients in state-dependent inverse problems is their \emph{state-aggregation structure}, discussed abstractly in
Section~\ref{sec:var_adj_grads}. In PDE-based applications, sensitivities with respect to the constitutive law at a given state value accumulate contributions
from spatial regions where the state attains that value. This structure appears explicitly in thermo-fluid viscosity reconstruction, where gradients involve
integration over temperature level sets~\cite{Bukshtynov2013}, and in inviscid vorticity reconstruction, where sensitivities are aggregated over streamfunction
level sets~\cite{DanailaProtas2015}. In reduced-order and turbulent-flow closure identification, analogous aggregation occurs along trajectories in state space
rather than over spatial manifolds~\cite{ProtasNoack2015,MatharuProtas2022}. While this structure enables infinite-dimensional control with adjoint costs
independent of control discretization, it introduces geometric and numerical challenges that have no counterpart in standard distributed-parameter optimization.

\textbf{Discretization, consistency, and gradient verification.}
Section~\ref{sec:comp} emphasizes that adjoint-based optimization is only as reliable as the consistency between the continuous adjoint formulation and its
discrete realization. This issue is particularly acute for state-dependent gradients, where level-set or trajectory-based aggregation must be approximated
numerically. In~\cite{Bukshtynov2011} and~\cite{Bukshtynov2013}, extensive $\kappa$-tests are used to verify gradient correctness across discretization levels,
highlighting the importance of systematic verification. In more complex settings, such as turbulent-flow closure reconstruction~\cite{MatharuProtas2022},
gradient evaluation itself becomes a dominant computational task, and small inconsistencies in discretization or time integration can lead to qualitatively
incorrect descent directions. Developing discretization strategies that preserve the structure of continuous adjoint sensitivities remains an open challenge,
especially in time-dependent and multiphysics contexts.

\textbf{Computational cost and scalability.}
Although adjoint formulations provide algorithmic scalability in the sense of one adjoint solve per optimization iteration, the absolute computational cost
can be substantial in state-dependent problems. In unsteady thermo-fluid and electrochemical models, each optimization iteration requires a full forward
simulation and a backward-in-time adjoint solve, often with additional overhead associated with storing or recomputing the forward trajectory, as demonstrated
in thermo-fluid viscosity reconstruction~\cite{Bukshtynov2013} and in electrochemical transport and kinetic-law identification problems
(Section~\ref{sec:apps_echem_2015_2026}). In turbulent-flow applications (Section~\ref{sec:apps_MP_2020_2024}), long-time integrations and statistically
averaged objectives further complicate adjoint computation and increase sensitivity to numerical errors, as demonstrated in large-eddy and turbulence-closure
reconstruction studies such as~\cite{MatharuProtas2022,MatharuProtas2024}. Balancing accuracy, memory usage, and computational efficiency in adjoint
implementations remains a major practical challenge for large-scale applications.

\textbf{Automated differentiation and algorithmic adjoints.}
Automated differentiation (AD) and algorithmic adjoint techniques offer an attractive alternative to hand-derived adjoints, particularly for complex multiphysics
models; see, e.g., \cite{GriewankWalther2008,GilesPierce2000,Naumann2011}. However, their deployment in state-dependent inverse problems is far from straightforward.
AD tools typically operate at the level of discrete code and generate sensitivities for a fully discretized problem, which may obscure the underlying function-space
structure emphasized in Section~\ref{sec:math}. In problems involving level-set integration, trajectory aggregation, or nonsmooth admissible-set constraints, naive
application of AD can yield gradients that are formally correct for the discrete problem yet inconsistent with the intended continuous formulation. Bridging the gap
between algorithmic adjoints and the function-space adjoint framework used throughout this review remains an open methodological challenge.

\textbf{Open challenges.}
Despite significant progress, adjoint and gradient computation for state-dependent parameters remains one of the most technically demanding aspects of this class
of inverse problems. Key open questions include how to systematically preserve state-aggregation structures under discretization, how to verify adjoint correctness
in large-scale multiphysics simulations, and how to combine hand-derived functional adjoints with automated tools without sacrificing interpretability or
robustness. Addressing these challenges is essential for extending state-dependent inverse problem methodologies to increasingly realistic and computationally
demanding applications.

\subsection{Handling Noisy, Sparse, or Indirect Observations}
\label{sec:challenges_observe}

A defining challenge in many state-dependent inverse problems is that available observations are often noisy, sparse, or indirect, rather than full-field
measurements of the state. As discussed in Section~\ref{sec:apps}, this situation is not exceptional but typical of realistic applications in physics,
engineering, and materials science. The resulting limitations interact in subtle ways with identifiability, nonconvexity, and regularization, and fundamentally
constrain what can be recovered.

\textbf{Observation operators and loss of information.}
From a mathematical perspective, the difficulty arises because the observation operator maps the state to a low-dimensional or indirect set of observables,
further compressing the information available for inversion. Section~\ref{sec:math} shows that in state-dependent inverse problems, identifiability is already
restricted to the range of state values attained by the system; sparse or indirect observations exacerbate this restriction by reducing sensitivity even within
the identifiable interval. As a result, the Fr\'echet derivative of the parameter--to--observation map may be severely rank-deficient, leading to flat directions
in the reduced cost functional and amplified nonuniqueness. This structural loss of information cannot be remedied by algorithmic improvements alone.

\textbf{Noise amplification and regularization dependence.}
Section~\ref{sec:comp} emphasizes that adjoint-based gradients propagate observational noise through the adjoint equations, often amplifying high-frequency
components of the reconstructed constitutive law. In applications with indirect data, such as boundary measurements or integral observables, this effect is
particularly pronounced. The heat-conduction example~\cite{Bukshtynov2011} already demonstrates that increasing noise levels induce oscillatory instabilities
unless explicit regularization is employed. Similar behavior is observed in thermo-fluid viscosity reconstruction~\cite{Bukshtynov2013} and in reduced-order
closure identification~\cite{ProtasNoack2015}, where noise and limited excitation lead to highly sensitive optimization landscapes. In such settings,
regularization does not merely improve stability but dominates the reconstruction outside a narrow state range.

\textbf{Indirect observations in electrochemical models.}
Electrochemical inverse problems provide a particularly clear illustration of the challenges posed by indirect data. In battery-model inversions
(Section~\ref{sec:apps_echem_2015_2026}), terminal voltage measurements integrate the effects of transport, reaction kinetics, and thermodynamics over the entire
cell, making it difficult to disentangle individual constitutive contributions. The study~\cite{Escalante2020} shows that voltage data alone may be insufficient
to uniquely identify nonlinear transport laws, even with perfect forward models, whereas incorporating internal concentration measurements~\cite{Sethurajan2015}
dramatically improves identifiability. These results highlight that data modality, not just data quantity, plays a decisive role in state-dependent inversion.

\textbf{Statistical observables and long-time averaging.}
In turbulent-flow closure reconstruction (Section~\ref{sec:apps_MP_2020_2024}), observations often take the form of time-averaged or statistical quantities
rather than instantaneous flow fields. While such observables are physically meaningful and experimentally accessible, they further reduce sensitivity to
pointwise variations in the constitutive law. The analysis in~\cite{MatharuProtas2022} demonstrates that, for certain natural statistical objectives, the
resulting inverse problem may fail to admit a minimizer, revealing a fundamental incompatibility between the observation operator and the chosen closure
ansatz. This phenomenon underscores that data reduction through averaging can qualitatively alter the variational structure of the inverse problem, not merely
degrade reconstruction accuracy.

\textbf{Algorithmic mitigation and its limits.}
A range of computational strategies discussed in Section~\ref{sec:comp}, including Sobolev regularization, admissible-set constraints, and experimental design,
can mitigate the adverse effects of noisy or indirect data. For example, identifiability shifting through modified boundary conditions in heat
conduction~\cite{Bukshtynov2011} and the use of multiple operating regimes in electrochemical models~\cite{Sethurajan2015} can enrich the information content
of observations. However, these strategies do not overcome fundamental information barriers imposed by the observation operator itself. When data are
insufficiently informative, optimization converges to regularization-dominated solutions whose predictive validity is limited.

\textbf{Open challenges.}
Despite substantial progress, there is no general framework for quantifying the information content of noisy, sparse, or indirect observations in
state-dependent inverse problems, nor for predicting when a given data modality suffices to identify a constitutive law of interest. Developing principled
connections between observation design, identifiability analysis, and regularization strategies remains a major open challenge. Addressing this gap is
essential for extending state-dependent inverse problem methodologies to realistic experimental and industrial settings where full-field data are rarely
available.

\subsection{Extension to Time-Dependent and Strongly Nonlinear PDEs}
\label{sec:challenges_time_nonlin}

While a substantial analytical theory exists for inverse problems governed by scalar elliptic and parabolic equations with state-dependent parameters,
cf.~Section~\ref{sec:math}, comparatively little is known about such problems in strongly nonlinear, time-dependent PDEs and coupled multiphysics systems.
This gap is not merely technical: it reflects fundamental difficulties associated with nonlinear dynamics, memory effects, and multiscale interactions that
are only partially captured by existing theory. As illustrated by the applications reviewed in Section~\ref{sec:apps}, these challenges become central in
realistic fluid-dynamical, turbulent, and electrochemical settings.

\textbf{Limits of existing analytical theory.}
Section~\ref{sec:math} establishes key properties of state-dependent inverse problems in relatively simple settings, including identifiability intervals,
differentiability of the parameter--to--state map, and adjoint-based gradient representations. However, many of these results rely on assumptions (such as
monotonicity, ellipticity, or compactness of solution operators) that do not extend naturally to strongly nonlinear or convective systems. In time-dependent
problems, additional difficulties arise from temporal nonlocality and the accumulation of sensitivity over long time intervals. As a result, rigorous
statements about identifiability, stability, or uniqueness are largely unavailable for inverse problems governed by nonlinear evolution equations such as the
Navier--Stokes system (Sections~\ref{sec:apps_NS_2013} and~\ref{sec:apps_MP_2020_2024}) or coupled electrochemical transport and kinetic models
(Section~\ref{sec:apps_echem_2015_2026}).

\textbf{Challenges in unsteady and nonlinear adjoint analysis.}
From a computational standpoint, Section~\ref{sec:comp} demonstrates that adjoint-based optimization remains feasible for time-dependent and nonlinear PDEs,
but the associated adjoint equations inherit the full complexity of the forward dynamics. In unsteady thermo-fluid problems~\cite{Bukshtynov2013}, the adjoint
system is backward in time, nonlinear, and coupled across multiple fields, amplifying both numerical sensitivity and implementation complexity. In turbulent-flow
applications (Section~\ref{sec:apps_MP_2020_2024}), adjoint equations must be integrated over long time horizons and are sensitive to chaotic dynamics, raising
questions about their reliability and interpretability when used to compute gradients of long-time or statistical objectives~\cite{MatharuProtas2022}. These
issues highlight a disconnect between the availability of formal adjoint equations and the existence of a robust analytical framework guaranteeing their
effectiveness.

\textbf{Nonlinearity, coupling, and multiscale effects.}
Strong nonlinearity and multiphysics coupling further complicate the extension of existing theory. In Navier--Stokes-based inverse
problems~\cite{Bukshtynov2013,MatharuProtas2020}, the state-dependent parameter influences the dynamics through nonlinear transport terms, making sensitivity
propagation highly state dependent and potentially intermittent in time. In electrochemical models (Section~\ref{sec:apps_echem_2015_2026}), state-dependent
transport and kinetic laws are coupled across bulk transport, interfacial reactions, and electric fields, often spanning disparate time and length
scales~\cite{Sethurajan2015,Escalante2020,Daniels2023}. These multiscale interactions challenge both analytical tractability and numerical robustness, and
complicate the interpretation of reconstructed constitutive laws as intrinsic material properties rather than effective, regime-dependent closures.

\textbf{Interaction with identifiability and regularization.}
The absence of rigorous theory in strongly nonlinear, time-dependent settings exacerbates the challenges discussed in earlier subsections. Identifiability
limitations become harder to diagnose a priori, regularization choices become more influential, and nonconvexity is more pronounced. As observed in
turbulent-flow closure reconstruction~\cite{MatharuProtas2022,MatharuProtas2024}, even the existence of minimizers may depend sensitively on the choice of
objective functional and admissible set. In such contexts, regularization and admissibility constraints play a dual role: they stabilize computations while
implicitly compensating for gaps in analytical understanding. Distinguishing between these two roles remains difficult without a supporting theoretical
framework.

\textbf{Open challenges.}
Extending the mathematical foundations of state-dependent inverse problems to strongly nonlinear, time-dependent, and coupled PDE systems remains a major
open challenge. Key questions include how to characterize identifiability and stability in nonlinear dynamical systems, how to interpret adjoint sensitivities
in chaotic or multiscale regimes, and how to design inverse formulations whose solutions are both mathematically well posed and physically meaningful. Progress
in this direction will likely require new analytical tools that bridge inverse-problem theory, nonlinear dynamics, and multiphysics modeling, as well as close
interaction between theory, computation, and application-driven insight.

\subsection{Benchmark Problems and Reproducibility}
\label{sec:challenges_bench_reproduct}

The diversity of governing equations, objective functionals, discretization choices, and regularization strategies employed across the applications reviewed in
Sections~\ref{sec:apps_heat_2011}--\ref{sec:apps_echem_2015_2026} highlights both the maturity and the fragmentation of the field. While this diversity reflects
the broad applicability of state-dependent inverse problem methodologies, it also makes systematic comparison of algorithms and theoretical claims difficult.
The absence of standardized benchmark problems tailored to state-dependent inverse problems remains a significant obstacle to reproducibility, validation, and
cumulative progress.

\textbf{Heterogeneity of problem formulations.}
As emphasized in Sections~\ref{sec:math} and~\ref{sec:comp}, state-dependent inverse problems differ fundamentally from classical coefficient-identification
problems: the unknown is an infinite-dimensional function of the state, identifiability is restricted to state-dependent intervals, and gradients exhibit
aggregation over level sets or trajectories. These structural features are instantiated in markedly different ways across applications. For example, the
canonical heat-conduction problem~\cite{Bukshtynov2011} involves a stationary elliptic PDE with full-field temperature data, whereas thermo-fluid viscosity
reconstruction~\cite{Bukshtynov2013} and turbulent-flow closure identification~\cite{MatharuProtas2022,MatharuProtas2024} involve unsteady, nonlinear dynamics
and objectives based on sparse or statistical observables. In electrochemical applications (Section~\ref{sec:apps_echem_2015_2026}), additional complexity
arises from coupled transport and kinetic laws and indirect voltage measurements. This heterogeneity complicates the design of benchmarks that are both
representative and broadly accessible.

\textbf{Verification of adjoints and gradients.}
Section~\ref{sec:verif_kappa} underscores that adjoint-based optimization is reliable only when gradients are verified and discretizations are consistent with the
continuous formulation. Among the applications reviewed, the heat-conduction and thermo-fluid studies~\cite{Bukshtynov2011,Bukshtynov2013} stand out for introducing
systematic $\kappa$-tests to verify adjoint-based gradients. Closely related gradient-verification procedures were subsequently adopted in inviscid vorticity
reconstruction~\cite{DanailaProtas2015} and in reduced-order and turbulent-flow closure identification~\cite{MatharuProtas2020,MatharuProtas2022,MatharuProtas2024},
where such tests are used both to confirm gradient correctness and to diagnose the onset of ill-posedness as regularization is reduced. More recently, analogous
verification has also been employed in electrochemical and kinetic-law identification problems~\cite{Ahmadi2026}, indicating that rigorous adjoint validation has
become standard practice across diverse classes of state-dependent inverse problems. By contrast, many large-scale or multiphysics applications necessarily rely
on more complex codes and workflows, where such verification is harder to perform and less frequently reported. The lack of common benchmark problems with known
reference solutions or controlled identifiability properties makes it difficult to assess whether discrepancies between methods reflect genuine algorithmic
differences or simply implementation artifacts.

\textbf{Role of regularization and admissible sets.}
Benchmarking is further complicated by the central role of regularization and admissible-set design in state-dependent inverse problems. As discussed in
Sections~\ref{sec:challenges_identif_ill_posed} and~\ref{sec:challenges_reg_stabiliz}, regularization does not merely stabilize computations but selects
a representative from an equivalence class of admissible constitutive laws. Different choices of Sobolev metrics, regularization parameters, or state
constraints can therefore lead to qualitatively different reconstructions even when applied to the same data. Without agreed-upon benchmark settings specifying
not only the forward model and data but also the admissible set and regularization structure, meaningful comparison of reconstruction quality remains elusive.

\textbf{Reproducibility across scales and model classes.}
The applications reviewed span a wide range of spatial and temporal scales, from one-dimensional model problems to large-eddy simulations of turbulent flows
and coupled electrochemical models. Reproducing such results often requires access to specialized solvers, meshes, and experimental or simulation data, which
are not always publicly available. This barrier is particularly acute for time-dependent and chaotic systems, where small differences in discretization or time
integration can lead to divergent trajectories and substantially different inverse outcomes (Sections~\ref{sec:apps_NS_2013} and~\ref{sec:apps_MP_2020_2024}).
Establishing benchmark problems that are computationally tractable yet retain the defining features of state-dependent inversion would lower the barrier to entry
and facilitate independent verification of results.

\textbf{Open challenges.}
Developing a suite of standardized benchmark problems for state-dependent inverse problems remains an open and pressing challenge. Such benchmarks should span
multiple levels of complexity, from analytically tractable scalar problems with known identifiability intervals to simplified time-dependent and multiphysics
models capturing essential nonlinear and state-dependent features. Equally important is the dissemination of reference implementations, verified adjoint
formulations, and clearly specified regularization and admissible-set choices. Addressing these needs would significantly enhance reproducibility, enable
objective comparison of emerging methods, and support the development of a more cumulative and transparent research culture in state-dependent inverse problems.

\subsection{Integration with Data-Driven or Hybrid Methods}
\label{sec:challenges_integr_method}

The increasing complexity of state-dependent inverse problems reviewed in Sections~\ref{sec:apps} and~\ref{sec:challenges_time_nonlin}, together with limited,
noisy, or indirect observations, has motivated growing interest in hybrid approaches that combine adjoint-based PDE-constrained optimization with data-driven
representations. As discussed in Section~\ref{sec:non_traditional_OPT}, such methods offer the potential to improve scalability, expressiveness, and computational
efficiency. At the same time, they introduce new challenges related to interpretability, stability, and the preservation of the mathematical structure underpinning
classical inverse-problem formulations.

\textbf{Motivation from state-dependent inverse problems.}
From the perspective of Sections~\ref{sec:math} and~\ref{sec:comp}, state-dependent inverse problems already exhibit severe nonlinearity, nonconvexity, and
identifiability limitations, even in purely deterministic settings. In many of the applications reviewed in Section~\ref{sec:apps}, the unknown constitutive law
is an infinite-dimensional function of the state, and its influence on the observations is mediated through complex, often time-dependent dynamics. These features
motivate the use of flexible function representations and surrogate models capable of capturing rich state-dependent behavior without resorting to overly restrictive
parametric ansatzes. However, introducing such flexibility raises fundamental questions about how learning-based components interact with identifiability intervals,
adjoint sensitivities, and regularization mechanisms.

\textbf{Hybrid representations within deterministic optimization.}
One principled route for integrating data-driven ideas into state-dependent inverse problems is to embed flexible function representations directly within
a classical adjoint-based optimization framework. The example of~\cite{Akerson2025}, discussed in Section~\ref{sec:non_traditional_OPT}, illustrates this approach:
neural networks are used as high-capacity parameterizations of state-dependent constitutive laws, while the governing PDE is enforced exactly and all network
parameters are optimized deterministically via reduced-space gradients. From the standpoint of Sections~\ref{sec:math} and~\ref{sec:challenges_adj_grad}, this
strategy preserves the function-space adjoint structure, identifiability constraints, and convergence properties of PDE-constrained optimization, while significantly
enriching the space of admissible constitutive relations. At the same time, it highlights an open challenge: understanding how classical notions of regularization,
admissible sets, and identifiability extend to learned function spaces with thousands or millions of degrees of freedom.

\textbf{Surrogates, reduced models, and algorithmic acceleration.}
A complementary class of hybrid methods replaces or augments expensive PDE solves with learned surrogates or reduced-order models embedded within optimization
loops~\cite{Xu2022}. Such approaches are particularly attractive in large-scale, time-dependent, or multiphysics settings, where the computational cost of repeated
forward and adjoint solves dominates the inversion process, as seen in turbulent-flow and electrochemical applications (Sections~\ref{sec:apps_MP_2020_2024}
and~\ref{sec:apps_echem_2015_2026}). While surrogate-assisted optimization can yield substantial speedups, it raises nontrivial challenges related to consistency:
surrogate errors may bias gradients, distort the optimization landscape, or obscure the identifiability structure of the original inverse problem. Ensuring that
surrogate models respect physical constraints and provide reliable sensitivities remains an open methodological issue.

\textbf{Interaction with uncertainty quantification and learning-based inference.}
As emphasized in Section~\ref{sec:non_traditional_OPT}, Bayesian and learning-based approaches offer a complementary lens on state-dependent inverse problems
by explicitly representing uncertainty and incorporating prior information. Hybrid strategies that combine adjoint-based optimization with probabilistic inference
or learning-based regularization blur the boundary between deterministic and statistical formulations. While such combinations can provide valuable uncertainty
estimates and improved robustness, they also complicate interpretation: it becomes less clear whether reconstructed constitutive laws reflect identifiable physical
properties, regularization-induced preferences, or artifacts of the learning architecture. Developing theoretical tools to disentangle these effects in state-dependent
settings remains an open challenge.

\textbf{Preserving interpretability and physical structure.}
A recurring theme across Sections~\ref{sec:challenges_identif_ill_posed}, \ref{sec:challenges_reg_stabiliz}, and~\ref{sec:challenges_adj_grad} is the importance
of preserving the physical and mathematical structure of state-dependent inverse problems. Hybrid and data-driven methods risk undermining this structure if
physical constraints, admissible sets, or adjoint consistency are treated only approximately. Conversely, approaches that embed learning components within
rigorously defined PDE-constrained optimization frameworks, as in~\cite{Akerson2025}, demonstrate that increased flexibility need not come at the expense of
interpretability or deterministic convergence guarantees. Identifying general principles for such structure-preserving integration is a key open problem.
From the perspective adopted in this review, learning-based and operator-learning approaches may be interpreted as flexible, data-driven parameterizations of
the admissible function space for state-dependent constitutive laws. However, such representations do not eliminate the intrinsic identifiability restrictions,
model-class dependence, or state-space coverage limitations emphasized throughout this paper; rather, they inherit these structural constraints through
the governing PDE and the available observations.

\textbf{Open challenges.}
Despite rapid methodological development, a coherent theoretical framework for hybrid and data-driven state-dependent inverse problems is still lacking. Open
questions include how to analyze identifiability and stability when constitutive laws are represented by learned models, how to design regularization and admissible
sets in high-capacity function spaces, and how to combine adjoint-based sensitivities with data-driven surrogates without compromising robustness. Addressing these
challenges will be essential for leveraging hybrid methods in large-scale, time-dependent, and multiphysics applications while retaining the rigor and interpretability
that characterize classical state-dependent inverse problem formulations.

\subsection{Error Analysis and Convergence Guarantees}
\label{sec:challenges_error_converge}

Despite extensive numerical evidence of successful reconstructions in Sections~\ref{sec:apps_heat_2011}--\ref{sec:apps_echem_2015_2026}, rigorous error estimates
and convergence guarantees for algorithms solving state-dependent inverse problems remain scarce. This gap limits predictive confidence, complicates uncertainty
assessment, and hinders principled comparison of competing methodologies, particularly in safety-critical or high-consequence applications.

\textbf{Distinction between algorithmic convergence and reconstruction accuracy.}
A fundamental difficulty in state-dependent inverse problems is that convergence of an optimization algorithm does not necessarily imply accuracy or uniqueness
of the reconstructed constitutive law. Section~\ref{sec:math} establishes that the inverse map is typically non-injective outside the identifiability interval,
so even exact minimization of the reduced objective may yield reconstructions that differ substantially from the true constitutive relation outside the explored
state range. As a result, classical notions of convergence to a ``true'' solution must be interpreted with care. This distinction is visible in applications
such as temperature-dependent heat conduction~\cite{Bukshtynov2011} and thermo-fluid viscosity reconstruction~\cite{Bukshtynov2013}, where iterative optimization
converges reliably, yet the recovered constitutive laws are accurate only on the identifiability interval and regularization dominated elsewhere.

\textbf{Lack of a priori error estimates.}
In contrast to classical coefficient-identification problems, where conditional stability estimates and convergence rates can sometimes be established under
source conditions or sparsity assumptions, comparable a~priori error bounds for state-dependent inverse problems are largely absent. The analytical results
reviewed in Section~\ref{sec:math} focus on existence of derivatives and structural properties of the inverse map, rather than quantitative error estimates.
In particular, the aggregation of sensitivities over state-level sets or trajectories precludes direct application of standard tools from linear inverse theory.
Consequently, most applications rely on empirical convergence diagnostics rather than rigorous bounds linking data noise, discretization error, regularization
strength, and reconstruction accuracy.

\textbf{Interplay between discretization, regularization, and convergence.}
Section~\ref{sec:comp} emphasizes that state-dependent inverse problems are solved in function space but realized numerically through coupled discretizations
of the physical domain, time, and state space. Discretization error interacts nontrivially with regularization and optimization, and mesh refinement does not
necessarily improve reconstruction accuracy in the absence of appropriate regularization. This behavior is observed in multiple applications, including heat
conduction~\cite{Bukshtynov2011}, thermo-fluid systems~\cite{Bukshtynov2013}, and reduced-order closure identification~\cite{ProtasNoack2015}, where
mesh-independent convergence of the optimization algorithm is achieved, yet convergence of the reconstructed constitutive law depends critically on the chosen
Sobolev metric and regularization parameters. Developing error analyses that account simultaneously for discretization, regularization, and state-dependent
identifiability remains an open challenge.

\textbf{Nonconvexity and local convergence.}
As discussed in Sections~\ref{sec:challenges_nonlin_nonconv_apps} and \ref{sec:challenges_identif_ill_posed}, the reduced optimization problems arising in
state-dependent inversion are generally nonconvex. Consequently, most existing algorithms admit, at best, local convergence guarantees to stationary points.
In applications such as turbulent-flow closure reconstruction~\cite{MatharuProtas2022,MatharuProtas2024}, even the existence of minimizers depends on
admissible-set design, making global convergence guarantees unattainable without additional structural assumptions. While quasi-Newton and gradient-based
methods perform well in practice, their convergence behavior is typically assessed empirically, and rigorous analysis remains limited to simplified settings.

\textbf{Role of gradient verification and empirical convergence assessment.}
In the absence of rigorous error theory, careful numerical verification plays an outsized role. Section~\ref{sec:verif_kappa} highlights the use of
$\kappa$-tests and related procedures to verify adjoint correctness and gradient consistency. Such tests are systematically employed in several applications,
including heat conduction and thermo-fluid problems~\cite{Bukshtynov2011,Bukshtynov2013}, inviscid vorticity reconstruction~\cite{DanailaProtas2015},
turbulent-flow closure identification~\cite{MatharuProtas2020,MatharuProtas2022,MatharuProtas2024}, and more recent electrochemical and kinetic-law
inversions~\cite{Ahmadi2026}. While gradient verification is essential for ensuring algorithmic correctness, it does not provide quantitative guarantees
on reconstruction error or predictive performance, underscoring the need for complementary analytical tools.

\textbf{Implications for predictive use and model validation.}
The lack of rigorous error bounds complicates the use of reconstructed constitutive laws in predictive simulations. In electrochemical and turbulent-flow
applications (Sections~\ref{sec:apps_echem_2015_2026} and \ref{sec:apps_MP_2020_2024}), reconstructed state-dependent laws are often interpreted as effective
closures rather than exact physical properties. Without error estimates or convergence guarantees, it is difficult to assess the reliability of predictions
made outside the calibration regime or under perturbed operating conditions. This limitation is particularly significant in applications where inverse results
inform design, control, or safety decisions.

\textbf{Open challenges.}
Developing a rigorous error and convergence theory for state-dependent inverse problems remains a major open challenge. Key questions include how to formulate
meaningful error metrics for function-valued reconstructions, how to relate data noise and discretization error to reconstruction accuracy within identifiability
intervals, and how to characterize convergence behavior in nonconvex, infinite-dimensional optimization problems. Progress in this direction will require new
analytical tools that integrate inverse-problem theory, nonlinear functional analysis, and numerical optimization, and will be essential for elevating
state-dependent inverse methods from powerful computational tools to fully predictive modeling frameworks.

\subsection{Application-Specific Challenges Across Physics and Engineering Domains}
\label{sec:challenges_phys_eng}

While Sections~\ref{sec:math} and~\ref{sec:comp} establish a general mathematical and computational framework for state-dependent inverse problems,
the applications reviewed in Section~\ref{sec:apps} make clear that the practical challenges of inversion are strongly shaped by the underlying physics and
modeling assumptions. Each application domain introduces its own sources of complexity, including multiscale coupling, incomplete physical descriptions,
experiment-model mismatch, and domain-specific constraints on observability and excitation. As a result, inverse reconstruction can recover at best
\emph{effective state-dependent constitutive laws within a prescribed model class}, rather than intrinsic material or transport properties valid across regimes.

\textbf{Model-class dependence and effective constitutive laws.}
From a mathematical perspective, Section~\ref{sec:math} emphasizes that identifiability and well-posedness are properties of the inverse problem \emph{relative
to a chosen forward model}. This dependence becomes particularly evident in application domains where the governing equations represent idealized or reduced
descriptions of complex physical processes. In turbulent-flow closure reconstruction (Section~\ref{sec:apps_MP_2020_2024}), the inferred state-dependent
eddy-viscosity laws depend critically on the LES formulation, the chosen state variable, and the targeted observables~\cite{MatharuProtas2022,MatharuProtas2024}.
Similarly, in electrochemical modeling (Section~\ref{sec:apps_echem_2015_2026}), reconstructed transport and kinetic laws are meaningful only within the specific
continuum battery model used, and may absorb neglected physics such as microstructural heterogeneity or side reactions~\cite{Escalante2020,Daniels2023}. These
observations reinforce a central theme of this review: state-dependent inversion identifies \emph{model-consistent closures}, not model-independent physical laws.

\textbf{Multiscale and multiphysics coupling.}
Many application domains involve strong coupling across disparate spatial and temporal scales, which complicates both analysis and interpretation. In
Navier--Stokes-based inverse problems~\cite{Bukshtynov2013,MatharuProtas2020}, state-dependent parameters influence the dynamics through nonlinear transport
terms, and small-scale effects can propagate intermittently into large-scale observables. In electrochemical systems, transport, reaction kinetics, and electric
fields are tightly coupled across bulk and interfacial processes, often spanning orders of magnitude in time and length scales~\cite{Sethurajan2015,Ahmadi2026}.
These multiscale interactions challenge the assumptions underlying identifiability analysis and error control, and make it difficult to interpret reconstructed
constitutive laws as intrinsic properties rather than regime-dependent effective descriptions.

\textbf{Experiment-model mismatch and observability constraints.}
Application-specific measurement limitations further shape the inverse problem. As discussed in Section~\ref{sec:challenges_observe}, many domains rely on
indirect or aggregated observations, such as voltage curves in batteries or statistical quantities in turbulent flows. In such settings, the inverse problem
may be dominated by model inadequacy rather than data noise. Independent physical analyses in electrochemical modeling~\cite{Richardson2018} demonstrate that
certain reconstructed transport anomalies cannot be reconciled within the adopted model class, signaling missing physics rather than inversion error. Analogous
issues arise in turbulent-flow modeling, where no state-dependent eddy-viscosity closure within a given ansatz can satisfy certain optimization
objectives~\cite{MatharuProtas2022}. These examples highlight that inverse methods cannot compensate for fundamental experiment-model mismatch.

\textbf{Relation to data-driven constitutive modeling.}
The application-specific nature of reconstructed constitutive laws is not unique to inverse problems formulated within PDE-constrained optimization. Recent
reviews of data-driven constitutive modeling in mechanics emphasize that learned or inferred laws are inevitably conditioned on the chosen representation,
training data, and admissible model class, and should be interpreted as effective rather than intrinsic material behavior~\cite{Fuhg2025}. This perspective
aligns closely with the conclusions drawn from the state-dependent inverse problems reviewed here, despite differences in methodology. Together, these bodies
of work point to a broader consensus across physics, engineering, and computational mechanics: constitutive identification is fundamentally model-relative.

\textbf{Implications for transferability and prediction.}
A direct consequence of application-specific challenges is limited transferability of reconstructed constitutive laws across regimes, geometries, or operating
conditions. In reduced-order and turbulent-flow closure identification (Sections~\ref{sec:apps_ROM_closure_2015} and~\ref{sec:apps_MP_2020_2024}), closures
calibrated on a given flow configuration may fail when applied outside the training regime. In electrochemical models, transport and kinetic laws inferred
under specific cycling conditions may not generalize to different temperatures or loading protocols~\cite{Daniels2023}. These limitations underscore the
importance of clearly defining the domain of validity of reconstructed laws and caution against overinterpretation of inverse results as universally predictive.

\textbf{Open challenges.}
Addressing application-specific challenges requires closer integration of inverse-problem methodology with domain knowledge, model development, and experimental
design. Key open questions include how to detect and diagnose model inadequacy within inverse formulations, how to design admissible sets and objectives that
reflect domain-specific physics, and how to quantify the regime of validity of reconstructed state-dependent laws. Progress in these directions will be essential
for translating advances in state-dependent inverse problem theory into reliable and trustworthy tools across physics and engineering applications.

\subsection{Future Directions}
\label{sec:challenges_future}

Looking ahead, several research directions appear particularly promising for advancing the theory and practice of inverse problems with state-dependent parameters.
Importantly, these directions do not represent departures from the variational and adjoint-based framework developed in Sections~\ref{sec:math} and~\ref{sec:comp},
but rather natural extensions motivated by the challenges synthesized throughout this section and illustrated by the applications in Section~\ref{sec:apps}.
\begin{itemize}
  \item \emph{Algorithmic innovation.}
    While adjoint-based optimization is now well established for state-dependent inverse problems, Sections~\ref{sec:challenges_nonlin_nonconv_apps},
    \ref{sec:challenges_adj_grad}, and \ref{sec:challenges_error_converge} make clear that existing algorithms remain vulnerable to nonconvexity, poor globalization,
    and sensitivity to initialization. Future progress will require the development of scalable optimization strategies with improved robustness, including adaptive
    regularization, continuation and homotopy methods, trust-region or filter-based globalization, and structure-aware preconditioning in function space. These needs
    are particularly acute in large-scale and time-dependent applications such as turbulent-flow closure reconstruction and electrochemical modeling
    (Sections~\ref{sec:apps_MP_2020_2024} and~\ref{sec:apps_echem_2015_2026}), where computational cost and landscape complexity remain dominant bottlenecks.
  \item \emph{Extension to broader classes of PDEs.}
    As discussed in Section~\ref{sec:challenges_time_nonlin}, much of the existing analytical theory for state-dependent inverse problems is confined to relatively
    simple elliptic and parabolic settings. Extending both theoretical and computational frameworks to strongly nonlinear, convective, multiphysics, and possibly
    stochastic PDEs remains a major open challenge. Progress in this direction is essential for applications involving turbulent flows, reactive transport, and
    coupled electrochemical systems, where nonlinear dynamics, multiscale coupling, and long-time behavior fundamentally shape identifiability and stability.
    Addressing these problems will likely require new analytical tools that bridge inverse-problem theory with nonlinear dynamics and multiscale analysis.
  \item \emph{Integration with uncertainty quantification and data-driven methods.}
    As emphasized in Section~\ref{sec:non_traditional_OPT} and revisited in Section~\ref{sec:challenges_integr_method}, combining deterministic PDE-constrained
    optimization with uncertainty quantification and carefully designed data-driven components offers a promising path toward robust inversion under realistic
    data limitations. Future work should aim to clarify how probabilistic inference, surrogate modeling, and learned representations can be integrated without
    undermining identifiability analysis, adjoint consistency, or physical interpretability. Developing hybrid frameworks that retain the rigor of variational
    optimization while leveraging data-driven flexibility remains a central methodological opportunity.
  \item \emph{Benchmark problems and standardized datasets.}
    Section~\ref{sec:challenges_bench_reproduct} highlights the lack of standardized benchmark problems as a major impediment to reproducibility and objective
    comparison. Establishing community benchmarks tailored specifically to state-dependent inverse problems (spanning analytically tractable model problems,
    simplified time-dependent systems, and representative multiphysics settings) would enable systematic assessment of algorithms, regularization strategies,
    and adjoint implementations. The dissemination of reference datasets, verified adjoint formulations, and open-source implementations would further accelerate
    cumulative progress across application domains.
  \item \emph{Cross-disciplinary applications.}
    Finally, extending state-dependent inverse problem methodologies to emerging application areas offers opportunities for both impact and theory development.
    Problems in soft matter, biological transport, energy systems, and complex materials share many of the defining features emphasized in this review: nonlinear
    dynamics, state-dependent constitutive behavior, indirect observations, and model uncertainty. Exploring such domains may not only broaden the reach of
    existing methods but also expose new structural challenges that motivate further advances in identifiability theory, regularization, and computational
    optimization.
\end{itemize}

Taken together, these directions suggest that the future of state-dependent inverse problems lies not in replacing established variational and adjoint-based
methodologies, but in extending and integrating them to address increasing model complexity, data limitations, and cross-disciplinary demands. Progress along
these lines will be essential for transforming state-dependent inverse problems from powerful computational tools into reliable, predictive components of
physics-based modeling and simulation.

\section{Conclusions}
\label{sec:conclusion}

This review has surveyed the mathematical foundations, computational methodologies, and representative applications of inverse problems with state-dependent
parameters formulated within PDE-constrained optimization frameworks. Beginning with the seminal variational formulation of Chavent and Lemonnier
in~\cite{ChaventLemonnier1974}, the theory establishes that identifying constitutive laws depending on the system state leads to inverse problems that are
intrinsically nonlinear, infinite dimensional, and structurally ill-posed. Key concepts such as identifiability intervals, state-aggregation structures in
adjoint sensitivities, and the non-injectivity of the parameter--to--observation map provide a unifying mathematical lens through which these problems can
be understood.

Building on this foundation, the review has highlighted how modern computational advances over the past decade have transformed state-dependent parameter
identification from a largely theoretical construct into a viable tool for large-scale, multiphysics applications across physics and engineering. Continuous
``optimize--then--discretize'' formulations, adjoint-based gradient evaluation in function space, Sobolev and admissible-set regularization strategies, and
careful verification of gradients have enabled robust reconstructions in settings ranging from heat transfer and fluid mechanics to electrochemical energy
systems and turbulent-flow modeling. The application studies surveyed in Section~\ref{sec:apps} demonstrate that these methods can be successfully deployed
for time-dependent, nonlinear, and coupled PDE models that were previously out of reach.

At the same time, the synthesis of open problems in Section~\ref{sec:challenges} underscores that many limitations are structural rather than algorithmic.
Nonconvexity induced by state-dependent identifiability, ill-posedness beyond noise amplification, objective-dependent nonexistence of minimizers,
sensitivity to data modality, and model-class dependence of reconstructed constitutive laws all place fundamental bounds on what can be inferred from data.
Regularization, admissible-set design, and experimental excitation can mitigate these issues, but cannot eliminate them. As a result, reconstructed
state-dependent laws should often be interpreted as effective, model-consistent closures rather than intrinsic, universally valid material properties.

Looking forward, progress in this field will likely be driven by a combination of deeper analytical understanding, algorithmic innovation, and carefully
designed hybrid methodologies. Extending theory to strongly nonlinear and time-dependent systems, developing structure-aware optimization algorithms,
integrating uncertainty quantification and data-driven components without sacrificing interpretability, and establishing standardized benchmark problems
all represent promising directions. More broadly, the growing range of applications (from energy systems to complex materials and biological transport)
suggests that state-dependent inverse problems will continue to play a central role in advancing physics-based modeling and predictive simulation.

In this sense, state-dependent parameter identification stands at a transitional point: no longer a niche extension of classical coefficient inversion,
but not yet supported by a complete theoretical and computational infrastructure. By consolidating existing knowledge and clarifying open challenges,
this review aims to provide both a reference framework and a foundation for future developments in PDE-constrained optimization for complex, state-dependent
material systems.

\bibliographystyle{spmpsci}
\bibliography{biblio_Bukshtynov,biblio_Protas,biblio_OPT,biblio_StateDepDiff}

@ARTICLE{ArbicBukshtynov2024,
author = {Arbic II, Paul R. and Bukshtynov, Vladislav},
title = {Efficient Gradient-based Optimization for Reconstructing Binary Images in Applications to Electrical Impedance Tomography},
journal = {Computational Optimization and Applications},
volume = {88},
pages = {379-403},
year  = {2024},
publisher = {Springer}
}

@ARTICLE{ChunEdwardsBukshtynov2024,
author = {Maria M.F.M. Chun and Briana L. Edwards and Vladislav Bukshtynov},
title = {Multiscale Optimization via Enhanced Multilevel {PCA}-based Control Space Reduction for Electrical Impedance Tomography Imaging},
journal = {Computers and Mathematics with Applications},
volume = {157},
pages = {215-234},
year = {2024}
}

@BOOK{BukshtynovBook2023,
author = {Bukshtynov, Vladislav},
title= {Computational Optimization: Success in Practice},
publisher = {Chapman and Hall/CRC},
isbn = {9781032229478},
pages = {416},
url = {https://www.routledge.com/Computational-Optimization/Bukshtynov/p/book/9781032229478},
year = {2023}
}

@ARTICLE{AbdullaBukshtynovSeif2021,
author={Abdulla, Ugur G. and Bukshtynov, Vladislav and Seif, Saleheh},
title={Cancer Detection through Electrical Impedance Tomography and Optimal Control Theory: Theoretical and Computational Analysis},
journal = {Mathematical Biosciences and Engineering},
publisher = {{AIMS} Press},
volume = {18},
number = {4},
pages = {4834-4859},
year = {2021},
month = {jun}
}

@ARTICLE{AbdullaBukshtynovHagverdiyev2019,
title={Gradient Method in {H}ilbert-{B}esov Spaces for the Optimal Control of Parabolic Free Boundary Problems},
author={Abdulla, Ugur G. and Bukshtynov, Vladislav and Hagverdiyev, Ali},
journal={Journal of Computational and Applied Mathematics},
volume={346},
number={},
pages={84-109},
year={2019}
}

@ARTICLE{VolkovBukshtynov2018,
author={Volkov, O. and Bukshtynov, V. and Durlofsky, L.J. and Aziz, K.},
title={Gradient-based {P}areto Optimal History Matching for Noisy Data of Multiple Types},
journal={Computational Geosciences},
year={2018},
volume={22},
number={6},
pages={1465-1485}
}

@ARTICLE{Bukshtynov2015,
author={Bukshtynov, V. and Volkov, O. and Durlofsky, L.J. and Aziz, K.},
title={Comprehensive Framework for Gradient-based Optimization in Closed-Loop Reservoir Management},
journal={Computational Geosciences},
year={2015},
volume={19},
number={4},
pages={877-897}
}

@ARTICLE{Bukshtynov2013,
author={Bukshtynov, V. and Protas, B.},
title={Optimal Reconstruction of Material Properties in Complex Multiphysics Phenomena},
journal={Journal of Computational Physics},
year={2013},
volume={242},
pages={889-914}
}

@BOOK{BukshtynovPhD2012,
author = {Bukshtynov, V.},
title= {Computational Methods for the Optimal Reconstruction of Material Properties in Complex Multiphysics Systems},
publisher = {Ph.D.~Dissertation, McMaster University, Open Access Dissertation and Theses: Paper 6795},
url = {http://hdl.handle.net/11375/11859},
year = {2012}
}

@ARTICLE{Bukshtynov2011,
author={Bukshtynov, V. and Volkov, O. and Protas, B.},
title={On Optimal Reconstruction of Constitutive Relations},
journal={Physica D: Nonlinear Phenomena},
year={2011},
volume={240},
number={16},
pages={1228-1244}
}

@BOOK{Luenberger1976,
  title={Optimization by Vector Space Methods},
  author={Luenberger, D.G.},
  series={Series in Decision and Control},
  year={1976},
  publisher={Wiley}
}

@BOOK{Nocedal2006,
 edition = {2nd},
 author = {J. Nocedal and S. J. Wright},
 title = {Numerical Optimization},
 publisher = {Springer},
 year = {2006}
}

@BOOK{Berger1977,
author={Berger, M. S},
title={Nonlinearity and Functional Analysis},
publisher = {Acad.~Press, New York},
year = {1977}
}

@BOOK{Gunzburger2003,
title={Perspectives in Flow Control and Optimization},
author={Gunzburger, Max D},
volume={5},
year={2003},
publisher={SIAM}
}

@BOOK{BoydVandenberghe2004,
title={Convex Optimization},
author={Boyd, S. and Vandenberghe, L.},
year={2004},
publisher={Cambridge University Press}
}

@BOOK{Ruszczynski2006,
title={Nonlinear Optimization},
author={Ruszczy\'nski, A.},
year={2006},
publisher={Princeton University Press}
}

@BOOK{Vogel2002,
title={Computational Methods for Inverse Problems},
author={Vogel, C. R.},
year={2006},
publisher={SIAM}
}

@BOOK{Tarantola2005,
author={Tarantola, A.},
title={Inverse Problem Theory and Methods for Model Parameter Estimation},
year={2005},
publisher={SIAM}
}

@BOOK{EnglHankeNeubauer1996,
author={Engl, H. and Hanke, M. and Neubauer, A.},
title={Regularization of Inverse Problems},
year={1996},
publisher={Kluwer}
}

@ARTICLE{Leeuwen2016,
author = {van Leeuwen, T. and Herrmann, F. J.},
title = {A Penalty Method for {PDE}-Constrained Optimization in Inverse Problems},
journal = {Inverse Problems},
publisher = {IOP Publishing},
year = {2016},
volume = {32},
number = {1},
pages = {015007}
}

@ARTICLE{MinGibou2008,
author={Min, C. and Gibou, F.},
title={Robust Second-Order Accurate Discretizations of the Multi-Dimensional {H}eaviside and {D}irac Delta Functions},
journal={Journal of Computational Physics},
year={2008},
volume={227},
number={22},
pages={9686-9695}
}

@ARTICLE{MinGibou2007,
author={Min, C. and Gibou, F.},
title={Geometric Integration over Irregular Domains with Application to Level-Set Methods},
journal={Journal of Computational Physics},
year={2007},
volume={226},
number={2},
pages={1432-1443}
}

@ARTICLE{Mayo1984,
author={Mayo, A.},
title={The Fast Solution of {P}oisson's and the Biharmonic Equations on Irregular Regions},
journal={SIAM Journal of Numerical Analysis},
year={1984},
volume={21},
pages={285-299}
}

@ARTICLE{ZahediTornberg2010,
author={Zahedi, S. and Tornberg, A.-K.},
title={Delta Function Approximations in Level Set Methods by Distance Function Extension},
journal={Journal of Computational Physics},
year={2010},
volume={229},
number={6},
pages={2199-2219}
}

@ARTICLE{Smereka2006,
author={Smereka, P.},
title={The Numerical Approximation of a Delta Function with Application to Level Set Methods},
journal={Journal of Computational Physics},
year={2006},
volume={211},
number={1},
pages={77-90}
}

@ARTICLE{Beale2008,
author={Beale, J.T.},
title={A Proof that a Discrete Delta Function is Second-Order Accurate},
journal={Journal of Computational Physics},
year={2008},
volume={227},
number={4},
pages={2195-2197}
}

@ARTICLE{Towers2009,
author={Towers, J.D.},
title={Discretizing Delta Functions via Finite Differences and Gradient Normalization},
journal={Journal of Computational Physics},
year={2009},
volume={228},
number={10},
pages={3816-3836}
}

@ARTICLE{Towers2007,
author={Towers, J.D.},
title={Two Methods for Discretizing a Delta Function Supported on a Level Set},
journal={Journal of Computational Physics},
year={2007},
volume={220},
number={2},
pages={915-931}
}

@BOOK{OliverReynoldsLiu2008,
title={Inverse Theory for Petroleum Reservoir Characterization and History Matching},
author={Oliver, Dean S. and Reynolds, Albert C and Liu, Ning},
year={2008},
publisher={Cambridge University Press}
}

@ARTICLE{ChaventLemonnier1974,
author = {Chavent, G. and Lemonnier, P.}, 
title = {Identification de la Non--Linearit\'e D'Une \'Equation Parabolique Quasilineaire},
volume = {1}, 
pages = {121-162}, 
year = {1974}, 
journal = {Applied Mathematics and Optimization}
}

@BOOK{BanksKunisch1989,
author = {Banks, H. and Kunisch, K.}, 
title = {Estimation Techniques for Distributed Parameter Systems},
publisher = {Birkh\"auser}, 
year = {1989}
}

@BOOK{NashedScherzer2002,
author = {Nashed, M. Z. and Scherzer, O.}, 
title = {Inverse Problems, Image Analysis, and Medical Imaging},
publisher = {AMS}, 
year = {2002}
}

@BOOK{GottliebDuChateau1996,
author = {Gottlieb, J. and DuCh{\^a}teau, P.}, 
title = {Parameter Identification and Inverse Problems in Hydrology, Geology, and Ecology},
publisher = {Kluwer Academic Publishers}, 
year = {1996}
}

@ARTICLE{Alifanov2004,
author = {Alifanov, O. M. and Nenarokomov, A. V. and Budnik, S. A. and Michailov, V. V. and Ydin, V. M.}, 
title = {Identification of Thermal Properties of Materials with Applications for Spacecraft Structures},
volume = {12}, 
pages = {579-594}, 
year = {2004}, 
journal = {Inverse Problems in Science and Engineering}
}

@ARTICLE{Alifanov2007,
author = {Alifanov, O. M. and Budnik, S. A. and Michaylov, V. V. and Nenarokomov, A. V. and Titov, D. M. and Ydin, V. M.}, 
title = {An Experimental--Computational System for Materials Thermal Properties Determination and its Application for
Spacecraft Structures Testing},
volume = {61}, 
pages = {341-451}, 
year = {2007}, 
journal = {Acta Astronautica}
}

@ARTICLE{Tai1995,
author = {Tai, X.-C. and K\"{a}rkk\"{a}inen, T.},
title = {Identification of a Nonlinear Parameter in a Parabolic Equation from a Linear Equation},
journal = {Computational and Applied Mathematics},
volume = {14},
year = {1995},
pages = {157-184}
}

@ARTICLE{DuChateau2004,
author = {DuCh{\^a}teau, P. and Thelwell, R. and Butters, G.},
title = {Problem Approach to the Identification of an Unknown Diffusion Coefficient},
journal = {Inverse Problems},
volume = {20},
year = {2004},
pages = {601-625}
}

@ARTICLE{JanickiKindermann2009,
author = {Janicki, M. and Kindermann, S.},
title = {Recovering Temperature Dependence of Heat Transfer Coefficient in Electronic Circuits},
journal = {Inverse Problems in Science and Engineering},
volume = {17},
year = {2009},
pages = {1129-1142}
}

@ARTICLE{HankeScherzer1999,
author = {Hanke, M. and Scherzer, O.},
title = {Error Analysis of an Equation Error Method for the Identification of the Diffusion Coefficient in a Quasi-Linear
Parabolic Differential Equation},
journal = {SIAM J. Appl. Math.},
volume = {59},
year = {1999},
pages = {1012-1027}
}

@ARTICLE{Luo2003,
author = {Luo, D. and He, L. and Lin, S. and Chen, T.--F. and Gao, D.},
title = {Determination of Temperature Dependent Thermal Conductivity by Solving {IHCP} in Infinite Regions},
journal = {Int.~Comm.~Heat Mass Transfer},
volume = {30},
year = {2003},
pages = {903-908}
}

@ARTICLE{DuChateau1997,
author = {DuCh{\^a}teau, P.},
title = {An Inverse Problem for the Hydraulic Properties of Porous Media},
journal = {SIAM J. Math. Anal.},
volume = {28},
year = {1997},
pages = {611-632}
}

@BOOK{Kugler2000,
author = {K\"ugler, Ph.},
title = {Identification of a Temperature Dependent Heat Conductivity by {T}ikhonov Regularization},
publisher = {Diploma Thesis, Johannes Kepler Universit\"at Linz},
year = {2000}
}

@ARTICLE{Kugler2003,
author = {K\"ugler, Ph.},
title = {Identification of a Temperature Dependent Heat Conductivity from Single Boundary Measurements},
journal = {SIAM J.~Numer.~Anal.},
volume = {41},
year = {2003},
pages = {1543-1563}
}

@ARTICLE{Neubauer2008,
author = {Neubauer, A.},
title = {Identification of a Temperature Dependent Heat Conductivity via Adaptive Grid Regularization},
journal = {Journal of Integral Equations and Applications},
volume = {20},
year = {2008},
pages = {229-242}
}

@BOOK{Isakov2017,
author = {Isakov, V.},
title = {Inverse Problems for Partial Differential Equations},
publisher = {Springer},
address   = {New York},
edition   = {3rd},
year = {2017}
}

@ARTICLE{ColemanNoll1963,
author = {Coleman, B. D. and Noll, W.},
title = {The Thermodynamics of Elastic Materials with Heat Conduction and Viscosity},
journal = {Arch.~Rat.~Mech.~Anal},
volume = {13},
year = {1963},
pages = {167-178}
}

@ARTICLE{Liu1972,
author = {Liu, I.-S.},
title = {Method of {L}agrange Multipliers for Exploitation of the Entropy Principle},
journal = {Arch.~Rat.~Mech.~Anal},
volume = {46},
year = {1972},
pages = {131-148}
}

@BOOK{TruesdellNoll2004,
author = {Truesdell, Clifford and Noll, Walter},
title = {The Non-Linear Field Theories of Mechanics},
series = {Springer Tracts in Natural Philosophy},
volume = {3},
edition = {3rd},
publisher = {Springer},
address = {Berlin, Heidelberg},
year = {2004}
}

@BOOK{TikhonovArsenin1977,
author = {Tikhonov, Andrei N. and Arsenin, Vasiliy Y.},
title = {Solutions of Ill-Posed Problems},
publisher = {Winston},
address = {Washington, D.C.},
year = {1977}
}

@BOOK{Biegler2007,
author = {Biegler, L. T. and Ghattas, O. and Heinkenschloss, M. and Keyes, D. and van Bloemen Waanders, B.},
title = {Real-Time {PDE}-Constrained Optimization},
publisher = {SIAM},
year = {2007}
}

@BOOK{Biegler2010,
author = {Biegler, Lorenz T.},
title = {Nonlinear Programming: Concepts, Algorithms, and Applications to Chemical Processes},
publisher = {Society for Industrial and Applied Mathematics (SIAM)},
address = {Philadelphia},
year = {2010}
}

@BOOK{Biegler2003conf,
author = {Biegler, L. T. and Heinkenschloss, M. and Ghattas, O. and Keyes, D. and van Bloemen Waanders, B.},
title = {Large-Scale PDE-Constrained Optimization},
publisher = {Springer},
year = {2003}
}

@BOOK{Muschik1989,
author = {Muschik, W.},
title = {Aspects of Non-Equilibrium Thermodynamics},
publisher = {World Scientific},
year = {1989}
}

@ARTICLE{Triano2008,
author = {Triano, V. and Papenfuss, Ch. and Cimmelli, V. A. and Muschik, W.},
title = {Exploitation of the {S}econd {L}aw: {C}oleman--{N}oll and {L}iu Procedure in Comparison},
journal = {J.~Non-Equilib.~Thermodyn.},
volume = {33},
year = {2008},
pages = {47-60}
}

@BOOK{Lions1968,
author = {Lions, J.-L.},
title = {Contr\^ole Optimal de Syst\`emes Gouvern\'es par des \'Equations aux D\'eriv\'ees Partielles},
publisher = {Dunod},
year = {1968}
}

@BOOK{Lions1969,
author = {Lions, J.-L.},
title = {Quelques M\'ethodes de R\'esolution des Probl\`emes aux Limites Non Lin\'eaires},
publisher = {Dunod--Gauthier-Villars},
year = {1969}
}

@BOOK{Lions1971,
author = {Lions, Jacques-Louis},
title = {Optimal Control of Systems Governed by Partial Differential Equations},
series = {Die Grundlehren der mathematischen Wissenschaften},
volume = {170},
publisher = {Springer-Verlag},
address = {Berlin--Heidelberg--New York},
year = {1971}
}

@BOOK{GlowinskiLionsHe2008,
author = {Glowinski, Roland and Lions, Jacques-Louis and He, Jiwen},
title = {Exact and Approximate Controllability for Distributed Parameter Systems: A Numerical Approach},
series = {Encyclopedia of Mathematics and its Applications},
volume = {117},
publisher = {Cambridge University Press},
address = {Cambridge},
year = {2008}
}

@BOOK{GlowinskiLionsTremolieres1981,
author = {Glowinski, Roland and Lions, Jacques-Louis and Tr{\'e}molieres, Raymond},
title = {Numerical Analysis of Variational Inequalities},
series = {Studies in Mathematics and Its Applications},
publisher = {North-Holland Publishing Company},
address = {Amsterdam--New York--Oxford},
year= {1981}
}

@ARTICLE{JosephLundgren1972,
author = {Joseph, D. D. and Lundgren, T. S.},
title = {Quasilinear {D}irichlet Problems Driven by Positive Sources},
journal = {Archive for Rational Mechanics and Analysis},
volume = {49},
year = {1972},
pages = {241-269}
}

@BOOK{Vazquez2007,
author = {V\'azquez, J. L.},
title = {The Porous Medium Equation: Mathematical Theory},
publisher = {Oxford University Press},
year = {2007}
}

@BOOK{Cea1971,
author = {Jean C\'ea},
title = {Optimisation: Th\'eorie et Algorithmes},
series = {M\'ethodes Math\'ematiques de l'Informatique},
publisher = {Dunod},
address = {Paris, France},
year = {1971}
}

@BOOK{Hinze2009,
author = {Michael Hinze and Rene Pinnau and Michael Ulbrich and Stefan Ulbrich},
title = {Optimization with {PDE} Constraints},
series = {Mathematical Modelling: Theory and Applications},
volume = {23},
publisher = {Springer},
address = {Dordrecht, The Netherlands},
year = {2009}
}

@BOOK{Troltzsch2010,
author = {Fredi Tr\"{o}ltzsch},
title = {Optimal Control of Partial Differential Equations: Theory, Methods and Applications},
series = {Graduate Studies in Mathematics},
volume = {112},
publisher = {American Mathematical Society},
address = {Providence, RI, USA},
year = {2010}
}

@ARTICLE{Egger2014,
author= {Herbert Egger and Jan-Frederik Pietschmann and Matthias Schlottbom},
title = {Identification of Nonlinear Heat Conduction Laws},
journal= {arXiv preprint},
volume= {arXiv:1404.2535},
year= {2014},
url= {https://arxiv.org/abs/1404.2535}
}

@ARTICLE{Egger2017,
author = {Egger, Herbert and Pietschmann, Jan-Frederik and Schlottbom, Matthias},
title = {On the Uniqueness of Nonlinear Diffusion Coefficients in the Presence of Lower Order Terms},
journal = {Inverse Problems},
year = {2017},
volume = {33},
number = {11},
pages = {115005}
}

@ARTICLE{BurgerOsher2004,
author = {Martin Burger and Stanley J. Osher},
title = {Convergence Rates of Convex Variational Regularization},
journal = {Inverse Problems},
volume = {20},
number = {5},
pages = {1411-1421},
year = {2004}
}

@BOOK{BorziSchulz2011,
author = {Alfio Borz{\`i} and Volker H. Schulz},
title = {Computational Optimization of Systems Governed by Partial Differential Equations},
publisher = {SIAM},
address = {Philadelphia, PA, USA},
year = {2011}
}

@ARTICLE{CannonDuChateau1973,
author = {Cannon, J. R. and DuCh{\^a}teau, P.},
title = {Determining Unknown Coefficients in a Nonlinear Heat Conduction Problem},
journal = {SIAM Journal on Applied Mathematics},
volume = {24},
number = {3},
year = {1973},
pages = {298-314}
}

@ARTICLE{CannonDuChateau1980a,
author = {Cannon, J. R. and DuCh{\^a}teau, P.},
title = {An Inverse Problem for a Nonlinear Diffusion Equation},
journal = {SIAM J. Appl. Math.},
volume = {39},
number = {2},
year = {1980},
pages = {272-289}
}

@ARTICLE{CannonDuChateau1980b,
author = {Cannon, J. R. and DuCh{\^a}teau, P.},
title = {An Inverse Problem for an Unknown Source in a Heat Equation},
journal = {Journal of Mathematical Analysis and Applications},
volume = {75},
number = {2},
year = {1980},
pages = {465-485}
}

@ARTICLE{CannonDuChateau1983,
author = {Cannon, J. R. and DuCh{\^a}teau, P.},
title = {An Inverse Problem for an Unknown Source Term in a Wave Equation},
journal = {SIAM J. Appl. Math.},
volume = {43},
number = {3},
year = {1983},
pages = {553-564}
}

@ARTICLE{CannonDuChateau1987,
author = {Cannon, J. R. and DuCh{\^a}teau, P.},
title = {Design of an Experiment for the Determination of an Unknown Coefficient in a Nonlinear Conduction-Diffusion Equation},
journal = {International Journal of Engineering Science},
volume = {25},
number = {8},
year = {1987},
pages = {1067-1078}
}

@ARTICLE{KohnVogelius1984,
author = {R. V. Kohn and M. Vogelius},
title = {Determining Conductivity by Boundary Measurements},
journal = {Communications on Pure and Applied Mathematics},
volume = {37},
number = {3},
pages = {289-298},
year = {1984}
}

@ARTICLE{KohnVogelius1987,
author = {R. V. Kohn and M. Vogelius},
title = {Relaxation of a Variational Method for Impedance Computed Tomography},
journal = {Communications on Pure and Applied Mathematics},
volume = {40},
number = {6},
pages = {745-777},
year = {1987}
}

@BOOK{Prilepko2000,
author = {Aleksey I. Prilepko and Dmitry G. Orlovsky and Igor A. Vasin},
title = {Methods for Solving Inverse Problems in Mathematical Physics},
series = {Monographs and Textbooks in Pure and Applied Mathematics},
volume = {231},
publisher = {Marcel Dekker, Inc.},
address = {New York, NY, USA},
year = {2000}
}

@ARTICLE{ChoulliYamamoto1997,
author = {Choulli, M. and Yamamoto, M.},
title = {An Inverse Parabolic Problem with Non-Zero Initial Condition},
journal = {Inverse Problems},
volume = {13},
number = {1},
pages = {19-27},
year = {1997}
}

@ARTICLE{ChoulliYamamoto1996,
author  = {Choulli, M. and Yamamoto, M.},
title = {Generic Well-Posedness of a Linear Inverse Parabolic Problem -- the {H\"older}-Space Approach},
journal = {Inverse Problems},
volume  = {12},
number  = {3},
pages   = {195-205},
year    = {1996}
}

@ARTICLE{EnglKunischNeubauer1989,
author = {Heinz W. Engl and Karl Kunisch and Andreas Neubauer},
title = {Convergence Rates for {T}ikhonov Regularization of Nonlinear Ill-Posed Problems},
journal = {Inverse Problems},
volume = {5},
number = {4},
pages = {523-540},
year = {1989}
}

@BOOK{Plato1995,
author = {Robert Plato},
title = {Iterative and Other Methods for Linear and Ill-Posed Equations},
publisher = {Technische Universit\"at Berlin},
year = {1995},
address = {Berlin, Germany}
}

@BOOK{DelfourZolesio2011,
author = {Michel C. Delfour and Jean-Paul Zol\'esio},
title = {Shapes and Geometries: Metrics, Analysis, Differential Calculus, and Optimization},
series = {Advances in Design and Control},
edition = {2nd},
publisher = {SIAM},
address = {Philadelphia, PA},
year = {2011}
}

@BOOK{SokolowskiZolesio1992,
author = {Jan Soko{\l}owski and Jean-Paul Zol\'esio},
title = {Introduction to Shape Optimization: Shape Sensitivity Analysis},
publisher = {Springer},
series = {Springer Series in Computational Mathematics},
volume = {16},
address = {Berlin},
year = {1992}
}

@BOOK{HaslingerNeittaanmaki1996,
author = {Jaroslav Haslinger and Pekka Neittaanm\"aki},
title = {Finite Element Approximation for Optimal Shape, Material and Topology Design},
publisher = {John Wiley \& Sons},
address = {Chichester},
edition = {2nd},
year = {1996}
}

@ARTICLE{HinzeTroltzsch2010,
author = {Michael Hinze and Fredi Tr\"oltzsch},
title = {Discrete Concepts Versus Error Analysis in {PDE}-Constrained Optimization},
journal = {GAMM-Mitteilungen},
publisher = {Wiley},
volume = {33},
number = {2},
pages = {148-162},
year = {2010}
}

@ARTICLE{Wilcox2015,
author = {Lucas C. Wilcox and Georg Stadler and Tan Bui-Thanh and Omar Ghattas},
title = {Discretely Exact Derivatives for Hyperbolic {PDE}-Constrained Optimization Problems Discretized by the Discontinuous {G}alerkin Method},
journal = {Journal of Scientific Computing},
publisher = {Springer},
volume = {63},
pages = {138-162},
year = {2015}
}

@BOOK{AschBocquetNodet2016,
author = {Asch, Mark and Bocquet, Marc and Nodet, Ma{\"e}lle},
title = {Data Assimilation: Methods, Algorithms, and Applications},
publisher = {Society for Industrial and Applied Mathematics (SIAM)},
address = {Philadelphia, PA},
year = {2016},
series = {Fundamentals of Algorithms}
}

@BOOK{ItoKunisch2008,
author = {Ito, Kazufumi and Kunisch, Karl},
title = {{L}agrange Multiplier Approach to Variational Problems and Applications},
publisher = {Society for Industrial and Applied Mathematics (SIAM)},
address = {Philadelphia, PA},
year = {2008},
series = {Advances in Design and Control},
volume = {15}
}

@BOOK{Neuberger1997,
author = {Neuberger, J.},
title = {Sobolev Gradients and Differential Equations},
publisher = {Springer},
year = {1997}
}

@BOOK{Bertsekas2016,
author = {Bertsekas, Dimitri P.},
title = {Nonlinear Programming},
edition = {3rd},
publisher = {Athena Scientific},
address = {Belmont, MA},
year = {2016}
}

@BOOK{Grisvard1985,
author = {Grisvard, P.},
title = {Elliptic Problems in Nonsmooth Domains},
publisher = {Pitman Publishing},
year = {1985}
}

@ARTICLE{Engquist2004,
author = {Engquist, B. and Tornberg, A.-K. and Tsai, R.},
title = {Discretization of {D}irac Delta Functions in Level Set Methods},
journal = {Journal of Computational Physics},
volume = {207},
pages = {28-51},
year = {2004}
}

@ARTICLE{Huntul2021a,
author = {Huntul, M. J.},
title = {Simultaneous Reconstruction of Time-Dependent Coefficients in the Parabolic Equation from Over-Specification Conditions},
journal = {Results in Applied Mathematics},
volume = {12},
year = {2021},
pages = {100197}
}

@ARTICLE{Hazanee2015,
author = {Hazanee, A. and Lesnic, D. and Ismailov, M. I. and Kerimov, N. B.},
title = {An Inverse Time-Dependent Source Problem for the Heat Equation with a Non-Classical Boundary Condition},
journal = {Applied Mathematical Modelling},
volume = {39},
number = {20},
year = {2015},
pages = {6258-6272}
}

@ARTICLE{Hazanee2019,
author = {Hazanee, A. and Lesnic, D. and Ismailov, M. I. and Kerimov, N. B.},
title = {Inverse Time-Dependent Source Problems for the Heat Equation with Nonlocal Boundary Conditions},
journal = {Applied Mathematics and Computation},
volume = {346},
year = {2019},
pages = {800-815}
}

@ARTICLE{Kanca2013,
author = {Kanca, F.},
title = {The Inverse Problem of the Heat Equation with Periodic Boundary and Integral Overdetermination Conditions},
journal = {Journal of Inequalities and Applications},
volume = {2013},
year = {2013},
pages = {108}
}

@ARTICLE{Huntul2021b,
author = {Huntul, M. J. and Abbas, M. and Baleanu, D.},
title = {An Inverse Problem of Reconstructing the Time-Dependent Coefficient in a One-Dimensional Hyperbolic Equation},
journal = {Advances in Difference Equations},
year = {2021},
volume = {2021},
pages = {452}
}

@ARTICLE{TroltzschYousept2012,
author = {Tr\"oltzsch, F. and Yousept, I.},
title = {{PDE}-Constrained Optimization of Time-Dependent {3D} Electromagnetic Induction Heating by Alternating Voltages},
journal = {ESAIM: Mathematical Modelling and Numerical Analysis},
volume = {46},
number = {4},
pages = {709-729},
year = {2012}
}

@ARTICLE{BarkerStoll2015,
author= {Barker, A. T. and Stoll, M.},
title= {Domain Decomposition in Time for {PDE}-Constrained Optimization},
journal= {Computer Physics Communications},
volume= {197},
pages= {136-143},
year = {2015}
}

@ARTICLE{StollBreiten2015,
author = {Stoll, M. and Breiten, T.},
title = {A Low‐Rank in Time Approach to {PDE}‐Constrained Optimization},
journal = {SIAM Journal on Scientific Computing},
volume = {37},
number = {1},
pages = {B1-B29},
year = {2015}
}

@ARTICLE{Nguyen2019,
author = {Nguyen, T. T. N},
title = {{L}andweber--{K}aczmarz for Parameter Identification in Time‑Dependent Inverse Problems: All‑at‑Once Versus Reduced Version},
journal = {Inverse Problems},
volume = {35},
number = {3},
pages = {035009},
year = {2019}
}

@ARTICLE{Iglesias2016,
author = {Iglesias, M. A.},
title = {A Regularizing Iterative Ensemble {K}alman Method for {PDE}-Constrained Inverse Problems},
journal = {Inverse Problems},
volume = {32},
number = {2},
pages = {025002},
year = {2016}
}

@ARTICLE{Iglesias2013,
author = {Iglesias, M. A. and Law, K. J. H. and Stuart, A. M.},
title = {The Ensemble {K}alman Filter for Inverse Problems},
journal = {Inverse Problems},
volume = {29},
number = {4},
pages = {045001},
year = {2013}
}

@ARTICLE{Stuart2010,
author = {Stuart, A. M.},
title = {Inverse Problems: A {B}ayesian Perspective},
journal = {Acta Numerica},
volume = {19},
pages = {451-559},
year = {2010}
}

@ARTICLE{DashtiStuart2017,
author = {Dashti, M. and Stuart, A. M.},
title = {The {B}ayesian Approach to Inverse Problems},
journal = {Handbook of Uncertainty Quantification},
pages = {311-428},
year = {2017},
publisher = {Springer International Publishing}
}

@BOOK{Evensen2009,
author = {Evensen, Geir},
title = {Data Assimilation: The Ensemble {K}alman Filter},
publisher = {Springer},
year = {2009},
series = {Springer Verlag},
address = {Berlin, Heidelberg}
}

@ARTICLE{HuntulTekin2023,
author = {Huntul, M. J. and Tekin, I.},
title = {An Inverse Problem of Identifying the Time-Dependent Potential and Source Terms in a Two-Dimensional Parabolic Equation},
journal = {Hacettepe Journal of Mathematics and Statistics},
volume = {52},
number = {6},
pages = {1578-1599},
year = {2023}
}

@ARTICLE{IbraheemHuntul2025,
author = {Ibraheem, Q. W. and Huntul, M. J. and Hussein, M. S.},
title = {Time-Dependent Term Identification in the Time-Space Fractional Derivative Diffusion Equation From Integral Over Specified Condition},
journal = {Mathematical Methods in the Applied Sciences},
volume = {n/a},
number = {n/a},
pages = {1-23},
year = {2025}
}

@ARTICLE{Krishnan2021,
author = {Krishnan, V. P. and Rakesh, R. and Senapati, S.},
title = {Stability for a Formally Determined Inverse Problem for a Hyperbolic {PDE} with Space and Time Dependent Coefficients},
journal = {SIAM Journal on Mathematical Analysis},
volume = {53},
number = {6},
pages = {6822-6846},
year = {2021}
}

@ARTICLE{BeranekReinholdUrban2023,
author = {Beranek, N. and Reinhold, M. A. and Urban, K.},
title = {A Space-Time Variational Method for Optimal Control Problems: Well‑Posedness, Stability and Numerical Solution},
journal = {Computational Optimization and Applications},
volume = {86},
pages = {767-794},
year = {2023}
}

@ARTICLE{Sunseri2020,
author = {Sunseri, I. and Hart, J. and van Bloemen Waanders, B. and Alexanderian, A.},
title = {Hyper-Differential Sensitivity Analysis for Inverse Problems Constrained by Partial Differential Equations},
journal = {Inverse Problems},
year = {2020},
publisher = {IOP Publishing},
volume = {36},
number = {12},
pages = {125001}
}

@ARTICLE{Amini2022,
author = {Amini, D. and Haghighat, E. and Juanes, R.},
title = {Physics-Informed Neural Network Solution of Thermo-Hydro-Mechanical Processes in Porous Media},
journal = {Journal of Engineering Mechanics},
publisher = {American Society of Civil Engineers (ASCE)},
volume = {148},
number = {11},
year = {2022},
month = {nov}
}

@ARTICLE{Amini2023,
author = {Amini, D. and Haghighat, E. and Juanes, R.},
title = {Inverse Modeling of Nonisothermal Multiphase Poromechanics Using Physics-Informed Neural Networks},
journal = {Journal of Computational Physics},
volume = {490},
pages = {112323},
year = {2023}
}

@ARTICLE{Yin2023,
author = {Yin, Z. and Orozco, R. and Louboutin, M. and Herrmann, F. J.},
title = {Solving Multiphysics-based Inverse Problems with Learned Surrogates and Constraints},
journal = {Advanced Modeling and Simulation in Engineering Sciences},
volume = {10},
pages = {14},
year = {2023}
}

@ARTICLE{Yamamoto2009,
author = {Yamamoto, M.},
title = {Carleman Estimates for Parabolic Equations and Applications},
journal = {Inverse Problems},
year = {2009},
volume = {25},
number = {12},
pages = {123013}
}

@ARTICLE{BarajasTartakovsky2019,
author = {Barajas-Solano, D. A. and Tartakovsky, A. M.},
title = {Approximate {B}ayesian Model Inversion for {PDE}s with Heterogeneous and State-Dependent Coefficients},
journal = {Journal of Computational Physics},
volume = {395},
pages = {247-262},
year = {2019}
}

@ARTICLE{Lan2023,
author = {Lan, S. and Li, S. and Pasha, M.},
title = {Bayesian Spatiotemporal Modeling for Inverse Problems},
journal = {Statistics and Computing},
year = {2023},
volume = {33},
number = {4},
pages = {89}
}

@ARTICLE{Bingham2024,
author  = {Bingham, D. and Butler, T. and Estep, D.},
title   = {Inverse Problems for Physics-Based Process Models},
journal = {Annual Review of Statistics and Its Application},
volume  = {11},
pages   = {461-482},
year    = {2024}
}

@ARTICLE{Reiser2025,
author = {Reiser, R. and Aguilar, J. E. and Guthke, A. and B\"urkner, P.-C.},
title = {Uncertainty Quantification and Propagation in Surrogate-based {B}ayesian Inference},
journal = {Statistics and Computing},
volume = {35},
pages = {66},
year = {2025}
}

@ARTICLE{Bai2024,
  author  = {Bai, T. and Teckentrup, A. L. and Zygalakis, K. C.},
  title   = {Gaussian Processes for {B}ayesian Inverse Problems Associated with Linear Partial Differential Equations},
  journal = {Statistics and Computing},
  volume  = {34},
  pages   = {139},
  year    = {2024}
}

@ARTICLE{Aarset2023,
author = {Aarset, C. and Holler, M. and Nguyen, T. T. N.},
title = {Learning-Informed Parameter Identification in Nonlinear Time-Dependent {PDE}s},
journal = { Applied Mathematics \& Optimization },
year = {2023},
volume = {88},
pages = {76}
}

@ARTICLE{Raissi2019,
author= {Raissi, M. and Perdikaris, P. and Karniadakis, G. E.},
title = {Physics-Informed Neural Networks: A Deep Learning Framework for Solving Forward and Inverse Problems Involving Nonlinear Partial Differential Equations},
journal = {Journal of Computational Physics},
volume = {378},
pages = {686-707},
year = {2019}
}

@ARTICLE{Raissi2020,
author = {Raissi, M. and Yazdani, A. and Karniadakis, G. E.},
title = {Hidden Fluid Mechanics: Learning Velocity and Pressure Fields from Flow Visualizations},
journal = {Science},
volume = {367},
number = {6481},
pages = {1026-1030},
year = {2020}
}

@ARTICLE{Lu2021,
author = {Lu, L. and Jin, P. and Pang, G. and Zhang, Z. and Karniadakis, G. E.},
title = {Learning Nonlinear Operators via {DeepONet} Based on the Universal Approximation Theorem of Operators},
journal = {Nature Machine Intelligence},
volume = {3},
pages = {218-229},
year = {2021}
}

@ARTICLE{Kovachki2023,
author  = {Kovachki, N. and Li, Z. and Liu, B. and Azizzadenesheli, K. and Bhattacharya, K. and Stuart, A. and Anandkumar, A.},
title   = {Neural Operator: Learning Maps Between Function Spaces With Applications to {PDE}s},
journal = {Journal of Machine Learning Research},
volume  = {24},
pages   = {1-97},
year    = {2023}
}

@ARTICLE{Karniadakis2021,
author = {Karniadakis, G. E. and Kevrekidis, I. G. and Lu, L. and Perdikaris, P. and Wang, S. and Yang, L.},
title = {Physics‐Informed Machine Learning},
journal = {Nature Reviews Physics},
volume = {3},
pages = {422-440},
year = {2021}
}

@ARTICLE{Luo2025,
author = {Luo, K. and Zhao, J. and Wang, Y. and Li, J. and Wen, J. and Liang, J. and Soekmadji, H. and Liao, S.},
title = {Physics-Informed Neural Networks for {PDE} Problems: a Comprehensive Review},
journal = {Artificial Intelligence Review},
volume = {58},
pages = {323},
year = {2025}
}

@INCOLLECTION{Kovachki2024,
author = {Kovachki, N. B. and Lanthaler, S. and Stuart, A. M.},
title = {Chapter 9 -- Operator learning: Algorithms and analysis},
editor = {Siddhartha Mishra and Alex Townsend},
series = {Handbook of Numerical Analysis},
publisher = {Elsevier},
volume = {25},
pages = {419-467},
year = {2024},
booktitle = {Numerical Analysis Meets Machine Learning}
}

@ARTICLE{Tartakovsky2020,
author = {Tartakovsky, A. M. and Ortiz Marrero, C. and Perdikaris, P. and Tartakovsky, G. D. and Barajas‐Solano, D.},
title = {Physics-Informed Deep Neural Networks for Learning Parameters and Constitutive Relationships in Subsurface Flow Problems},
journal = {Water Resources Research},
number = {5},
volume = {56},
publisher = {American Geophysical Union (AGU)},
year = {2020},
month = {04}
}

@ARTICLE{Xu2022,
author = {Xu, K. and Darve, E.},
title = {Physics Constrained Learning for Data-Driven Inverse Modeling from Sparse Observations},
journal = {Journal of Computational Physics},
year = {2022},
volume = {453},
pages = {110938}
}

@ARTICLE{Sun2023,
author = {Sun, Y. and Sengupta, U. and Juniper, M.},
title = {Physics‐Informed Deep Learning for Simultaneous Surrogate Modelling and {PDE}‐Constrained Optimization of an Airfoil Geometry},
journal = {Computer Methods in Applied Mechanics and Engineering},
volume = {411},
pages = {116042},
year = {2023}
}

@ARTICLE{MowlaviNabi2023,
author = {Mowlavi, S. and Nabi, S.},
title = {Optimal Control of {PDE}s Using Physics-Informed Neural Networks},
journal = {Journal of Computational Physics},
volume = {473},
pages = {111731},
year = {2023}
}

@ARTICLE{Louboutin2023,
author = {Louboutin, M. and Yin, Z. and Orozco, R. L. and Grady, T. J., II and Siahkoohi, A. and Rizzuti, G. and Witte, P. A. and M{\o}yner, O. and Gorman, G. J. and Herrmann, F. J.},
title = {Learned Multiphysics Inversion with Differentiable Programming and Machine Learning},
journal = {The Leading Edge},
volume = {42},
number = {7},
pages = {474-486},
year = {2023},
month = {07}
}

@ARTICLE{Keil2022,
author = {Keil, T. and Kleikamp, H. and Lorentzen, R. J. and Oguntola, M. B. and Ohlberger, M.},
title = {Adaptive Machine Learning-based Surrogate Modeling to Accelerate {PDE}-Constrained Optimization in Enhanced Oil Recovery},
journal = {Advances in Computational Mathematics},
volume = {48},
pages = {73},
year = {2022}
}

@ARTICLE{Guth2024,
author = {Guth, P. A. and Schillings, C. and Weissmann, S.},
title = {One-Shot Learning of Surrogates in {PDE}-Constrained Optimization under Uncertainty},
journal = {SIAM/ASA Journal on Uncertainty Quantification},
volume = {12},
number = {2},
pages = {614-645},
year = {2024}
}

@ARTICLE{Homescu2002,
author = {Homescu, C. and Navon, I. M. and Li, Z.},
title = {Suppression of Vortex Shedding for Flow Around a Circular Cylinder Using Optimal Control},
journal = {Int. J. Numer. Meth. Fluids},
volume = {38},
pages = {43-69},
year = {2002}
}

@BOOK{Kirchhoff1894,
author = {Kirchhoff, G.},
title = {Vorlesungen \"uber die Theorie der Warme},
publisher = {Barth},
address = {Leipzig},
year = {1894}
}

@ARTICLE{Zhao2021,
author = {Zhao, H. and Braatz, R. D. and Bazant, M. Z.},
title = {Image Inversion and Uncertainty Quantification For Constitutive Laws of Pattern Formation},
journal = {Journal of Computational Physics},
volume = {436},
pages = {110279},
year = {2021}
}

@BOOK{Danaila2021,
author = {Danaila, I. and Kaplanski, F. and Sazhin, S. S.},
title = {Vortex Ring Models},
series = {Mathematical Engineering},
publisher = {Springer},
year = {2021}
}

@ARTICLE{Brunk2023,
author = {Brunk, A. and Egger, H. and Habrich, O.},
title = {On Uniqueness and Stable Estimation of Multiple Parameters in the {C}ahn--{H}illiard Equation},
journal = {Inverse Problems},
year = {2023},
month = {apr},
publisher = {IOP Publishing},
volume = {39},
number = {6},
pages = {065002}
}

@ARTICLE{Akerson2025,
author = {Akerson, A. and Rajan, A. and Bhattacharya, B.},
title = {Learning Constitutive Relations from Experiments: 1. {PDE} Constrained Optimization},
journal = {Journal of the Mechanics and Physics of Solids},
volume = {201},
pages = {106128},
year = {2025}
}

@BOOK{GriewankWalther2008,
author = {Griewank, A. and Walther, A.},
title = {Evaluating Derivatives: Principles and Techniques of Algorithmic Differentiation},
publisher = {SIAM},
year = {2008},
edition = {2}
}

@ARTICLE{GilesPierce2000,
author = {Giles, M. B. and Pierce, N. A.},
title = {An Introduction to the Adjoint Approach to Design},
journal = {Flow, Turbulence and Combustion},
year = {2000},
volume = {65},
number = {3},
pages = {393-415}
}

@BOOK{Naumann2011,
author = {Naumann, U.},
title = {The Art of Differentiating Computer Programs: An Introduction to Algorithmic Differentiation},
publisher = {SIAM},
year = {2011}
}

@ARTICLE{Fuhg2025,
author = {Fuhg, J. N. and Anantha Padmanabha, G. and Bouklas, N. and Bahmani, B. and Sun, W. and Vlassis, N. N. and Flaschel, M. and Carrara, P. and De Lorenzis, L.},
title = {A Review on Data-Driven Constitutive Laws for Solids},
journal = {Archives of Computational Methods in Engineering},
year = {2025},
volume = {32},
pages = {1841-1883}
}

@ARTICLE{Richardson2018,
author = {Richardson, G. and Foster, J. M. and Sethurajan, A. K. and Krachkovskiy, S. A. and Halalay, I. C. and Goward, G. R. and Protas, B.},
title = {The Effect of Ionic Aggregates on the Transport of Charged Species in {L}ithium Electrolyte Solutions},
journal = {Journal of The Electrochemical Society},
year = {2018},
month = {jun},
publisher = {The Electrochemical Society},
volume = {165},
number = {9},
pages = {H561}
}

@ARTICLE{Ahmadi2026,
author = {Ahmadi, A. and Sanders, K. J. and Goward, G. R. and Protas, B.},
title = {Data-Driven Approach to Learning Optimal Forms of Constitutive Relations in Models Describing {L}ithium Plating in Battery Cells},
journal = {Computers \& Chemical Engineering},
volume = {204},
pages = {109252},
year = {2026}
}

@ARTICLE{Daniels2023,
author = {Daniels, L. and Sahu, S. and Sanders, K. J. and Goward, G. R. and Foster, J. M. and Protas, B.},
title = {Learning Optimal Forms of Constitutive Relations Characterizing Ion Intercalation from Data in Mathematical Models of {L}ithium-Ion Batteries},
journal = {The Journal of Physical Chemistry C},
volume = {127},
number = {35},
pages = {17508-17523},
year = {2023}
}

@ARTICLE{Escalante2020,
author = {{Morales Escalante}, J. M. and Ko, W. and Foster, J. M. and Krachkovskiy, S. and Goward, G. and Protas, B.},
title = {Discerning Models of Phase Transformations in Porous Graphite Electrodes: Insights from Inverse Modelling Based on {MRI} Measurements},
journal = {Electrochimica Acta},
volume = {349},
pages = {136290},
year = {2020}
}

@ARTICLE{MatharuProtas2020,
author = {Matharu, P. and Protas, B.},
title = {Optimal Closures in a Simple Model for Turbulent Flows},
journal = {SIAM Journal on Scientific Computing},
volume = {42},
number = {1},
pages = {B250-B272},
year = {2020}
}

@ARTICLE{MatharuProtas2022,
title = {Optimal Eddy Viscosity in Closure Models for Two-Dimensional Turbulent Flows},
author = {Matharu, P. and Protas, B.},
journal = {Phys. Rev. Fluids},
volume = {7},
issue = {4},
pages = {044605},
numpages = {21},
year = {2022},
month = {Apr},
publisher = {American Physical Society}
}

@ARTICLE{MatharuProtas2024,
author = {Matharu, P. and Protas, B.},
title = {Adjoint-based Enforcement of State Constraints in {PDE} Optimization Problems},
journal = {Journal of Computational Physics},
volume = {517},
pages = {113298},
year = {2024}
}

@ARTICLE{Sethurajan2019,
author = {Sethurajan, A. K. and Krachkovskiy, S. and Goward, G. R. and Protas, B.},
title = {Bayesian Uncertainty Quantification in Inverse Modeling of Electrochemical Systems},
journal = {Journal of Computational Chemistry},
year = {2019},
volume = {40},
issue = {5},
pages = {740-752}
}

@ARTICLE{Sethurajan2015,
author={Sethurajan, A.K. and Krachkovskiy, S.A. and Halalay, I.C. and Goward, G.R. and Protas, B.},
title={Accurate Characterization of Ion Transport Properties in Binary Symmetric Electrolytes Using In Situ
{NMR} Imaging and Inverse Modeling},
journal={Journal of Physical Chemistry {B}},
year={2015},
volume={119},
number={37},
pages={12238-12248}
}

@ARTICLE{ProtasNoack2015,
author={Protas, B. and Noack, B.R. and \"{O}sth, J.},
title={Optimal Nonlinear Eddy Viscosity in {G}alerkin Models of Turbulent Flows},
journal={Journal of Fluid Mechanics},
year={2015},
volume={766},
pages={337-367}
}

@ARTICLE{DanailaProtas2015,
author={Danaila, I. and Protas, B.},
title={Optimal Reconstruction of Inviscid Vortices},
journal={Proceedings of the Royal Society A: Mathematical, Physical and Engineering Sciences},
year={2015},
volume={471},
number={2180},
art_number={20150323}
}

@ARTICLE{ProtasNoack2014,
author={Protas, B. and Noack, B.R. and Morzy\'{n}ski, M.},
title={An Optimal Model Identification for Oscillatory Dynamics with a Stable Limit Cycle},
journal={Journal of Nonlinear Science},
year={2014},
volume={24},
number={2},
pages={245-275}
}

@ARTICLE{ProtasBewleyHagen2004,
author={Protas, B. and Bewley, T.R. and Hagen, G.},
title={A Computational Framework for the Regularization of Adjoint Analysis in Multiscale {PDE} Systems},
journal={Journal of Computational Physics},
year={2004},
volume={195},
number={1},
pages={49-89}
}

@ARTICLE{VolkovProtasLiaoGlander2009,
author={Volkov, O. and Protas, B. and Liao, W. and Glander, D.},
title={Adjoint-based Optimization of Thermo-Fluid Phenomena in Welding Processes},
journal={Journal of Engineering Mathematics},
year={2009},
volume={65},
pages={201-220}
}

@ARTICLE{Protas2008,
author={Protas, B.},
title={Adjoint-Based Optimization of {PDE} Systems with Alternative Gradients},
journal={Journal of Computational Physics},
year={2008},
volume={227},
pages={6490-6510}
}

@BOOK{AdlerHolder2022,
author = {Adler, A. and Holder, D.},
title = {Electrical Impedance Tomography. Methods, History, and Applications},
edition = {2nd},
publisher = {CRC Press},
address = {Boca Raton, FL},
year = {2022}
}

@BOOK{Batchelor1988,
author = {Batchelor, G. K.},
title = {An Introduction to Fluid Dynamics},
edition = {7th},
publisher = {Cambridge University Press},
address = {Cambridge, UK},
year = {1988}
}

@BOOK{Saffman1992,
author = {Saffman, P. G.},
title = {Vortex Dynamics},
publisher = {Cambridge University Press},
address = {Cambridge, UK},
year = {1992}
}

@BOOK{Holmes2012,
author = {Holmes, P. and Lumley, J. L. and Berkooz, G. and Rowley, C. W.},
title = {Turbulence, Coherent Structures, Dynamical Systems and Symmetry},
series = {Cambridge Monographs on Mechanics},
edition = {2nd},
publisher={Cambridge University Press},
year={2012}
}

@BOOK{Noack2011,
author = {Noack, B. R. and Morzy{\'n}ski, M. and Tadmor, G.},
title = {Reduced-Order Modelling for Flow Control},
publisher = {Springer},
series = {CISM International Centre for Mechanical Sciences},
volume = {528},
address = {Vienna},
year = {2011}
}

\end{document}